\input ./style/arxiv-general.cfg
\documentclass[aop,MSNbibl,seceqn,dvips]{arximspdf}
\makeatletter
   \@ifpackageloaded{graphicx}{}{\usepackage{graphicx}}
\makeatother

\doi{10.1214/14-AOP960} 
\volume{43}
\issue{6}
\pubyear{2015}
\firstpage{3279}
\lastpage{3336}
\docsubty{FLA}

\makeatletter
\newproclaim{example}{Example}
\newcommand{\eqref}[1]{(\ref{#1})}
\newtheorem{thmm}{Theorem}
\newtheorem{coro}{Corollary}
\newtheorem{lemma}{Lemma}[section]
\newtheorem{prop}{Proposition}
\newcommand{\bd}{\mathbf}
\newcommand{\bdd}{\bolds}
\def\Tr{\operatorname{Tr}}
\makeatother

\begin{document}
\begin{frontmatter}

\title{Moments of traces of circular beta-ensembles}
\runtitle{Moments of traces of circular beta-ensembles}

\begin{aug}
\author[A]{\fnms{Tiefeng}~\snm{Jiang}\thanksref{T1}\ead[label=e1]{jiang040@umn.edu}}
\and
\author[B]{\fnms{Sho} \snm{Matsumoto}\corref{}\thanksref{T2}\ead[label=e2]{sho-matsumoto@math.nagoya-u.ac.jp}\ead
[label=e3]{shom@sci.kagoshima-u.ac.jp}}
%
\thankstext{T1}{Supported in part by NSF Grants DMS-04-49365,
DMS-12-08982 and DMS-14-06279.}
\thankstext{T2}{Supported in part by
JSPS Grant-in-Aid for Young Scientists (B) 25800062.}
\runauthor{T. Jiang and S. Matsumoto}
\affiliation{University of Minnesota and Nagoya University}
%
%
\address[A]{School of Statistics\\
University of Minnesota\\
224 Church Street SE\\
Minneapolis, Minnesota 55455\\
USA\\
\printead{e1}}%
%
\address[B]{Graduate School of Mathematics\\
Nagoya University\\
Furocho, Chikusaku, Nagoya\\
Japan\\
and\\
Graduate School of Science and Engineering\\
Kagoshima University\\
1-21-35, Korimoto, Kagoshima\\
Japan \\
\printead{e2}\\
\phantom{E-mail:\ }\printead*{e3}}
\end{aug}
%

\received{\smonth{3} \syear{2013}}
\revised{\smonth{8} \syear{2014}}

%
\begin{abstract}
Let $\theta_1, \ldots, \theta_n$ be random variables from Dyson's
circular $\beta$-ensemble with probability density function
$\operatorname{Const}\cdot\prod_{1\leq j< k\leq n}|e^{i\theta_j} - e^{i\theta
_k}|^{\beta}$. For each $n\geq2$ and $\beta>0$, we obtain some
inequalities on $\mathbb{E}  [ p_{\mu}(Z_n)\overline{p_{\nu
}(Z_n)}  ]$, where $Z_n=(e^{i\theta_1}, \ldots, e^{i\theta_n})$
and $p_{\mu}$ is the power-sum symmetric function for partition $\mu
$. When $\beta=2$, our inequalities recover an identity by Diaconis
and Evans for Haar-invariant unitary matrices. Further, we have
the following:
$\lim_{n\to\infty}\mathbb{E} [p_{\mu}(Z_n)\overline{p_{\nu
}(Z_n)}  ] = \delta_{\mu\nu} (\frac{2}{\beta}
)^{l(\mu)}z_{\mu}$
for any $\beta>0$ and partitions $\mu,\nu$;
$\lim_{m\to\infty}\mathbb{E} [ |p_m(Z_n)|^2  ] = n$ for
any $\beta>0$ and $n\geq2$,
where $l(\mu)$ is the length of $\mu$ and $z_{\mu}$ is explicit on
$\mu$.
These results apply to the three important ensembles: COE ($\beta=1$),
CUE ($\beta=2$) and CSE ($\beta=4$). We further examine the
nonasymptotic behavior of $\mathbb{E} [ |p_m(Z_n)|^2  ]$
for $\beta=1, 4$. The central limit theorems of $\sum_{j=1}^ng(e^{i\theta_j})$ are obtained when (i) $g(z)$ is a polynomial
and $\beta>0$ is arbitrary, or (ii) $g(z)$ has a Fourier expansion and
$\beta=1,4$. The main tool is the Jack function.
\end{abstract}



%
\end{frontmatter}

\section{Introduction}\label{intro}
Let $M_n$ be an $n\times n$ Haar-invariant unitary matrix, that is, the
entries of unitary matrix $M_n$ are random variables satisfying that
the probability distribution of the entries of $M_n$ is the same as
that of $UM_n$ and that of $M_nU$ for any $n\times n$ unitary matrix
$U$. Diaconis and Evans (Theorem~2.1 from \cite{Evans}) proved that

(a) \textit{Consider $a=(a_1, \ldots, a_k)$ and $b=(b_1, \ldots, b_k)$
with $a_j, b_j\in\{0,1,2,\ldots\}$. Then for $n\geq\sum_{j=1}^{k}ja_j \vee\sum_{j=1}^kjb_j$},
%
\begin{equation}
\label{DiaEvans} \mathbb{E} \Biggl[\prod_{j=1}^k
\bigl(\Tr\bigl(M_n^j\bigr)\bigr)^{a_j}\overline{
\bigl(\Tr \bigl(M_n^j\bigr)\bigr)^{b_j}} \Biggr]
=\delta_{ab}\prod_{j=1}^kj^{a_j}a_j!,
\end{equation}
\textit{where $\delta_{ab}$ is Kronecker's delta.}

\textup{(b)} \textit{For any positive integers $j$ and $k$},
%
\begin{equation}
\label{weapon} \mathbb{E} \bigl[ \Tr\bigl(M_n^j\bigr)
\overline{\Tr\bigl(M_n^k\bigr)} \bigr]=\delta
_{jk} \cdot j \wedge n.
\end{equation}

The idea of the proof is based on the group representation theory of
unitary group $U(n)$. Some other derivations for (\ref{DiaEvans}) and
(\ref{weapon}) are given in \cite{Shahshahani,Pastur,Rains,Stolz}.
The right-hand side\vspace*{1pt} of (\ref{DiaEvans}) is evidently equal to $\mathbb
{E} [\prod_{j=1}^k\xi_j^{a_j}\bar{\xi}_j^{b_j} ]$ where
$\xi_j$'s are independent complex-normal random variables with $\xi
_j\sim\mathbb{C}N(0, j)$ for each $j$.

\begin{figure}

\includegraphics{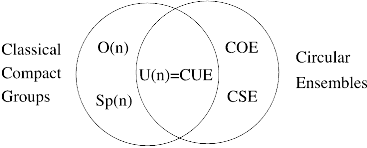}

\caption{Circular Ensembles and Haar-invariant matrices from classical
compact groups.}\label{figure-ensemble}
\end{figure}

Notice an $n\times n$ Haar-invariant unitary matrix is also called a
CUE, which belongs to the Circular Ensembles of three members: the
Circular Orthogonal Ensemble (COE), the Circular Unitary Ensemble (CUE)
and the Circular Symplectic Ensemble (CSE); see Figure~\ref{figure-ensemble}
for the relationship, where the left circle consists
of matrices which induce the Haar probability measure on the orthogonal
group~$O(n)$, Haar probability measure on the unitary group $U(n)$ and
Haar probability measure on the real symplectic group $Sp(n)$, respectively.

Let $e^{i\theta_1}, \ldots, e^{i\theta_n}$ be the eigenvalues of an
$n\times n$ Haar-invariant unitary matrix, or equivalently, an $n\times
n$ CUE, it is known (see, e.g., \cite{Forrester,Mehta}) that the
density function of $\theta_1, \ldots,\theta_n$ is $f(\theta_1,
\ldots, \theta_n|\beta)$ with $\beta=2$, where
%
\begin{equation}
\label{bensemble} f(\theta_1, \ldots, \theta_n|\beta)=(2
\pi)^{-n}\cdot\frac{\Gamma
(1+\beta/2)^n}{\Gamma(1+\beta n/2)}\prod_{1\leq j<k \leq
n}\bigl|e^{i\theta_j}-e^{i\theta_k}\bigr|^{\beta}
\end{equation}
with $\beta>0$ and $\theta_i\in[0, 2\pi)$ for $1\leq i \leq n$. The
density function of $\theta_1, \ldots,\theta_n$ for the COE is
$f(\theta_1, \ldots, \theta_n|\beta)$ with $\beta=1$, and that for
the CSE is $f(\theta_1, \ldots, \theta_n|\beta)$ with $\beta=4$.

The purpose of this paper is to study the analogues of (\ref
{DiaEvans}) and (\ref{weapon}) for the circular $\beta$-ensembles
with density function $f(\theta_1, \ldots, \theta_n|\beta)$ in
(\ref{bensemble}) for any $\beta>0$. Further, we develop the central
limit theorems for functions of $(e^{i\theta_1}, \ldots, e^{i\theta
_n})$. Before stating the main results, we next introduce some
background about the circular $\beta$-ensembles.

The circular ensembles were first introduced by physicist Dyson \cite{Dyson62a,Dyson62b,Dyson62c} for the study of nuclear
scattering data. In fact, as studied in \cite{Dyson62a}, Dyson shows
that the consideration of time reversal symmetry leading to the three
Gaussian ensembles behaves equally well to unitary matrices. A time
reversal symmetry requires that $U = U^T$, no time reversal symmetry
has no constraint, and a time reversal symmetry for a system with an
odd number of spin $1/2$ particles requires $U = U^D$, where $D$ denotes
the quaternion dual. Choosing such matrices with a uniform probability
then gives COE, CUE and CSE, respectively (see, e.g., \cite{Forresterbook,Mehta}). The entries of COE and CUE are asymptotically
complex normal random variables when the sizes of the matrices are
large \cite{Jiang09a,Jiang09c,Jiang06}.

Let $U$ be an $n\times n$ Haar-invariant unitary matrix. As mentioned
earlier, $U$ is also a CUE; the matrix $U^TU$ gives a COE.
Furthermore, the matrix $U^D U$ gives a CSE when $n$ is even; see
Chapter~9 from \cite{Mehta}. For the relations among the zonal
polynomials, the Schur functions, the Gelfand pairs and the three
circular ensembles; see, for example, Chapter VII in \cite{Macdonald1998} or Section~2.7 in \cite{Blower} for reference.

Now we consider the moments in (\ref{DiaEvans}) and (\ref{weapon})
for the circular $\beta$-ensembles. Taking $\beta=1$ in (\ref
{bensemble}), that is, choosing $W_n$ such that it is an $n\times n$
COE, by an elementary check in Lemma~\ref{blue}, we have
%
\begin{equation}
\label{glue} \mathbb{E} \bigl[\bigl |\Tr(W_n)\bigr|^2 \bigr]=
\frac{2n}{n+1}
\end{equation}
for all $n\geq2$. This suggests that, unlike the right-hand sides of
(\ref{DiaEvans}) or (\ref{weapon}) that are free of $n$, the moments
for the general circular $\beta$-ensemble may depend on $n$ for $\beta
\ne2$. In fact, by using the Jack functions, we will soon see from
(\ref{hehe}) below that the second moment in (\ref{glue}) does depend
on $n$ except $\beta=2$, in which case $W_n$ is an $n\times n$ CUE.

In this paper, we will first prove some inequalities on the moments in
(\ref{DiaEvans}) and (\ref{weapon}) for the circular $\beta
$-ensembles with arbitrary $\beta>0$. In particular, some of our
inequalities for $\beta=2$ recover the equality in (\ref{DiaEvans})
by Diaconis and Evans~\cite{Evans}. Further, we evaluate the limiting
behavior by letting $n\to\infty$ for the left-hand side in (\ref
{DiaEvans}) and $k\to\infty$ for the left-hand side in (\ref
{weapon}), respectively. Their limits exist and look quite similar to
the right-hand sides of (\ref{DiaEvans}) and (\ref{weapon}). Finally,
we spend much effort to study the central limit theorems of $\sum_{j=1}^ng(e^{i\theta_j})$ for two situations: (a) $g(x)$ is a
polynomial and $\beta>0$ is arbitrary; (b) $g(x)$ has a Fourier
expansion and $\beta=1,4$. The key to obtain (b) is the nonasymptotic
behavior of $\mathbb{E}|\sum_{j=1}^ne^{im\theta_j}|^2$ for any $n$
and $m$, which are analyzed in detail.

The method of the proof is the Jack functions. The main results are
obtained by using their orthogonal properties and combinatorial structures.

From the studies in this paper, it is obvious to see the importance of
understanding the circular $\beta$-ensembles through the Jack
functions. Realizing that the Jack functions are a special class of the
Macdonald polynomials, we have obtained the analogue of the results in
this paper in the setting of the Macdonald polynomials. These will be
published elsewhere in the future.

The organization of the rest of the paper is as follows. We present the
moment inequalities in Section~\ref{moment_inequality} and their
proofs are given in Section~\ref{Proof_moment_inequality}; the
nonasymptotic behavior of $\mathbb{E}|\sum_{j=1}^ne^{im\theta_j}|^2$
and the central limit theorems are stated in Section~\ref{CLT_introduction} and their
proofs are arranged in Section~\ref{Proof_CLT}. In the \hyperref[noodle_eat]{Appendix}, we prove (\ref{glue}) by
two ways different from the method of the Jack functions. Some other
explicit formulas of moments are also given in the same section.

\section{Moment inequalities for circular beta-ensembles}\label
{moment_inequality}

Let $\lambda=(\lambda_1, \lambda_2, \ldots)$ be a partition, that
is, the sequence is in nonincreasing order and only finite of $\lambda
_i$'s are nonzero. The weight of $\lambda$ is $|\lambda|=\lambda
_1+\lambda_2 + \cdots$. Denote by $m_i(\lambda)$ the multiplicity of
$i$ in $(\lambda_1, \lambda_2, \ldots)$ for each $i$, and $l(\lambda
)$ the length of $\lambda$: $l(\lambda)=m_1(\lambda)+ m_2(\lambda)
+ \cdots.$ Recall the convention $0!=1$. Set
%
\begin{equation}
\label{insect} z_{\lambda}=\prod_{i\geq1}i^{m_i(\lambda)}m_i(
\lambda)!.
\end{equation}
Let $\rho=(\rho_1, \rho_2, \ldots)$ be a partition, and
%
\begin{equation}
\label{remote} p_{\rho}=\prod_{i=1}^{l(\rho)}p_{\rho_i}\qquad
\mbox{where } p_k(x_1, x_2,
\ldots)=x_1^k + x_2^k + \cdots
\end{equation}
for integer $k\geq1$ and indeterminates $x_i$'s. The function $p_{\rho
}$ is called the power-sum symmetric function. For real number $\alpha
>0$, integers $K\geq1$ and $n\geq1$, define
two constants $A=A(n,K,\alpha)$ and $B=B(n,K,\alpha)$ by
%
\begin{eqnarray}
\label{meso} A &=& \biggl(1-\frac{|\alpha-1|}{n-K+\alpha} \delta(\alpha\geq1)
\biggr)^K \quad\mbox{and}
\nonumber
\\[-8pt]
\\[-8pt]
\nonumber
B &=& \biggl(1+\frac{|\alpha-1|}{n-K+\alpha} \delta(\alpha< 1) \biggr)^K,
\end{eqnarray}
%
where $\delta(\alpha\geq1)=1-\delta(\alpha<1)$ is $1$ if $\alpha
\geq1$, or $0$ otherwise.
With these notation, we have one of main results as follows.

\begin{thmm}\label{face} Let $\beta>0$ and $\theta_1, \ldots,
\theta_n$ have density $f(\theta_1, \ldots, \theta_n|\beta)$ as in~(\ref{bensemble}). Set $Z_n=(e^{i\theta_1}, \ldots, e^{i\theta_n})$
and $\alpha=2/\beta$. For partitions $\mu$ and $\nu$, the following hold:
\begin{longlist}[(a)]
\item[(a)] If $n\geq K=|\mu|$, then
\[
A\leq \frac{\mathbb{E}  [ |p_{\mu}(Z_n)|^2  ]}{\alpha
^{l(\mu)}z_{\mu}}\leq B.
\]
\item[(b)] If $|\mu|\ne|\nu|$, then $\mathbb{E}  [p_{\mu
}(Z_n)\overline{p_{\nu}(Z_n)}  ]=0$.
If $\mu\ne\nu$ and $n\geq K=|\mu| \vee|\nu|$, then
\[
\bigl|\mathbb{E} \bigl[ p_{\mu}(Z_n)\overline{p_{\nu}(Z_n)}
\bigr] \bigr| \leq \max\bigl\{|A-1|,|B-1|\bigr\}\cdot\alpha^{(l(\mu)+l(\nu))/2}(z_{\mu
}z_{\nu})^{1/2}.
\]
\item[(c)] There exists a constant $C$ depending only on $\beta$ such
that for any $m\geq1$ and $n\geq2$, we have
\[
\bigl| \mathbb{E} \bigl[ \bigl|p_m(Z_n)\bigr|^2 \bigr]-n\bigr |
\leq C\frac
{n^32^{n\beta}}{m^{1\wedge\beta}}.
\]
\end{longlist}
\end{thmm}

Take $\beta=2$ in (a) and (b) of Theorem~\ref{face}, then $A=1$ and
$B=1$. The two results recover the result of Diaconis and Evans in
(\ref{DiaEvans}). Further, letting $n\to\infty$ in (a) and (b) of
Theorem~\ref{face}, we see that $A$ and $B$ (depending on $n$)
converge to $1$; letting $m\to\infty$ in (c) of the theorem, then the
last term in (c) goes to $0$. So we obviously have the following results.

\begin{coro}\label{main} Let the conditions be as in Theorem~\ref
{face}. Then, for any $\beta>0$,
\begin{eqnarray*}
&& \mathrm{(a)} \quad \lim_{n\to\infty}\mathbb{E} \bigl[p_{\mu}(Z_n)
\overline {p_{\nu}(Z_n)} \bigr] = \delta_{\mu\nu}
\biggl(\frac{2}{\beta
} \biggr)^{l(\mu)}z_{\mu};
\\
&& \mathrm{(b)} \quad \lim_{m\to\infty}\mathbb{E} \bigl[ \bigl|p_m(Z_n)\bigr|^2
\bigr] = n \qquad\mbox{for any $n\geq2$}.
\end{eqnarray*}
\end{coro}

Part (b) of the above corollary says that, as $m\to\infty$, the limit
of $\mathbb{E} \times\break  [|p_m(Z_n)|^2  ] $ does not depend on
parameter $\beta$, which is consistent with (\ref{weapon}). We
further take a careful examination on $\mathbb{E} [ |p_m(Z_n)|^2
 ]$ as $\beta=1$ and $4$. Some upper bounds of $\mathbb{E} [
|p_m(Z_n)|^2  ]$ are given in Propositions \ref{mouth} and \ref
{search}. By studying $A$ and $B$ in (\ref{meso}), we have the
following corollary from Theorem~\ref{face}.

\begin{coro}\label{theory} Let $\beta>0$ and $f(\theta_1, \ldots,
\theta_n|\beta)$ be as in (\ref{bensemble}). Set $\alpha=2/\beta$
and $Z_n=(e^{i\theta_1}, \ldots, e^{i\theta_n})$. Let $\mu$ and
$\nu$ be partitions with $\mu\ne\nu$ and $K= |\mu|\vee|\nu|$. If
$n\geq2K$, then
\begin{eqnarray*}
&& \mathrm{(a)}\quad  \biggl| \frac{\mathbb{E} [|p_{\mu}(Z_n)|^2 ]}{\alpha^{l(\mu
)}z_{\mu}}-1 \biggr|\leq\frac{6|1-\alpha|K}{n};
\\
&& \mathrm{(b)}\quad  \bigl|\mathbb{E} \bigl[p_{\mu}(Z_n)\overline{p_{\nu
}(Z_n)}
\bigr] \bigr| \leq\frac{6|1-\alpha|K}{n}\cdot\alpha ^{(l(\mu)+l(\nu))/2}(z_{\mu}z_{\nu})^{1/2}.
\end{eqnarray*}
\end{coro}

The above results are in the forms of inequalities or limits. We
actually derive an exact formula
in Proposition~\ref{exact}
to compute $\mathbb{E}  [ |p_{\mu}(Z_n)|^2  ]$ for every
partition $\mu$. In general, it is not easy to evaluate this quantity
for arbitrary $\mu$, however, we are able to do so when $\mu$ is
special. For instance, by using the exact formula we calculate the moment
in (\ref{glue}) for any $\beta>0$ as follows.

\begin{example*}For any $n\geq1$,
%
\begin{equation}
\label{hehe} \mathbb{E} \bigl[\bigl|p_1(Z_n)\bigr|^2
\bigr] =\frac{2}{\beta}\frac{n}{n-1+2\beta^{-1}}= \cases{ \displaystyle\frac{2n}{n+1},&\quad $\mbox{if
$\beta=1$;}$\vspace*{2pt}
\cr
1, &\quad $\mbox{if $\beta=2$}$;\label{cry}
\vspace*{2pt}
\cr
\displaystyle\frac{n}{2n-1}, &\quad $\mbox{if $\beta=4$.}$}
\end{equation}
The verification of this formula through Proposition~\ref{exact}
is provided in the \hyperref[noodle_eat]{Appendix}. We also give $\mathbb
{E}[|p_1(Z_n)|^4]$, $\mathbb{E}[|p_2(Z_n)|^2]$ and $\mathbb
{E}[p_2(Z_n)\overline{p_1(Z_n)^2}]$ in closed forms in the \hyperref[noodle_eat]{Appendix}.
\end{example*}

The main tool used in our proofs is the Jack functions. Diaconis and
Evans~\cite{Evans} and Diaconis and Shahshahani \cite{Shahshahani}
use the group representation theory to study (\ref{DiaEvans}) and
(\ref{weapon}) because $U(n)$ is a compact Lie group. The situations
for the Circular Orthogonal Ensembles ($\beta=1$) and the Circular
Symplectic Ensembles ($\beta=4$) are different. In fact, the two
ensembles are not groups.

The proofs of (\ref{DiaEvans}) and (\ref{weapon}) involve with the
Schur functions. The connection is that the irreducible characters of
the unitary groups, when seen as symmetric functions in the
eigenvalues, are given by Schur functions. Looking at Figure~\ref{figure-ensemble}, an Haar-invariant unitary matrix is also a CUE.
From the perspective of symmetric functions, the COE is connected to
the zonal polynomials, and the CSE to symplectic zonal polynomials.
The three functions are special cases of the Jack polynomial
$J_{\lambda}^{(\alpha)}$ with $\alpha=1, 2$ and $1/2$, respectively,
where $\lambda$ is a partition. See Section~\ref{proofs} for this or
\cite{Macdonald1998} for general properties of the Jack polynomials.
By using the Jack functions, we are able to prove (a) and (b) in
Theorem~\ref{face}. Part (c) in the theorem is proved by evaluating
the expectation/integral with respect to $f(\theta_1, \ldots, \theta
_n|\beta)$ in~(\ref{bensemble}) directly.

Treating $n$ as a variable, the bound $n^32^{n\beta}m^{-(1\wedge\beta
)}$ in (c) of Theorem~\ref{face} seems quite large. It is possibly to
be improved. However, as $\beta=4$, we show in Proposition~\ref
{search} in the next section that $\mathbb{E}  [ |p_{m}(Z_n)|^2
 ]$ has the scale of $m\log m$ when $n$ and $m$ are not far from
each other. This partially explains why the bound is large.

\section{Central limit theorems for circular beta-ensembles}\label
{CLT_introduction}
For the sake of precision, we replace $Z_n$ appeared earlier with
$Z_n^{\alpha}$. Specifically, let $Z_n^{\alpha}=(e^{i\theta_1},
\ldots, e^{i\theta_n})$ follow the $\beta$-circular ensemble with
$\alpha=\frac{2}{\beta}$ and the density function $f(\theta_1,
\ldots, \theta_n|\beta)$ as in (\ref{bensemble}). According to our
notation in previous sections, $p_m(Z_n^\alpha)=\sum_{j=1}^ne^{im\theta_j}$ for any integer $m\geq0$. In the paper, the
symbol $\mathbb{C}N(0, \sigma^2)$ stands for the complex normal
distribution generated by $\sigma\cdot(\xi_1+i\xi_2)/\sqrt{2}$,
where $\xi_1$ an $\xi_2$ are i.i.d. real random variables with the
standard normal distribution $N(0,1)$. The first result is a CLT for
general circular $\beta$-ensemble.

\begin{thmm}[(CLT for any $\beta$-circular ensemble)]\label{Nevada}
Let $Z_n^{\alpha}=(e^{i\theta_1}, \ldots, e^{i\theta_n})$ follow
the $\beta$-circular ensemble. Then, for fixed $m\geq1$, the random
vector $(p_1(Z_n^{\alpha}), p_2(Z_n^{\alpha}), \ldots,
p_m(Z_n^{\alpha}))$ converges weakly to $(\xi_1, \ldots, \xi_m)$ as
$n\to\infty$, where $\xi_j$'s are independent random variables with
$\xi_j \sim\mathbb{C}N(0, \frac{2j}{\beta})$ for each $j$.
\end{thmm}

An immediate consequence of the theorem is as follows.

\begin{coro}\label{campus} Let $(e^{i\theta_1}, \ldots, e^{i\theta
_n})$ follow the $\beta$-circular ensemble. Let $g(z)=\sum_{k=0}^mc_kz^k$ with fixed $m\geq1$ and $c_k\in\mathbb{C}$ for all
$k$. Set $X_n=\sum_{j=1}^ng(e^{i\theta_j})$. Then $X_n-\mu_n$
converges weakly to $\mathbb{C}N(0,\sigma^2)$ as $n\to\infty$, where
\[
\mu_n=n c_0 \quad\mbox{and}\quad \sigma^2=
\frac{2}{\beta}\sum_{k=1}^mk|c_k|^2.
\]
\end{coro}

We next study the central limit theorem when the function $g(z)$ is not
a polynomial. To avoid a lengthier paper, we only focus on the cases
$\beta=1$ and $\beta=4$. A~discussion on the general case will be
given later in this section. We first need to understand the variance
of $p_m(Z_n^\alpha)$.

\begin{prop}[(Bound of variance on COE)]\label{mouth} For all $m\geq1$,
$n \geq2$ and \mbox{$\beta=1$}, there exists a universal constant $K>0$ such that
\begin{eqnarray*}
\mathbb{E} \bigl[ \bigl|p_{m}\bigl(Z_n^2
\bigr)\bigr|^2 \bigr] \leq \cases{ 2m,&\quad $\mbox{if $1\leq m\leq n$}$;
\vspace*{2pt}
\cr
Kn,&\quad  $\mbox{if $m> n$}$.}
\end{eqnarray*}
\end{prop}

\begin{prop}[(Bound of variance on CSE)]\label{search}  Let $\beta=4$.
Then there exists a universal constant $K>0$ such that the following
hold:
\begin{longlist}[(iii)]
\item[(i)] $\mathbb{E} [|p_m(Z_n^{1/2})|^2] \leq K \delta^{-1}n$ for all
$m\geq(1+\delta)n$ and $\delta\in(0,1]$.

\item[(ii)] $\mathbb{E} [|p_m(Z_n^{1/2})|^2] \leq K m\log(m+1)$ for all
$m\geq1$ and $n\geq2$.

\item[(iii)] $\mathbb{E} [|p_m(Z_n^{1/2})|^2] \geq K(w+1)^{-2} m\log m$ for
all $12\leq n\leq m \leq2n$ where $w=m-n\geq0$.
\end{longlist}
\end{prop}

From (ii) and (iii), we see that $\mathbb{E} [|p_m(Z_n^{1/2})|^2]$ is
of the scale ``$m\log m$'' when $m$ and $n$ are not far from each other.
It is known from (\ref{weapon}) and Proposition~\ref{mouth} that
$\mathbb{E}|p_m(Z_n^\alpha)|^2 \leq Kn$ for any $n\geq2$, $m\geq1$
and $\beta=1, 2$, where $K$ is a universal constant. This together
with (b) of Corollary~\ref{main} seems to suggest that the second
moment for $\beta=4$ is always bounded by $Kn$. Proposition~\ref
{search} tells us a different story. However, (b) of Corollary~\ref
{main} is indeed consistent with (i).

The proofs of Propositions \ref{mouth} and \ref{search} are very
involved. We use the combinatorial structure (\ref{please}) to
understand the second moments. Major effort is devoted to analyzing
(\ref{please}) through (\ref{earth}) and (\ref{lantern}).

Another way to calculate above variance is through the covariance of
$e^{im\theta_1}$ and $e^{im\theta_2}$ by symmetry [see (\ref
{her_him})], which again can be computed by using the two-point
correlation function $\rho_{(2)}(\theta_1, \theta_2)$. The explicit
form of $\rho_{(2)}(\theta_1, \theta_2)$ is given in Proposition~13.2.2 from \cite{Forresterbook}. It seems very hard to estimate the
variance by using the proposition. But it is possible in principle.

\begin{thmm}[(CLT for COE)]\label{circular_case} Let $(e^{i\theta_1},
\ldots, e^{i\theta_n})$ follow the circular orthogonal ensemble
($\beta=1$).
Let $\{a_j, b_j \in\mathbb{C};   j=1,2,\ldots\}$
satisfy $\sum_{j=1}^{\infty}j(|a_j|^2+|b_j|^2)=\sigma^2\in(0,
\infty)$.
Then, $\sum_{j=1}^{\infty} (a_{j}p_j(Z_n^2)+ b_{j}\overline
{p_j(Z_n^2)}  )$ converges weakly to the law of $U+iV$ as $n\to
\infty$, where $(U, V)\in\mathbb{R}^2$ has the law $ N_2(\bd{0},
\bdd{\Sigma})$ with
\[
\bdd{\Sigma}= \pmatrix{ \displaystyle\sum_{j=1}^{\infty}j|a_j+
\bar{b}_j|^2 &\displaystyle 2\cdot\operatorname{Im}\Biggl(\sum
_{j=1}^{\infty}ja_jb_j\Biggr)
\vspace*{2pt}\cr
\displaystyle 2\cdot\operatorname{Im}\Biggl(\sum_{j=1}^{\infty}ja_jb_j
\Biggr) & \displaystyle\sum_{j=1}^{\infty
}j|a_j-
\bar{b}_j|^2 } .
\]
%
\end{thmm}

Obviously, if $b_j=0$ for all $j$, then $\bdd{\Sigma}=\sigma^2\bd
{I}_2$ with $\sigma^2=\sum_{j=1}^{\infty}j|a_j|^2$, and hence
$U+iV\sim\mathbb{C}N(0, 2\sigma^2)$.

\begin{thmm}[(CLT for CSE)]\label{symplectic_case} Let $(e^{i\theta
_1}, \ldots, e^{i\theta_n})$ follow the circular symplectic ensemble
($\beta=4$).
Let $\{a_j, b_j \in\mathbb{C};   j=1,2,\ldots\}$
satisfy $\sum_{j=1}^{\infty}(j\log j)\times\break (|a_j|^2+|b_j|^2)\in(0, \infty
)$. Set $\sigma^2=\sum_{j=1}^{\infty}j(|a_j|^2+|b_j|^2)$.
Then $\sum_{j=1}^{\infty} (a_{j}p_j(Z_n^{1/2})+ b_{j}\overline
{p_j(Z_n^{1/2})}  )$ converges weakly to the law of $U+iV$ as
$n\to\infty$, where $(U, V)\in\mathbb{R}^2$ has the law $ N_2(\bd
{0}, \bdd{\Sigma})$ with
\begin{eqnarray*}
\bdd{\Sigma}=\frac{1}{4} \pmatrix{\displaystyle \sum_{j=1}^{\infty}j|a_j+
\bar{b}_j|^2 & \displaystyle2\cdot\operatorname{Im}\Biggl(\sum
_{j=1}^{\infty}ja_jb_j\Biggr)
\vspace*{2pt}\cr
\displaystyle 2\cdot\operatorname{Im}\Biggl(\sum_{j=1}^{\infty}ja_jb_j
\Biggr) & \displaystyle\sum_{j=1}^{\infty
}j|a_j-
\bar{b}_j|^2 } .
\end{eqnarray*}
\end{thmm}

Similar to the comment below Theorem~\ref{circular_case}, if $b_j=0$
for all $j$, then $\bdd{\Sigma}=\sigma^2\bd{I}_2$ with $\sigma
^2=\frac{1}{4}\sum_{j=1}^{\infty}j|a_j|^2$, and hence $U+iV\sim
\mathbb{C}N(0, 2\sigma^2)$.

Though Proposition~\ref{search} says that $\mathbb{E}
[|p_m(Z_n^{1/2})|^2]$ is of scale ``$m\log m$'' when $m$ and $n$ are not
far from each other, the variance of the limiting distribution in
Theorem~\ref{symplectic_case} is not affected by this fact. The
variance is similar to those in the circular orthogonal and unitary
ensemble ($\beta=1, 4$).

Diaconis and Evans \cite{Evans} obtains the CLTs for the orthogonal
groups, the unitary groups and the symplectic groups. Their tool is the
identities in (\ref{DiaEvans}) and (\ref{weapon}).
Reviewing Corollary~\ref{theory}, we no longer have identities for any
$\beta\ne2$; this increases much difficulty to get the corresponding
CLTs. It is understandable because after all the three members in the
classical compact groups have group structures in addition to their
combinatorial ones. So the group representation theory can be possibly
used in the paper by Diaconis and Evans. The general circular $\beta
$-ensemble loses the former property and has only the combinatorial structure.

Johansson in \cite{Johansson} further explores the convergence speed
of $\Tr(M_n^m)$ to a normal distribution, where $m$ is fixed and $M_n$
is an Haar-invariant orthogonal, unitary or symplectic random matrix.
He shows that the convergence rate is exponentially fast.

By Proposition~\ref{search}, the condition ``$\sum_{j=1}^{\infty
}(j\log j)(|a_j|^2+|b_j|^2)<\infty$'' in Theorem~\ref{symplectic_case}
can be slightly relaxed. For simplicity, we just leave it as it is.
Also, the conditions ``$\sum_{j=1}^{\infty}j(|a_j|^2+|b_j|^2)< \infty
$'' and ``$\sum_{j=1}^{\infty}(j\log j)(|a_j|^2+|b_j|^2)< \infty$''
can be easily satisfied. For instance, the first condition is satisfied
if $a_j$ and $b_j$ are of the order $\frac{1}{j(\log j)^{(1/2)+\delta
}}$ for some $\delta>0$, and the second one is satisfied if $a_j$ and
$b_j$ are of the order $\frac{1}{j(\log j)^{1+\delta}}$ for some
$\delta>0$.

To study the number of eigenvalues falling in an arc of the unit circle
in the complex plane, namely, $\sum_{j=1}^n I(e^{i\theta_j} \in A)$
with $A$ being a subset of $S^1=\{z\in\mathbb{C};  |z|=1\}$, one
needs to handle the Fourier expansion of the indicator function
$I_{[a,b]}(x)$ with $[a, b]\subset[0, 2\pi]$. It is known from \cite{Evans} that the coefficients $a_j$ and $b_j$ in the contexts of
Theorems \ref{circular_case} and \ref{symplectic_case} are of scale
$\frac{1}{j}$. Our theorems do not cover this special case. By using a
construction of the circular $\beta$-ensemble, Killip \cite{Killip}
specifically considers this situation and obtains a CLT. The author
does not investigate the general CLTs as treated in our Theorems \ref
{Nevada}, \ref{circular_case} or \ref{symplectic_case}.

Finally, we provide some examples which satisfy the condition
\[
\sum_{j=1}^{\infty}(j\log j)
\bigl(|a_j|^2+|b_j|^2\bigr)<
\infty.
\]
They are the solutions of some classical partial differential
equations. We leave readers for the trivial calculations of the means
and the variances of the limiting normal distributions.

\begin{example*} Let $u=u(x,y)$ be defined on $\mathbb
{R}^2$ and satisfy the Laplace equation
\[
\cases{ \Delta u=0, &\quad $\mbox{$x^2+y^2 <
a^2$}$;\vspace*{2pt}
\cr
u=h(\theta), &\quad $\mbox{$x^2+y^2
= a^2$},$}
\]
where $h(\theta)$ is a known function and $a>0$ is given.
Let $(x,y)=(r\cos\theta, r\sin\theta)$.
The solution has a Poisson's formula. It can also be expressed in the
following Fourier series:
%
\begin{equation}
\label{little_light} u(r,\theta)=\frac{1}{2}A_0 + \sum
_{j=1}^{\infty}r^j(A_j \cos j
\theta+ B_j \sin j\theta)
\end{equation}
for $r\in(0, a)$ and $\theta\in[0, 2\pi]$,
where $A_j$'s and $B_j$'s are obtained from the Fourier series of
$h(\theta)$ so that
\[
A_j=\frac{1}{\pi a^j}\int_0^{2\pi}h(
\phi)\cos j\phi \,d\phi\quad \mbox{and}\quad B_j=\frac{1}{\pi a^j}\int
_0^{2\pi}h(\phi)\sin j\phi \,d\phi. 
\]
See, for example, more details on page 160 from \cite{Strauss}.
Clearly, if $C:= \break \sup_{\phi\in[0, 2\pi]}|h(\phi)|<\infty$, then
$|A_j|\leq\frac{2C}{a^j}$ and $|B_j|\leq\frac{2C}{a^j}$. And the
coefficients $|r^jA_j|$ and $|r^jB_j|$ in (\ref{little_light}) are
bounded by $ 2C(\frac{r}{a})^j$ for $0<r<a$. Then use the formulas
$\cos j\theta=\frac{e^{ij\theta}+ e^{-ij\theta}}{2}$ and $\sin
j\theta=\frac{e^{ij\theta}- e^{-ij\theta}}{2i}$ to transfer $u(r,
\theta)$ in~(\ref{little_light}) to the form of
$a_0+\sum_{j=1}^{\infty}(a_je^{ij\theta} + b_je^{-ij\theta})$,
where $a_j$'s and $b_j$'s are complex numbers. Fix $r<a$. It is easy to
see that $|a_j|=O((ra^{-1})^j)$ and $b_j=O((ra^{-1})^j)$ as $j\to
\infty$. Theorems \ref{circular_case} and \ref{symplectic_case} can
then be applied to get the CLT for $a_0+\sum_{j=1}^{\infty
}(a_jp_j(Z_n^{\alpha}) + b_j\overline{p_j(Z_n^{\alpha})} )$ for
$\alpha=2$ and $\alpha=\frac{1}{2}$, respectively.
\end{example*}

\begin{example*}Let $u(x, t)$ be a function defined on
$[0, \pi]\times[0, \infty)$. Consider the following heat equation
with boundary conditions defined by
%
\begin{eqnarray}
\label{hot_pot} \cases{ u_t=ku_{xx}, &\quad\mbox{$x\in(0, \pi),
t>0$};\vspace*{2pt}
\cr
u(0,t)=u(\pi, t)=0;&\vspace*{2pt}
\cr
u(x,0)=\phi(x),&}
\end{eqnarray}
where $k>0$ is a constant. Suppose $\phi(x)=\sum_{j=1}^{\infty
}A_j\sin jx$ for all $x\in[0, \pi]$. Then the solution of (\ref
{hot_pot}) is given by
\[
\label{ali} u(x, t)=\sum_{j=1}^{\infty}A_je^{-j^2kt}
\sin jx.
\]
See, for example, page 85 from \cite{Strauss}.
If $\sup_{j\geq0}|A_j|<\infty$, then $A_je^{-j^2kt}=O(e^{-j^2kt})$
as $j\to\infty$. Similar to the previous example, we can write $u(x, t)$ in the
form of
$a_0+\sum_{j=1}^{\infty}(a_je^{ij\theta} + b_je^{-ij\theta})$,
where $a_j$'s and $b_j$'s are complex numbers with $|a_j| \vee
|b_j|=O(e^{-j^2kt})$ as $j\to\infty$. Theorems \ref{circular_case}
and \ref{symplectic_case} can then be applied to obtain the CLT for
$a_0+\sum_{j=1}^{\infty}(a_jp_j(Z_n^{\alpha}) + b_j\overline
{p_j(Z_n^{\alpha})} )$ with $\alpha=2$ and $\alpha=\frac{1}{2}$,
respectively.
\end{example*}

To get the analogues of Theorems \ref{circular_case} and \ref
{symplectic_case} for any $\beta\ne1,2,4$, one needs to get upper
bounds for $\mathbb{E} [|p_m(Z_n^{\alpha})|^2]$ as in Propositions
\ref{mouth} and \ref{search}. It will be even more involved because
of the lack of classifications of partitions as in (\ref{romantic})
for general $\beta>0$, particularly for irrational $\beta>0$.
However, by using our method, it is possible to get upper bounds for
any $\beta=\ldots, \frac{1}{4},\frac{1}{3}, \frac{1}{2},
1,2,3,4,\ldots .$

\section{Proofs of moment inequalities in Section \texorpdfstring{\protect\ref
{moment_inequality}}{2}}\label{Proof_moment_inequality}
This section is divided into two parts. In Section~\ref{proofs}, the
necessary background of the Jack functions including their orthogonal
properties and combinatorial structures are given. With this
preparation, we prove parts (a) and (b) of Theorem~\ref{face} and
Corollary~\ref{theory}. In Section~\ref{party}, we prove part (c) of
Theorem~\ref{face} by analysis.

\subsection{\texorpdfstring{Proofs of \textup{(a)} and \textup{(b)} of Theorem \protect\ref{face} and
Corollary \protect\ref{theory}}{Proofs of (a) and (b) of Theorem 1 and
Corollary 2}}\label{proofs}

For a partition $\lambda$, the notation $\lambda'=(\lambda_1',
\lambda_2',\ldots)$ represents the conjugate partition of $\lambda$,
whose Young diagram is obtained by transposing the Young diagram of
$\lambda$.

Let us review Jack symmetric functions briefly.
We do not need the exact definition of Jack functions.
In fact, their orthogonal properties are actively used here.
For any real number $\alpha>0$ and each integer $k \ge1$,
we denote by $\Lambda^k(\alpha)$ the algebra of
symmetric functions of degree $k$
over the field $\mathbb{Q}(\alpha)$.
Recall power-sum symmetric function $p_{\rho}$ in (\ref{remote}).
The family of $p_\rho$ over partitions $\rho$ of $k$ forms
a basis on $\Lambda^k(\alpha)$.
A scalar product on $\Lambda^k(\alpha)$ is defined by
%
\begin{equation}
\label{scalar} \langle p_{\lambda}, p_{\mu}\rangle_{\alpha}=
\delta_{\lambda\mu
} \alpha^{l(\lambda)}z_{\lambda}
\end{equation}
for any partitions $\lambda$ and $\mu$ of $k$,
where $z_{\lambda}$ is as in (\ref{insect}).
Set
%
\begin{equation}\label{ski}
C_{\lambda}(\alpha)=\prod_{(i,j)\in\lambda} \bigl\{
\bigl(\alpha(\lambda_i -j ) + \lambda_j'-i+1
\bigr) \bigl(\alpha(\lambda_i -j ) + \lambda
_j'-i+\alpha\bigr) \bigr\},
\end{equation}
%
where $(i,j)$ runs over all cells of the Young diagram of $\lambda$.
By definition,
Jack functions $\{J_\lambda^{(\alpha)}\}$ form an orthogonal basis on
$\Lambda^k(\alpha)$ and
satisfy
%
\begin{equation}
\label{relation} \bigl\langle J_{\lambda}^{(\alpha)}, J_{\mu}^{(\alpha)}
\bigr\rangle_{\alpha
}= \delta_{\lambda\mu} C_{\lambda}(\alpha);
\end{equation}
see, for example, Chapter VI from \cite{Macdonald1998} or \cite{Forresterbook}.

Since both power-sum symmetric functions and Jack functions
form a basis of~$\Lambda^k(\alpha)$,
they can be mutually expanded.
Let $\theta_{\rho}^{\lambda}(\alpha)$ denote the coefficient of
$p_{\rho}$ in $J_{\lambda}^{(\alpha)}$, that is,
%
\begin{equation}
\label{Christmas} J_{\lambda}^{(\alpha)}=\sum_{\rho: |\rho|=|\lambda|}
\theta_{\rho
}^{\lambda}(\alpha)p_{\rho}.
\end{equation}
The $\theta_{\rho}^{\lambda}(\alpha)$'s are real numbers.
Inversely, let $\Theta_{\rho}^{\lambda}(\alpha)$ be the coefficient
of $J_{\lambda}^{(\alpha)}$ in $p_{\rho}$, that is,
%
\begin{equation}
\label{eve} p_{\rho}=\sum_{\lambda: |\lambda|=|\rho|}
\Theta_{\rho}^{\lambda
}(\alpha) J_{\lambda}^{(\alpha)}.
\end{equation}

\begin{lemma}\label{summer} Recalling $\theta_{\rho}^{\lambda
}(\alpha)$ in (\ref{Christmas}) and $\Theta_{\rho}^{\lambda
}(\alpha)$ in (\ref{eve}). Then, for any partitions $\lambda$ and
$\rho$ with $|\lambda|=|\rho|$, we have
%
\begin{equation}
\label{winter} \Theta_{\rho}^{\lambda}(\alpha)= \frac{\alpha^{l(\rho)}z_{\rho
}}{C_{\lambda}(\alpha)}
\theta_{\rho}^{\lambda}(\alpha).
\end{equation}
\end{lemma}

\begin{pf}
It follows from (\ref{Christmas}) and (\ref{scalar}) that
\[
\bigl\langle J_{\lambda}^{(\alpha)}, p_{\rho}\bigr
\rangle_{\alpha} = \biggl\langle \sum_{v}
\theta_{v}^{\lambda}(\alpha)p_v, p_{\rho}
\biggr\rangle_{\alpha
} = \sum_{v}
\theta_{v}^{\lambda}(\alpha)\langle p_v,
p_{\rho
}\rangle_{\alpha}= \theta_{\rho}^{\lambda}(
\alpha) \alpha ^{l(\rho)}z_{\rho}.
\]
Similarly, by (\ref{eve}) and (\ref{relation}),
\[
\bigl\langle J_{\lambda}^{(\alpha)}, p_{\rho}\bigr
\rangle_{\alpha} = \biggl\langle J_{\lambda}^{(\alpha)}, \sum
_{v}\Theta_{\rho}^{v}(\alpha
)J_v^{(\alpha)}\biggr\rangle_{\alpha} = \sum
_{v} \Theta_{\rho
}^{v}(\alpha)\bigl
\langle J_{\lambda}^{(\alpha)}, J_{v}^{(\alpha
)}\bigr
\rangle_{\alpha}= \Theta_{\rho}^{\lambda}(\alpha)
C_{\lambda
}(\alpha).
\]
These two equalities lead to (\ref{winter}).
\end{pf}

The coefficients $\theta_{\rho}^{\lambda}$'s satisfy the following
orthogonality relations ((10.31) and~(10.32) from \cite{Macdonald1998}):
%
\begin{eqnarray}
\label{Markov} \sum_{\rho}z_{\rho}
\alpha^{l(\rho)}\theta_{\rho}^{\lambda
}(\alpha)
\theta_{\rho}^{\mu}(\alpha) &=&\delta_{\lambda\mu}C_{\lambda}(
\alpha);
\nonumber
\\[-8pt]
\\[-8pt]
\nonumber
\sum_{\lambda}\frac{1}{C_{\lambda}(\alpha)}
\theta_{\rho
}^{\lambda}(\alpha)\theta_{\sigma}^{\lambda}(
\alpha) &=&\delta_{\rho\sigma}z_{\rho}^{-1}
\alpha^{-l(\rho)}.
\end{eqnarray}
In other words, if $a_{\lambda\rho}:= (z_{\rho}\alpha^{l(\rho
)}/C_{\lambda}(\alpha) )^{1/2}\theta_{\rho}^{\lambda}(\alpha)$,
then $A_m=(a_{\lambda\rho})_{|\lambda|=|\rho|=m}$ is an orthogonal
matrix of size $p(m)$
for $m\geq1$. Here, $p(m)$ is the number of partitions of $m$. The
following are some special cases of the Jack polynomials.

In other words, if $a_{\lambda\rho}:= (z_{\rho}\alpha^{l(\rho
)}/C_{\lambda}(\alpha) )^{1/2}\theta_{\rho}^{\lambda}(\alpha)$,
then $A_m=(a_{\lambda\rho})_{|\lambda|=|\rho|=m}$ is an orthogonal
matrix of size $p(m)$
for $m\geq1$. Here, $p(m)$ is the number of partitions of $m$. The
following are some special cases of the Jack polynomials.

\begin{example*} Let $\alpha=1,  s_{\lambda}$ be the
Schur polynomial and $\chi_{\mu}^{\lambda}$ the character value for
the irreducible representation of the symmetric groups. It is well
known that $J_{\lambda}^{(1)}=h(\lambda)s_{\lambda}$ with $h(\lambda
)=\sqrt{C_{\lambda}(1)}$ as the hook-length product. Further, by (7.8)
from Chapter VI of \cite{Macdonald1998} and (\ref{eve}) that
\[
\theta_{\mu}^{\lambda}(1)=\frac{h(\lambda)\chi_{\mu}^{\lambda
}}{z_{\mu}} \quad\mbox{and}\quad
\Theta_{\mu}^{\lambda}(1)=\frac{\chi_{\mu}^{\lambda}}{h(\lambda
)}.
\]
\end{example*}

\begin{example*} Let $\alpha=2$. Then $J_\lambda^{(2)}$
coincides with the zonal polynomial $Z_\lambda$.
By~(2.13) and (2.16) from Chapter VII of \cite{Macdonald1998}, we have
\[
\theta_{\mu}^{\lambda}(2)=\frac{2^k k!}{2^{l(\mu)} z_\mu} \omega
_{\mu}^{\lambda} \quad\mbox{and}\quad \Theta_{\mu}^{\lambda}(2)=
\frac{2^k k! }{h(2\lambda)}\omega_{\mu
}^{\lambda},
\]
with $k=|\lambda|=|\mu|$,
where $h(2\lambda)=C_\lambda(2)$ is the hook-length product of
$2\lambda=(2\lambda_1,2\lambda_2,\ldots)$ and
$\omega^\lambda_\mu$ is the value of the zonal spherical function of
a Gelfand pair
$(\mathfrak{S}_{2k}, \mathfrak{B}_k)$.
Here, $\mathfrak{S}_{2k}$ is the symmetric group and $\mathfrak{B}_k$
is the hyperoctahedral group in $\mathfrak{S}_{2k}$.
\end{example*}

\begin{example*}[(Example~1(a) on page~383 from \cite{Macdonald1998})]
For each partition $\rho$ of $k$, we have
%
\begin{equation}
\label{piano} \theta_{\rho}^{(k)}(\alpha)=\frac{k!}{z_{\rho}}
\alpha^{k- l(\rho
)} \quad\mbox{and}\quad \theta_{\rho}^{(1^k)}(\alpha)=
\frac
{k!}{z_{\rho}}(-1)^{k-l(\rho)}.
\end{equation}

For each partition $\lambda$ with $l(\lambda) \le n$, we define
\[
\mathcal{N}_{\lambda}^{\alpha}(n)= \prod
_{(i,j)\in\lambda}\frac
{n+(j-1)\alpha-(i-1)}{n+j\alpha-i},
\]
which is a positive real number.
As we saw in \eqref{relation},
Jack functions are orthogonal with respect to
the scalar product $\langle\cdot, \cdot\rangle_\alpha$.
We next need the second orthogonal property for them.
\end{example*}

\begin{lemma}\label{nice} Let $\lambda$ and $\mu$ be two partitions.
Let $\alpha>0$ and $n\geq1$.
Then
\begin{eqnarray*}
&& \frac{1}{(2\pi)^n} \int_{[0,2\pi)^n} J_{\lambda}^{(\alpha)}
\bigl(e^{i\theta_1},\ldots ,e^{i\theta_n}\bigr) J_{\mu}^{(\alpha)}
\bigl(e^{-i\theta_1},\ldots,e^{-i\theta_n}\bigr)
\\
&&\hspace*{50pt}\quad{} \times\prod_{1\leq p<q\leq n}\bigl|e^{i \theta
_p}-e^{i\theta_q}\bigr|^{2/\alpha}
\,d\theta_1 \cdots d \theta_n
\\
&&\qquad=  \delta_{\lambda\mu}\cdot\delta\bigl(l(\lambda)\leq n\bigr)\cdot
\frac
{\Gamma({n}/{\alpha}+1)}{\Gamma(1+{1}/{\alpha
})^n}C_\lambda(\alpha) \mathcal{N}_{\lambda}^{\alpha}(n).
\end{eqnarray*}
\end{lemma}

\begin{pf}
Since $J_\lambda^{(\alpha)}(x_1,\ldots,x_n)=0$ if $l(\lambda) > n$,
we assume $l(\lambda) \leq n$ in the following discussion. It is known
(e.g., Theorem~12.1.1 from \cite{Mehta}) that
%
\begin{equation}
\label{mehta} \frac{1}{(2\pi)^n} \int_{[0,2\pi)^n}\prod
_{1 \le p < q \le n}\bigl|e^{i \theta_p} -e^{i\theta_q}\bigr|^{2/\alpha} \,d
\theta_1 \cdots d \theta_n =\frac{\Gamma({n}/{\alpha}+1)}{\Gamma(1+{1}/{\alpha})^n}.
\end{equation}
From (10.22), (10.35) and (10.37) in \cite{Macdonald1998}, we see that
\begin{eqnarray*}
&& \frac{1}{(2\pi)^n n!   C_{\lambda}(\alpha)} \int_{[0,2\pi
)^n}J_{\lambda}^{(\alpha)}
\bigl(e^{i\theta_1},\ldots,e^{i\theta
_n}\bigr)J_{\mu}^{(\alpha)}
\bigl(e^{-i\theta_1},\ldots,e^{-i\theta_n}\bigr)
\\
&&\hspace*{88pt}\quad{} \times \prod_{1 \le p < q \le n}\bigl|e^{i \theta_p}
-e^{i\theta_q}\bigr|^{2/\alpha} \,d\theta_1 \cdots d
\theta_n
\\
&&\qquad=  \delta_{\lambda\mu}\cdot c_n \mathcal{N}_\lambda^\alpha(n),
\end{eqnarray*}
where
\[
c_n:=\frac{1}{(2\pi)^n n!}\int_{[0,2\pi)^n} \prod
_{1 \le p < q \le n}\bigl|e^{i \theta_p} -e^{i\theta_q}\bigr|^{2/\alpha} \,d
\theta_1 \cdots \,d\theta_n =\frac{1}{n!}\cdot
\frac{\Gamma(
{n}/{\alpha}+1)}{\Gamma(1+{1}/{\alpha})^n}
\]
by (\ref{mehta}).
Hence, the desired conclusion follows.
\end{pf}

\begin{prop}\label{exact} Let $\beta>0$ be a constant. Suppose
$\theta_1, \ldots, \theta_n$ have a joint density as in (\ref
{bensemble}). Let $Z_n=(e^{i\theta_1}, \ldots, e^{i\theta_n})$.
Given partitions $\mu$ and $\nu$ of weight~$K$, then
\[
\mathbb{E} \bigl[p_{\mu}(Z_n)
\overline{p_{\nu}(Z_n)} \bigr]= \alpha^{l(\mu)+l(\nu)}z_{\mu}z_{\nu}
\sum_{\lambda\vdash K:
l(\lambda)\leq n}\frac{\theta_{\mu}^{\lambda}(\alpha)\theta_{\nu
}^{\lambda}(\alpha)}{C_{\lambda}(\alpha)} \mathcal{N}_\lambda^\alpha(n).
\]
\end{prop}

\begin{pf}
Reviewing (\ref{bensemble}), by (\ref{eve}) and Lemma~\ref{nice}, we have
\[
\mathbb{E} \bigl[p_{\mu}(Z_n)\overline{p_{\nu}(Z_n)}
\bigr] =\sum_{\lambda\vdash K:  l(\lambda)\leq n}\Theta_{\mu}^{\lambda
}(
\alpha)\Theta_{\nu}^{\lambda}(\alpha)C_\lambda(\alpha)
\mathcal {N}_{\lambda}^{\alpha}(n),
\]
where $\alpha=2/\beta$. By Lemma~\ref{summer}, the above is
identical to
\[
 \alpha^{l(\mu)+l(\nu)}z_{\mu}z_{\nu}\sum
_{\lambda\vdash K:
l(\lambda)\leq n} \frac{\theta_{\mu}^{\lambda}(\alpha)\theta_{\nu}^{\lambda
}(\alpha)}{
C_{\lambda}(\alpha)}\mathcal{N}_{\lambda}^{\alpha}(n).
\]
The proof is completed.
\end{pf}

For positive integers $n$ and $K$ and real number $\alpha>0$.
Define
%
\begin{eqnarray}
\label{er} \Gamma_{n, K}^{\alpha} &=&\max_{\lambda\vdash K:  l(\lambda)\leq n}
\mathcal{N}_\lambda ^\alpha(n)
\nonumber
\\[-8pt]
\\[-8pt]
\nonumber
&=&\max_{\lambda\vdash K:  l(\lambda)\leq n}\prod_{(i,j)\in\lambda
}
\frac{n+(j-1)\alpha-(i-1)}{n+j\alpha-i};
\\
\label{three} \gamma_{n, K}^{\alpha} &=&\min_{\lambda\vdash K:  l(\lambda)\leq n}
\mathcal{N}_\lambda ^\alpha(n)
\nonumber
\\[-8pt]
\\[-8pt]
\nonumber
&=& \min_{\lambda\vdash K:  l(\lambda)\leq n}\prod_{(i,j)\in
\lambda}
\frac{n+(j-1)\alpha-(i-1)}{n+j\alpha-i}.
\end{eqnarray}

\begin{lemma}\label{sheep} Let $\alpha>0,  K\geq1$ and $\Gamma_{n,
K}^{\alpha}$ be as in (\ref{er}) and $\gamma_{n, K}^{\alpha}$ be as
in~(\ref{three}). If $n\geq K$, then $ A\leq\gamma_{n, K}^{\alpha
}\leq \Gamma_{n, K}^{\alpha} \leq B$
where $A$ and $B$ are as in (\ref{meso}). Further, if $n\geq K$, then
%
\begin{equation}
\label{star} \max_{\lambda\vdash K}\bigl| \mathcal{N}_\lambda^\alpha(n)-1
\bigr| \leq \max\bigl\{|A-1|, |B-1|\bigr\}.
\end{equation}
\end{lemma}

\begin{pf} For $\lambda\vdash K$ such that $l(\lambda)\leq n$ and
$(i,j)\in\lambda$, it is easy to check that
%
\begin{equation}
\label{silent1} 1\leq i\leq\min\{n, K\} \quad\mbox{and}\quad 1\leq j \leq K.
\end{equation}
Thus, $n+(j-1)\alpha-(i-1)\geq n-i+1>0 \mbox{ and } n+j\alpha
-i\geq j\alpha>0$. It follows that
%
\begin{equation}
\label{mama1} b_{i,j}(\alpha):=\frac{n+(j-1)\alpha-(i-1)}{n+j\alpha-i}>0.
\end{equation}
Write
%
\begin{equation}
\label{ding1} b_{i,j}(\alpha)= 1+ \frac{1-\alpha}{n+j\alpha-i}.
\end{equation}

\textit{Case} 1: $\alpha\geq1$. By (\ref{mama1}) and (\ref
{ding1}), we see that $b_{i,j}(\alpha) \in[0, 1]$ for all $\lambda
\vdash K$ such that $l(\lambda)\leq n$ and $(i,j)\in\lambda$, which
concludes $\Gamma_{n, K}^{\alpha}\leq1$.

Further, by (\ref{silent1}), $n+j\alpha-i\geq n-K + \alpha>0$ for
all $\lambda\vdash K$ such that $l(\lambda)\leq n$ and $(i,j)\in
\lambda$. Thus, noticing $1-\alpha\leq0$, we get
\[
b_{i,j}(\alpha)\geq1+\frac{1-\alpha}{n-K+\alpha}=1-\frac
{|1-\alpha|}{n-K+\alpha}>0
\]
for all $n\geq K$. This yields
\[
\gamma_{n, K}^{\alpha}\geq \biggl(1-\frac{|1-\alpha|}{n-K+\alpha
}
\biggr)^K.
\]
The above two conclusions lead to that
%
\begin{equation}
\label{grass} \biggl(1-\frac{|1-\alpha|}{n-K+\alpha} \biggr)^K \leq
\gamma_{n,
K}^{\alpha}\leq\Gamma_{n, K}^{\alpha}
\leq1
\end{equation}
for all $n\geq K$ and $\alpha\geq1$.

\textit{Case} 2: $\alpha\in(0,1]$. By (\ref{ding1}),
$b_{i,j}(\alpha)\geq 1$ for all $\lambda\vdash K$ such that
$l(\lambda)\leq n$ and $(i,j)\in\lambda$, which shows $\gamma_{n,
K}^{\alpha}\geq1$.

Moreover, by (\ref{silent1}) again, $n+j\alpha-i\geq n-K+\alpha$ for
all $\lambda\vdash K$ such that $l(\lambda)\leq n$ and $(i,j)\in
\lambda$. Thus, with $1-\alpha> 0$, we have from (\ref{ding1}) that
\[
b_{i,j}(\alpha)\leq1+\frac{1-\alpha}{n-K+\alpha}.
\]
By the definition of $\Gamma_{n, K}^{\alpha}$ and the earlier
conclusion, we get
\[
1\leq\gamma_{n, K}^{\alpha}\leq\Gamma_{n, K}^{\alpha}
\leq \biggl(1+\frac{1-\alpha}{n-K+\alpha} \biggr)^K
\]
for all $n\geq K$ and $\alpha\in(0,1]$. This and (\ref{grass}) prove
the first part of the lemma.

Finally, by the definitions in (\ref{er}) and (\ref{three}),
\[
\gamma_{n, K}^{\alpha} \leq \mathcal{N}_\lambda^\alpha(n)=
\prod_{(i,j)\in\lambda}b_{i,j}(\alpha) \leq
\Gamma_{n, K}^{\alpha}
\]
for all $\lambda\vdash K$ since $l(\lambda)\leq n$ holds
automatically if $n\geq K$. By the proved conclusion,
\[
A-1 \leq \mathcal{N}_\lambda^\alpha(n) -1\leq B-1
\]
for all $\lambda\vdash K$. This implies (\ref{star}).
\end{pf}

\begin{pf*}{Proof of (a) and (b) of Theorem~\ref{face}}
(a) By Proposition~\ref{exact}, take $\mu=\nu$ with weight $K$ to have
\[
\mathbb{E} \bigl[ \bigl|p_{\mu}(Z_n)\bigr|^2 \bigr]=
\alpha^{2l(\mu)}z_{\mu}^2\sum
_{\lambda\vdash K: l(\lambda)\leq
n}\frac{\theta_{\mu}^{\lambda}(\alpha)^2}{C_{\lambda}(\alpha
)}\mathcal{N}_\lambda^\alpha(n).
\]
Lemma~\ref{sheep} says that $\Gamma_{n, K}^{\alpha}>0$ and $\gamma
_{n, K}^{\alpha}>0$ for all $n\geq K$. By the definitions of $\Gamma
_{n, K}^{\alpha}$ in (\ref{er}) and $\gamma_{n, K}^{\alpha}$ in
(\ref{three}), since $C_{\lambda}(\alpha)>0$ for any partition
$\lambda$ and $\alpha>0$,
%
\begin{eqnarray}
\label{west2} \gamma_{n, K}^{\alpha}\cdot\alpha^{2l(\mu)}z_{\mu}^2
\sum_{\lambda\vdash K: l(\lambda)\leq n}\frac{\theta_{\mu}^{\lambda
}(\alpha)^2}{C_{\lambda}(\alpha)} &\leq &\mathbb{E} \bigl[
\bigl|p_{\mu}(Z_n)\bigr|^2 \bigr]
\nonumber
\\[-8pt]
\\[-8pt]
\nonumber
& \leq& \Gamma_{n, K}^{\alpha}\cdot\alpha^{2l(\mu)}z_{\mu}^2
\sum_{\lambda\vdash K: l(\lambda)\leq n }\frac{\theta_{\mu}^{\lambda
}(\alpha)^2}{C_{\lambda}(\alpha)}.
\nonumber
\end{eqnarray}
From assumption $n\geq K$, if $\lambda\vdash K$, we know $l(\lambda
)\leq n$ automatically. Therefore, from (\ref{Markov}) the two sums in
(\ref{west2}) are both equal to $z_{\mu}^{-1}\alpha^{-l(\mu)}$.
Consequently,
\[
\gamma_{n, K}^{\alpha} \leq\frac{\mathbb{E}  [ |p_{\mu
}(Z_n)|^2  ]}{\alpha^{l(\mu)}z_{\mu}} \leq
\Gamma_{n,
K}^{\alpha}.
\]
The conclusion (a) then follows from the first part of Lemma~\ref{sheep}.

(b) First, assume $|\mu|\ne|\nu|$. Notice
\begin{eqnarray*}
&& \mathbb{E} \bigl[p_{\mu}(Z_n)\overline{p_{\nu}(Z_n)}
\bigr]
\\
&&\qquad= \operatorname{Const}\cdot\int_0^{2\pi}\cdots\int
_{0}^{2\pi} p_{\mu}\bigl(e^{i\theta
_1},
\ldots, e^{i\theta_n}\bigr)\overline{p_{\nu}\bigl(e^{i\theta_1},
\ldots , e^{i\theta_n}\bigr)}
\\
&&\hspace*{88pt}\qquad\quad{} \times\prod_{1\leq j<k \leq n}\bigl|e^{i\theta_j}-e^{i\theta
_k}\bigr|^{\beta}
\,d\theta_1\cdots d\theta_n.
\end{eqnarray*}
For an integrable function $h(x)$, we know $\int_0^{2\pi}h(e^{ix})
\,dx=\int_b^{b+2\pi}h(e^{ix}) \,dx$ for any $b\in\mathbb{R}$. Using
the induction and the Fubini theorem, we see that
\begin{eqnarray*}
&& \mathbb{E} \bigl[p_{\mu}(Z_n)\overline{p_{\nu}(Z_n)}
\bigr]
\\
&&\qquad= \operatorname{Const}\cdot\int_b^{b+2\pi}\cdots\int
_{b}^{b+2\pi} p_{\mu
}\bigl(e^{i\theta_1},
\ldots, e^{i\theta_n}\bigr)\overline{p_{\nu
}\bigl(e^{i\theta_1},
\ldots, e^{i\theta_n}\bigr)}
\\
&&\hspace*{108pt}\qquad\quad{} \times\prod_{1\leq j<k \leq n}\bigl|e^{i\theta_j}-e^{i\theta
_k}\bigr|^{\beta}
\,d\theta_1\cdots d\theta_n.
\end{eqnarray*}
Making transform\vspace*{1pt} $\eta_j=\theta_j-b$ for $1\leq j \leq n$, noting
that $p_{\mu}(e^{ib+i\eta_1}, \ldots,\break e^{ib+i\eta_n})=e^{ib|\mu
|}p_{\mu}(e^{i\eta_1}, \ldots,   e^{i\eta_n})$ for any $b\in\mathbb
{R}$, we obtain that
\[
\mathbb{E} \bigl[p_{\mu}(Z_n)\overline{p_{\nu}(Z_n)}
\bigr]=e^{ib(|\mu|-|\nu|)}\mathbb{E} \bigl[p_{\mu}(Z_n)
\overline{p_{\nu
}(Z_n)} \bigr]
\]
for any $b\in\mathbb{R}$.
If $|\mu|\ne|\nu|$, since $b$ is arbitrary,
we then conclude
\[
\mathbb{E} \bigl[p_{\mu}(Z_n)\overline{p_{\nu}(Z_n)}
\bigr]=0
\]
for all $n\geq2$.

To prove the second part of (b), by the first part, it suffices to
prove the conclusion for $n\geq|\mu|=|\nu|=K$. Observe that
$l(\lambda)\leq n$ if $\lambda\vdash K$. Thus, it follows from
Proposition~\ref{exact} that
\begin{eqnarray*}
& & \mathbb{E} \bigl[p_{\mu}(Z_n)\overline{p_{\nu}(Z_n)}
\bigr]
\nonumber
\\
&&\qquad=  \alpha^{l(\mu)+l(\nu)}z_{\mu}z_{\nu}\sum
_{\lambda\vdash K}\frac
{\theta_{\mu}^{\lambda}(\alpha)\theta_{\nu}^{\lambda}(\alpha
)}{C_{\lambda}(\alpha)}\mathcal{N}_\lambda^\alpha(n)
\\
&&\qquad =  \alpha^{l(\mu)+l(\nu)}z_{\mu}z_{\nu} \biggl[ \sum
_{\lambda
\vdash K}\frac{\theta_{\mu}^{\lambda}(\alpha)\theta_{\nu
}^{\lambda}(\alpha)}{C_{\lambda}(\alpha)} + \sum
_{\lambda\vdash
K}\frac{\theta_{\mu}^{\lambda}(\alpha)\theta_{\nu}^{\lambda
}(\alpha)}{C_{\lambda}(\alpha)} \bigl( \mathcal{N}_\lambda^\alpha(n)-1
\bigr) \biggr]
\\
&&\qquad =  \alpha^{l(\mu)+l(\nu)}z_{\mu}z_{\nu} \sum
_{\lambda\vdash
K}\frac{\theta_{\mu}^{\lambda}(\alpha)\theta_{\nu}^{\lambda
}(\alpha)}{C_{\lambda}(\alpha)} \bigl( \mathcal{N}_\lambda^\alpha
(n)-1 \bigr),
\end{eqnarray*}
where the last identity comes from the orthogonal property in (\ref
{Markov}). Therefore,
\begin{eqnarray*}
& & \bigl|\mathbb{E} \bigl[p_{\mu}(Z_n)\overline{p_{\nu}(Z_n)}
\bigr]\bigr |
\\
&&\qquad\leq \max_{\lambda\vdash K} \bigl| \mathcal{N}_\lambda^\alpha(n)-1
\bigr| \cdot\alpha^{l(\mu)+l(\nu)}z_{\mu}z_{\nu} \sum
_{\lambda\vdash
K}\frac{|\theta_{\mu}^{\lambda}(\alpha)|\cdot|\theta_{\nu
}^{\lambda}(\alpha)|}{C_{\lambda}(\alpha)}.
\end{eqnarray*}
Now, by the Cauchy--Schwarz inequality the sum above is bounded by
\[
\biggl(\sum_{\lambda\vdash K}\frac{|\theta_{\mu}^{\lambda}(\alpha
)|^2}{C_{\lambda}(\alpha)}
\biggr)^{1/2}\cdot \biggl(\sum_{\lambda
\vdash K}
\frac{|\theta_{\nu}^{\lambda}(\alpha)|^2}{C_{\lambda
}(\alpha)} \biggr)^{1/2}=z_{\mu}^{-1/2}z_{\nu}^{-1/2}
\alpha^{-(l(\mu
)+l(\nu))/2}
\]
by (\ref{Markov}). The above two inequalities imply
\begin{eqnarray*}
\bigl|\mathbb{E} \bigl[p_{\mu}(Z_n)\overline{p_{\nu}(Z_n)}
\bigr]\bigr | & \leq& \max_{\lambda\vdash K}\bigl| \mathcal{N}_\lambda^\alpha
(n)-1 \bigr| \cdot\alpha^{(l(\mu)+l(\nu))/2}(z_{\mu}z_{\nu})^{1/2}
\\
& \leq& \max\bigl\{|A-1|,|B-1|\bigr\}\cdot\alpha^{(l(\mu)+l(\nu))/2}(z_{\mu
}z_{\nu})^{1/2}
\end{eqnarray*}
by (\ref{star}).
\end{pf*}

\begin{lemma}\label{vote} Let $A$ and $B$ be as in (\ref{meso}) with
$\beta>0$. Set $\alpha=2/\beta$. If $n\geq2K$, then
\[
\max\bigl\{|A-1|, |B-1|\bigr\} \leq\frac{6|1-\alpha|K}{n}.
\]
\end{lemma}

\begin{pf} By the definitions of $A$ and $B$, it suffices to show that,
as $n\geq2K$,
%
\begin{eqnarray}
 1- \biggl(1-\frac{\alpha-1}{n-K+\alpha} \biggr)^K &\le&\frac
{6|1-\alpha|K}{n}\qquad
\mbox{for $\alpha\geq1$;} \label{stick}
\\
 \biggl(1+\frac{1-\alpha}{n-K+\alpha} \biggr)^K-1 &\leq&\frac
{6|1-\alpha|K}{n}\qquad
\mbox{for $\alpha\in(0,1)$.} \label{pogo}
\end{eqnarray}
First, if $\alpha\geq1$, then $(\alpha-1)/(n-K+\alpha)\in[0, 1)$.
Notice $(1+x)^K\geq1+Kx$ for all $x\geq-1$ (see, e.g., Theorem~42 on
page~40 from \cite{Hardy}), we have
\[
1- \biggl(1-\frac{\alpha-1}{n-K+\alpha} \biggr)^K \leq\frac{K(\alpha
-1)}{n-K+\alpha}\leq
\frac{2K|1-\alpha|}{n}
\]
since $n-K+\alpha\geq n/2$ as $n\geq2K$. This proves (\ref{stick}).

Second, for $\alpha\in(0,1)$, it is easy to verify that $(1-\alpha
)/(n-K+\alpha) \leq1/K$ provided $n\geq2K$. By the fact that
$(1+x)^K \leq1+ 3Kx$ for all $0\leq x\leq1/K$, we obtain
\[
\biggl(1+\frac{1-\alpha}{n-K+\alpha} \biggr)^K-1\leq\frac{3(1-\alpha
)K}{n-K+\alpha} \leq
\frac{6|1-\alpha|K}{n}
\]
since $n-K+\alpha\geq n/2$ if $n\geq2K$ as used earlier. This
concludes (\ref{pogo}).
\end{pf}

\begin{pf*}{Proof of Corollary~\ref{theory}}
(a) By Theorem~\ref{face}
\[
A-1\leq \frac{\mathbb{E}  [|p_{\mu}(Z_n)|^2 ] }{\alpha
^{l(\mu)}z_{\mu}}-1\leq B-1.
\]
Thus,
\[
\biggl| \frac{\mathbb{E}  [|p_{\mu}(Z_n)|^2  ]}{\alpha
^{l(\mu)}z_{\mu}}-1 \biggr|\leq\max\bigl\{|A-1|, |B-1|\bigr\}.
\]
The conclusion (a) then follows from Lemma~\ref{vote}.

(b) The conclusion obviously holds if $|\mu|\ne|\nu|$ by (b) of
Theorem~\ref{face}. If $|\mu|=|\nu|=K$, by (b) of Theorem~\ref
{face} and Lemma~\ref{vote}, we get the desired result.
\end{pf*}

\subsection{\texorpdfstring{Proof of \textup{(c)} of Theorem \protect\ref{face}}
{Proof of (c) of Theorem 1}}\label{party}
We start the proof through a series of lemmas.

\begin{lemma}\label{Jelly} Let $\beta>0$. For positive integers $m$
and $k$ and real numbers $a_1, \ldots, a_k$, define
\[
D=\int_{0}^{\pi}\cos(2mt) \Biggl|\prod
_{i=1}^k\sin (t+a_i )
\Biggr|^{\beta} \,dt.
\]
Then $|D|\leq6(1+\beta)(\frac{k}{m})^{1\wedge\beta}$.
\end{lemma}

\begin{pf} First, since $|D|\leq\int_0^{\pi}1 \,dt=\pi$, the
conclusion obviously holds for $m=1$. Now we assume $m\geq2$. Set
$s=mt$. Then
%
\begin{eqnarray}
\label{identity} D 
&= & \frac{1}{m}\int
_0^{m\pi}\cos(2s)\Biggl |\prod
_{i=1}^k\sin \biggl(\frac{s}{m}+a_i
\biggr) \Biggr|^{\beta} \,ds
\nonumber
\\
& = & \frac{1}{m}\sum_{j=0}^{m-1}
\int_{j\pi}^{(j+1)\pi}\cos(2s) \Biggl|\prod
_{i=1}^k\sin \biggl(\frac{s}{m}+a_i
\biggr)\Biggr |^{\beta} \,ds
\nonumber
\\[-8pt]
\\[-8pt]
\nonumber
& = & \frac{1}{m}\sum_{j=0}^{m-1}
\int_{0}^{\pi}\cos(2s) \Biggl|\prod
_{i=1}^k\sin \biggl(\frac{s+j\pi}{m}+a_i
\biggr) \Biggr|^{\beta} \,ds
\\
&= &\int_0^{\pi}L_m(s)\cos(2s) \,ds,\nonumber
\end{eqnarray}
where we make a transform: $s\to s-j\pi$ in the second identity to get
the third one, and
\[
\label{bao} L_m(s)=\frac{1}{m}\sum
_{j=0}^{m-1} \Biggl|\prod_{i=1}^k
\sin \biggl(b_{ij} + \frac{s}{m} \biggr) \Biggr|^{\beta}
\]
for $0\leq s \leq\pi$ and $b_{ij}=a_i +\frac{j\pi}{m}$. Since
$(a+b)^{\beta}\leq a^{\beta} + b^{\beta}$ for any $a\geq0,  b\geq
0, \beta\in(0,1]$, and $|c^{\beta}-d^{\beta}| \leq\beta|c-d|$ for
any $c, d\in[-1,1], \beta>1$, it is not difficult to see that $|
|x|^{\beta} -|y|^{\beta}|\leq(1+\beta) | |x|-|y|  |^{1
\wedge\beta}\leq(1+\beta)|x-y|^{1 \wedge\beta}$ for any $\beta
>0$ and $x, y \in[-1, 1]$. Therefore,
%
\begin{eqnarray}
\label{sichuan} & &\Biggl |L_m(s)-\frac{1}{m}\sum
_{j=0}^{m-1} \Biggl|\prod_{i=1}^k
\sin b_{ij} \Biggr|^{\beta}\Biggr |
\nonumber
\\
&&\qquad \leq \frac{1}{m}\sum_{j=0}^{m-1}
\Biggl| \Biggl|\prod_{i=1}^k\sin
\biggl(b_{ij} + \frac{s}{m} \biggr) \Biggr|^{\beta}- \Biggl|\prod
_{i=1}^k\sin b_{ij}
\Biggr|^{\beta} \Biggr|
\\
&&\qquad \leq \frac{1+\beta}{m}\sum_{j=0}^{m-1}
\Biggl|\prod_{i=1}^k\sin \biggl(b_{ij}
+ \frac{s}{m} \biggr)-\prod_{i=1}^k
\sin b_{ij} \Biggl|^{\beta\wedge1}.\nonumber
\end{eqnarray}
Now, by the product rule, $ (\prod_{i=1}^k\sin(b_{ij} + t)
)'=\sum_{l=1}^k\cos(b_{lj}+t) \times\break \prod_{1\leq i\leq k,  i\ne l}\sin
(b_{ij} + t)$ for any $t\in\mathbb{R}$. Thus, the absolute value of
the derivative is bounded by $k$ for ant $t\in\mathbb{R}$. By the
mean-value theorem,
\[
\Biggl|\prod_{i=1}^k\sin \biggl(b_{ij}
+ \frac{s}{m} \biggr)-\prod_{i=1}^k
\sin b_{ij}\Biggr | \leq\frac{ks}{m}.
\]
This implies that the last term in (\ref{sichuan}) is controlled by
\[
\frac{1+\beta}{m}\sum_{j=0}^{m-1} \biggl(
\frac{ks}{m} \biggr)^{1\wedge
\beta}=(1+\beta) \biggl(\frac{ks}{m}
\biggr)^{1\wedge\beta}.
\]
It follows from (\ref{sichuan}) that
\[
\Biggl|L_m(s)-\frac{1}{m}\sum_{j=0}^{m-1}
\Biggl|\prod_{i=1}^k\sin b_{ij}
\Biggr|^{\beta} \Biggr| \leq(1+\beta) \biggl(\frac{ks}{m} \biggr)^{1\wedge\beta}.
\]
Set $C=\frac{1}{m}\sum_{j=0}^{m-1} |\prod_{i=1}^k\sin b_{ij}
|^{\beta}$. Notice $\int_0^{\pi}\cos(2s) \,ds=0$. From the above, we
use the simple fact that $|\cos(2s)|\leq1$ to have
\begin{eqnarray*}
\biggl|\int_0^{\pi}L_m(s)\cos(2s) \,ds \biggr| &
= & \biggl|\int_0^{\pi}C\cos(2s) \,ds+ \int
_0^{\pi}\bigl(L_m(s)-C\bigr)\cos
(2s) \,ds \biggr|
\\
& \leq& (1+\beta) \biggl(\frac{k}{m} \biggr)^{1\wedge\beta}\int
_0^{\pi}s^{1 \wedge\beta} \,ds. 
\end{eqnarray*}
Now the last integral above is bounded by $\int_0^1 1 \,ds + \int_1^{\pi}s \,ds=(\pi^2 + 1)/2\leq6$. The proof is completed by using
(\ref{identity}).
\end{pf}

\begin{lemma}\label{button} For $\beta>0$, let $f(\theta_1, \ldots,
\theta_n|\beta)$ be as in (\ref{bensemble}). Define
\begin{eqnarray*}
I(m,n)&=&\int_0^{2\pi}\cdots\int
_0^{2\pi}\cos\bigl(m(\theta_2 -\theta
_1)\bigr)
\\
&&\hspace*{54pt}{} \times f(\theta_1, \ldots, \theta_n|\beta) \,d
\theta_1\cdots d\theta_n\qquad (m\geq0, n\geq2).
\end{eqnarray*}
Then, for some constant $K=K(\beta)$, we have $|I(m,n)|\leq(K
n2^{n\beta})m^{-(1\wedge\beta)}$ for all $m\geq1$ and $n\geq2$.
\end{lemma}

\begin{pf} Evidently, since $f(\theta_1, \ldots, \theta_n|\beta)$
is a probability density function, we know
%
\begin{equation}
\label{Fateev} I(0,n)=1
\end{equation}
for all $n\geq2$. Since $|e^{ix} -e^{iy}|^2=|1-e^{i(x-y)}|^2=(1-\cos
(x-y))^2+\sin^2(x-y)=2(1-\cos(x-y))=4\sin^2 ((x-y)/2)$ for any $x, y
\in\mathbb{R}$, the probability density function in (\ref
{bensemble}) becomes
\[
\label{Einstein} f(\theta_1, \ldots, \theta_n|
\beta)=C_n\prod_{1\leq j<k \leq
n} \biggl|\sin \biggl(
\frac{\theta_j - \theta_k}{2} \biggr) \biggr|^{\beta},
\]
where $\theta_1, \ldots, \theta_n\in[0, 2\pi]$ and
\[
C_n=2^{n(n-1)\beta/2}(2\pi)^{-n}\cdot\frac{\Gamma(1+\beta
/2)^n}{\Gamma(1+\beta n/2)}.
\]
Now,
\begin{eqnarray*}
I(m,n)
&=& \int_0^{2\pi}\cdots\int_0^{2\pi}
\cos\bigl(m(\theta_2 -\theta _1)\bigr)f(
\theta_1, \ldots, \theta_n|\beta) \,d\theta_1
\cdots d\theta _n
\\
&=& C_n\int_0^{2\pi}\cdots\int
_0^{2\pi}\cos\bigl(m(\theta_2 -\theta
_1)\bigr)\\
&&\hspace*{69pt}{}\times\prod_{1\leq j<k \leq n} \biggl|\sin \biggl(
\frac{\theta_j -
\theta_k}{2} \biggr) \biggr|^{\beta} \,d\theta_2\cdots d
\theta_n \,d\theta_1.
\nonumber
\end{eqnarray*}
Making transforms $x_i=\theta_i-\theta_1$ for $i=2,3,\ldots,n$, we
obtain that
\[
I(m,n)=C_n\int_0^{2\pi}\int
_{-\theta_1}^{2\pi-\theta_1}\cdots \int_{-\theta_1}^{2\pi-\theta_1}
\cos(mx_2)\cdot G_n(x) \,dx_2\cdots
dx_n \,d\theta_1
\]
with
\[
\label{useful} G_n(x)=\prod_{i=2}^n
\biggl|\sin \biggl(\frac{x_i}{2} \biggr) \biggr|^{\beta
}\cdot\prod
_{2\leq j<k \leq n} \biggl|\sin \biggl(\frac{x_j -
x_k}{2} \biggr)
\biggr|^{\beta},
\]
where the second product is understood to be $1$ if $n=2$.
For a periodic and integrable function $h(x)$ with period $2\pi$, we
know that $\int_b^{b+2\pi}h(x) \,dx=\int_0^{2\pi}h(x) \,dx$. By
induction and the Fubini theorem, we have
%
\begin{eqnarray}
I(m,n)&= & C_n\int_0^{2\pi}\cdots\int
_{0}^{2\pi}\cos(mx_2)\cdot
G_n(x) \,dx_2\cdots dx_n \,d\theta_1
\nonumber
\\
&= & (2\pi)C_n\int_{0}^{2\pi}\cdots\int
_{0}^{2\pi}\cos(mx_2)\cdot
G_n(x) \,dx_2\cdots dx_n\label{Lamar}
\\
&= & (2\pi)C_n\int_{0}^{2\pi}\cdots\int
_{0}^{2\pi}\cos (mx_2)J_n(x)
H_n(x) \,dx_2\cdots dx_n,\label{hotwater}
\end{eqnarray}
where $G_n(x)=J_n(x)H_n(x)$ and
\begin{eqnarray*}
H_n(x)=\cases{ \displaystyle\prod_{i=3}^n
\biggl|\sin \biggl(\frac{x_i}{2} \biggr) \biggr|^{\beta}\cdot \prod
_{3\leq j<k \leq n} \biggl|\sin \biggl(\frac{x_j - x_k}{2} \biggr)
\biggr|^{\beta}, &\quad  $\mbox{if $n\geq4$}$;\vspace*{2pt}
\cr
\displaystyle\biggl|\sin \biggl(
\frac{x_3}{2} \biggr) \biggr|^{\beta}, & \quad $\mbox{if $n= 3$}$;\vspace*{2pt}
\cr
1, &\quad  $\mbox{if $n=2$},$}
\end{eqnarray*}
and
\begin{eqnarray*}
J_n(x)=\cases{ \displaystyle\biggl|\sin \biggl(\frac{x_2}{2} \biggr)
\biggr|^{\beta}\prod_{i=3}^n \biggl|\sin
\biggl(\frac{x_2 - x_i}{2} \biggr) \biggr|^{\beta}, &\quad $\mbox{if $n\geq3$}$;
\vspace*{2pt}
\cr
\displaystyle\biggl|\sin \biggl(\frac{x_2}{2} \biggr) \biggr|^{\beta}, &\quad $
\mbox{if $n=2$.}$}
\end{eqnarray*}
In particular,
%
\begin{equation}
\label{Santafe} I(m, 2)=2\pi C_2\int_0^{2\pi}
\cos(mx_2)J_2(x) \,dx_2.
\end{equation}
Taking $m=0$ in (\ref{Lamar}), we know from (\ref{Fateev}) that
\[
\int_0^{2\pi}\cdots\int_{0}^{2\pi}
\prod_{i=2}^n \biggl|\sin \biggl(
\frac{x_i}{2} \biggr) \biggr|^{\beta}\cdot\prod
_{2\leq j<k \leq n} \biggl|\sin \biggl(\frac{x_j - x_k}{2} \biggr)
\biggr|^{\beta} \,dx_2 \,dx_3\cdots dx_n=
\frac{1}{2\pi C_n}
\]
for all $n\geq2$, where the second product above is understood to be
$1$ if $n=2$. This implies
%
\begin{equation}
\label{back} \int_{0}^{2\pi}\cdots\int
_0^{2\pi} H_n(x) \,dx_3
\cdots dx_n=\frac
{1}{2\pi C_{n-1}}
\end{equation}
for all $n\geq3$. Now, recalling the definition of $J_n(x)$, let
$t=x_2/2$, we have
\[
\int_{0}^{2\pi}\cos(mx_2)J_n(x)
\,dx_2=2\int_{0}^{\pi}\cos (2mt) \Biggl|\prod
_{i=1}^{n-1}\sin (t+a_i )
\Biggr|^{\beta} \,dt
\]
for all $n\geq2$, where $a_1=0,  a_i=-x_{i+1}/2$ for $i=2,\ldots
,n-1$. By Lemma~\ref{Jelly},
%
\begin{equation}
\label{regular} \biggl|\int_{0}^{2\pi}
\cos(mx_2)J_n(x) \,dx_2\biggr |\leq12(1+\beta )
\biggl(\frac{n}{m} \biggr)^{1\wedge\beta}
\end{equation}
for all $n\geq2$. Therefore, this and (\ref{Santafe}) imply that for
some constant $K_1=K_1(\beta)$,
%
\begin{equation}
\label{bar} \bigl|I(m,2)\bigr| \leq\frac{K_1}{m^{(1\wedge\beta)}}.
\end{equation}
Now assume $n\geq3$. By (\ref{hotwater}) and (\ref{regular}), and
then (\ref{back}), we obtain
%
\begin{eqnarray}
\label{dong} \bigl|I(m,n)\bigr| &\leq& 24\pi(1+\beta)C_n \biggl(
\frac{n}{m} \biggr)^{1\wedge
\beta}\int_{0}^{2\pi}
\cdots\int_0^{2\pi} H_n(x)
\,dx_3\cdots dx_n
\nonumber
\\[-8pt]
\\[-8pt]
\nonumber
& = & 12(1+\beta) \biggl(\frac{n}{m} \biggr)^{1\wedge\beta}
\frac
{C_n}{C_{n-1}}
\end{eqnarray}
for all $n\geq3$. Now,
%
\begin{equation}
\label{li} \frac{C_n}{C_{n-1}}=\frac{\Gamma(1+\beta/2)}{2\pi}\cdot\frac
{\Gamma(1+\beta n/2-\beta/2)}{\Gamma(1+\beta n/2)}\cdot
2^{(n-1)\beta}
\end{equation}
for all $n\geq3$. By Lemma~2.4 from \cite{Dong_Jiang_Li}, there
exists a constant $K_2=K_2(\beta)$ such that
\[
\frac{\Gamma(1+\beta n/2-\beta/2)}{\Gamma(1+\beta n/2)} \leq\frac
{K_2}{n^{\beta/2}}
\]
for all $n\geq1$. This, (\ref{dong}) and (\ref{li}) imply that there
exists a constant $K=K(\beta)$ such that
\[
\bigl|I(m,n)\bigr| \leq K\cdot \biggl(\frac{n}{m} \biggr)^{1\wedge\beta}\cdot
\frac{1}{n^{\beta/2}}\cdot2^{n\beta}= K n^{(1\wedge\beta) -\beta
/2}\frac{2^{n\beta}}{m^{1\wedge\beta}}\leq
K\frac{n2^{n\beta
}}{m^{1\wedge\beta}}
\]
for all $n\geq3$. This together with (\ref{bar}) proves the lemma.
\end{pf}

\begin{pf*}{Proof of (c) of Theorem~\ref{face}}
Observe that, for any
real numbers $x_1, \ldots, x_n$,
\begin{eqnarray*}
\Biggl|\sum_{j=1}^ne^{ix_j}
\Biggr|^2 & = & \sum_{j=1}^ne^{ix_j}
\cdot \sum_{j=1}^ne^{-ix_j}
\\
&= & n + \sum_{j\ne k}e^{i(x_j-x_k)}= n + \sum
_{1\leq j< k\leq n} \bigl(e^{i(x_j-x_k)} + e^{-i(x_j-x_k)}
\bigr)
\\
& = & n+ 2\sum_{1\leq j< k\leq n}\cos(x_j-x_k).
\end{eqnarray*}
Thus, by the symmetry of $f(\theta_1, \ldots, \theta_n|\beta)$,
%
\begin{eqnarray}
\label{her_him} \mathbb{E} \bigl[ \bigl|p_m(Z_n)\bigr|^2
\bigr]&=&\mathbb{E} \Biggl[ \Biggl|\sum_{j=1}^ne^{im\theta_j}
\Biggr|^2 \Biggr]
\nonumber
\\[-8pt]
\\[-8pt]
\nonumber
&=& n + n(n-1)\cdot\mathbb{E}\bigl[\cos\bigl\{m(\theta_1-
\theta_2)\bigr\}\bigr].
\end{eqnarray}
The conclusion then follows from Lemma~\ref{button}.
\end{pf*}

\section{Proofs of central limit theorems in Section \texorpdfstring{\protect\ref
{CLT_introduction}}{3}}\label{Proof_CLT}

Before proving the central limit theorems, we will spend a lot efforts
in studying the second moments, which enable us to reduce the infinite
Fourier series in Theorems \ref{circular_case} and \ref
{symplectic_case} to finite sums, and hence we can apply the moment
inequalities stated in Section~\ref{moment_inequality}. We will prove
Proposition~\ref{mouth} in Section~\ref{red-blue}, and Proposition~\ref{search} in Section~\ref{black_yellow}. All of the central limit
theorems will be proved in Section~\ref{green_brown}. We start with
the combinatorial structure of the second moment.

Review that $Z_n^{\alpha}=(e^{i\theta_1}, \ldots, e^{i\theta_n})$
follow the $\beta$-circular ensemble with $\alpha=\frac{2}{\beta}$.
Its probability density function is given in (\ref{bensemble}).
Following our notation, $p_m(Z_n^\alpha)=\sum_{j=1}^ne^{im\theta_j}$
for any integer $m\geq0$. We know from Proposition~\ref{exact} that
%
\begin{equation}
\label{please} \mathbb{E} \bigl[\bigl|p_m\bigl(Z_n^\alpha
\bigr)\bigr|^2\bigr]= \alpha^2 m^2 \sum
_{\lambda\vdash m
:l(\lambda) \le n} \frac{\theta^\lambda_{(m)}(\alpha)^2}{C_\lambda(\alpha)}
 \mathcal{N}_\lambda^\alpha(n),
\end{equation}
%
where
%
\begin{eqnarray}
\label{earth} \mathcal{N}_{\lambda}^{\alpha}(n) & = & \prod
_{(i,j)\in\lambda
}\frac{n+(j-1)\alpha-(i-1)}{n+j\alpha-i}
\nonumber
\\[-8pt]
\\[-8pt]
\nonumber
& = & \prod_{(i,j)\in\lambda} \biggl(1+ \frac{1-\alpha}{n+j\alpha
-i}
\biggr).
\end{eqnarray}
We also know the following formula (page~383 from \cite{Macdonald1998}):
For each $\lambda\vdash m$,
%
\begin{equation}
\theta^\lambda_{(m)}(\alpha)= \mathop{\prod_{(i,j) \in\lambda}}_{ (i,j) \neq(1,1)} \bigl(\alpha(j-1)-(i-1)\bigr),
\end{equation}
where the product runs over all boxes of Young diagram $\lambda$,
except the $(1,1)$-box.

\subsection{Proof of Proposition \texorpdfstring{\protect\ref{mouth}}{1}}\label{red-blue}
Let us first evaluate $\theta^\lambda_{(m)}(2)$ and $C_{\lambda
}(2)$. Suppose $\alpha=2$.
The $(3,2)$th box in the Young diagram $\lambda$
gives $\alpha(j-1)-(i-1)=2 \cdot(2-1)-(3-1)=0$,
and hence $\theta^\lambda_{(m)}(2)$ vanishes
if $\lambda$ has the $(3,2)$-box.
In other words,
$\theta^\lambda_{(m)}(2)$ vanishes unless
$\lambda_3 \le1$.
Denote by $\mathcal{P}_m^{(2)}(n)$
the set of such partitions of $m$ with lengths $\le n$:
%
\begin{equation}
\label{romantic} \mathcal{P}_m^{(2)}(n)=\bigl\{\lambda=(
\lambda_1,\lambda_2,\ldots ,\lambda_n) \vdash
m; \lambda_3 \le1\bigr\}.
\end{equation}
The elements in $\mathcal{P}_m^{(2)}(n)$ can be classified into the
following three categories.
\begin{longlist}[1.]
\item[1.] The one-row partition $(m)$;
\item[2.] A two-row partition $(m-r,r)$ with
$r=1,2, \ldots, [\frac{m}{2}]$;
\item[3.]$\lambda=(r,s,1^{m-r-s})$ with
$r \ge s \ge1$ and $3\leq l(\lambda) =m-r-s+2 \le n$.
\end{longlist}
For each case, the quantity $\theta^\lambda_{(m)}(2)
= \prod_{\scriptsize
\begin{array}{c} (i,j) \in\lambda\\ (i,j) \neq(1,1)
\end{array}
}
(2j-i-1)$
is computed as follows:
%
\begin{eqnarray}
\theta^{(m)}_{(m)}(2)&=& 2 \cdot4 \cdots(2m-2)
=2^{m-1} \cdot (m-1)!;\label{bridge}
\\
\theta^{(m-r,r)}_{(m)}(2) &=& (-1) 2^{m-2r} (m-r-1)!
\cdot\frac
{(2r-2)!}{(r-1)!}; \label{cub}
\\
\theta^{(r,s,1^{m-r-s})}_{(m)}(2)&=& (-1)^{m-r-s+1}
\cdot2^{r-s} \cdot(r-1)! \label{pillow}
\nonumber
\\[-8pt]
\\[-8pt]
\nonumber
&&{}\times\frac{(2s-2)!}{(s-1)!} \cdot (m-r-s+1)!.
\end{eqnarray}


Now we study $C_\lambda(2)$. Note that $C_\lambda(2)$ coincides with
the hook-length product of $2\lambda=(2\lambda_1,2\lambda_2,\ldots)$.
The hook-length product of
$\lambda$ is computed in Section~6
from \cite{Matsumoto_Collins}:
\begin{longlist}[1.]
\item[1.]$C_{(m)}(2)= (2m)!$;
\item[2.]$C_{(m-r,r)}(2)= \frac{(2r)!   (2m-2r+1)!}{2m-4r+1}$;
\item[3.]$C_{(r,s,1^{m-r-s})}(2)=
(m+r-s+1)(m+r-s)(m-r+s)(m-r+s-1) \cdot(m-r-s+1)!   (m-r-s)!
\frac{(2r-1)!   (2s-2)!}{2r-2s+1}$.
\end{longlist}

Hence, the term
$[(\alpha^{l(\mu)} z_\mu)^2
\frac{\theta^{\lambda}_{\mu}(\alpha)^2}{C_\lambda(\alpha)}]_{
\mu=(m),  \alpha=2}=
4 m^2 \frac{\theta^{\lambda}_{(m)}(2)^2}{C_\lambda(2)}$
is given below.
%
\begin{eqnarray}
&&4 m^2 \frac{\theta^{(m)}_{(m)}(2)^2}{C_{(m)}(2)} =\frac{2^{2m}(m!)^2}{(2m)!};\label{light}
\\
& &4 m^2 \frac{\theta^{(m-r,r)}_{(m)}(2)^2}{C_{(m-r,r)}(2)}= \frac{{2r\choose r}}{{2(m-r)\choose m-r}} \cdot
\frac{2^{2m-4r} m^2 (2m-4r+1)}{
(m-r)^2 (2r-1)^2 (2m-2r+1)};\label{mouse}
\\
&&4 m^2 \frac{\theta^{(r,s,1^{m-r-s})}_{(m)}(2)^2}{C_{(r,s,1^{m-r-s})}(2)} \nonumber\label{lantern}
\\
&&\qquad =\frac{4m^2 (m-r-s+1)}{(m+r-s+1)(m+r-s)(m-r+s)(m-r+s-1)}
\\
&&\qquad\quad{} \cdot\frac{2^{2r-2s} [(r-1)!]^2 (2r-2s+1)
{\cdot (2s-2)!}}{[(s-1)!]^2 (2r-1)!}.
\nonumber
\end{eqnarray}

\emph{Note: through the rest of the paper, $C$ stands for a
generic
constant which may change from line to line.}

\begin{lemma}\label{mutton} Recall $\mathcal{N}_{\lambda}^{2}(n)$ as
in (\ref{earth}). Then there exists a universal constant $K\in
(0,\infty)$ such that $\mathcal{N}_{\lambda}^{2}(n) \leq K \sqrt
{\frac{n}{m}}$ uniformly for all $m, n$ and all $\lambda$ satisfying:
\begin{longlist}[(iii)]
\item[(i)] $\lambda=(m)$ and $m\geq n\geq1$,

\item[(ii)] $\lambda=(m-r, r)$ with $1\leq r \leq m/2$ and $m\geq
n\geq2$ or

\item[(iii)] $\lambda=(r,s,1^{m-r-s})$ with $r \ge s \ge1$, $3\leq
m-r-s+2 \le n$ and $m\geq n$.
\end{longlist}
\end{lemma}

\begin{pf}
The following basic estimate will be used several times.
%
\begin{equation}
\label{shadow} \log\frac{l}{k}\leq\sum_{j=k}^{l}
\frac{1}{j} \leq1+\log\frac{l}{k}
\end{equation}
for all $1\leq k\leq l$. It is obviously true if $k=l$. Now, for $1\leq
k< l$,
\[
\sum_{j=k}^{l}\frac{1}{j}\leq1+
\sum_{j=k+1}^{l}\int_{j-1}^j
\frac
{1}{x} \,dx=1+\int_k^l
\frac{1}{x} \,dx=1+\log\frac{l}{k}.
\]
Similarly,
\[
\sum_{j=k}^{l}\frac{1}{j}\geq
\sum_{j=k}^{l}\int_{j}^{j+1}
\frac
{1}{x} \,dx=\int_k^{l+1}
\frac{1}{x} \,dx \geq\log\frac{l}{k}.
\]

(i) Since $\lambda=(m)$, we have from (\ref{earth}) and the fact $1-x
\leq e^{-x}$ for all $x\in\mathbb{R}$ that
%
\begin{eqnarray}\label
{library}
\mathcal{N}_{\lambda}^{2}(n) &=&\prod
_{(i,j)\in\lambda} \biggl(1-\frac{1}{n+2j-i} \biggr)
 =  \prod_{j=1}^m \biggl(1-
\frac{1}{n+2j-1} \biggr)
\nonumber
\\[-8pt]
\\[-8pt]
\nonumber
& \leq& \exp \Biggl(-\frac{1}{2}\sum
_{j=1}^m\frac{1}{n-1+j} \Biggr)
\nonumber
\end{eqnarray}
since $n+2j-1\leq2(n-1+j)$. From (\ref{shadow}), we get that
\[
\label{sum} \sum_{j=1}^m
\frac{1}{n-1+j}=\sum_{j=n}^{n+m-1}
\frac{1}{j}\geq \log\frac{n+m-1}{n}
 \geq \log\frac{m}{n}
\]
for all $m\geq n\geq1$. This gives that $\mathcal{N}_{\lambda
}^{2}(n)\leq\sqrt{\frac{n}{m}}$ for any $m\geq n\geq1$.

(ii) Now, $\lambda=(m-r, r)$ with $1\leq r \leq m/2$ and $n\geq2$.
Recall (\ref{library}). We have
\begin{eqnarray*}
\mathcal{N}_{\lambda}^{2}(n) &=&\prod
_{(i,j)\in\lambda} \biggl(1-\frac{1}{n+2j-i} \biggr)
\\
& =& \prod_{j=1}^{m-r} \biggl(1-
\frac{1}{n+2j-1} \biggr)\cdot\prod_{j=1}^{r}
\biggl(1-\frac{1}{n+2j-2} \biggr)
\\
& \leq& \exp \Biggl(-\frac{1}{2}\sum
_{j=1}^{m-r}\frac{1}{n-1+j} -\frac{1}{2}
\sum_{j=1}^{r}\frac{1}{n-2
+j} \Biggr)
\end{eqnarray*}
by the inequality $n+2j-i\leq2(n-i+j)$ for $i=1,2$. Hence,
\begin{eqnarray*}
-2\log\mathcal{N}_{\lambda}^{2}(n) & \geq& \sum
_{j=n}^{m+n-r-1}\frac{1}{j} +\sum
_{j=n-1}^{n+r-2}\frac{1}{j}
\\
& \geq& \log \biggl(\frac{m+n-r-1}{n}\cdot\frac{n+r-2}{n-1} \biggr)
\end{eqnarray*}
for any $1\leq r \leq m/2$ by (\ref{shadow}). Notice $\frac
{m+n-r-1}{n}\geq\frac{m}{2n}$ and $\frac{n+r-2}{n-1}\geq1$ since
$1\leq r \leq m/2$. We then have
\[
\mathcal{N}_{\lambda}^{2}(n) \leq2 \sqrt{\frac{n}{m}}.
\]

(iii) In this case, $\lambda=(r,s,1^{m-r-s})$ with
$r \ge s \ge1$ and $3\leq l(\lambda) =m-r-s+2 \le n$ and $m\geq n$.
First, these restrictions imply
%
\begin{equation}
\label{tulip} r\geq\frac{m-n}{2}+1,\qquad m-r\geq2 \mbox{ and } n\geq3.
\end{equation}
Now,
\begin{eqnarray*}
\mathcal{N}_{\lambda}^{2}(n) &=&\prod
_{(i,j)\in\lambda} \biggl(1-\frac{1}{n+2j-i} \biggr)
\\
& = & \prod_{j=1}^r \biggl(1-
\frac{1}{n+2j-1} \biggr)\cdot\prod_{j=1}^s
\biggl(1-\frac{1}{n+2j-2} \biggr)\\
&&{}\times\prod_{i=3}^{m-r-s+2}
\biggl(1-\frac{1}{n+2-i} \biggr)
\\
& \leq& \exp \Biggl(-\frac{1}{2}\sum
_{j=1}^r\frac{1}{n-1+j}-\frac
{1}{2}
\sum_{j=1}^s\frac{1}{n-2+j}-
\frac{1}{2}\sum_{i=3}^{m-r-s+2}
\frac{1}{n-i+1} \Biggr)
\end{eqnarray*}
by the inequality $n+2j-i \leq2 (n+j-i)$ for all $j\geq1$ and $i\leq
m-r-s+2 \le n$.
Rearranging the indices in the sums and using (\ref{shadow}), we
obtain that
\begin{eqnarray*}
-2\log\mathcal{N}_{\lambda}^{2}(n) & \geq& \sum
_{j=n}^{n+r-1}\frac{1}{j} + \sum
_{j=n-1}^{n+s-2}\frac
{1}{j} + \sum
_{j=n+r+s-m-1}^{n-2}\frac{1}{j}
\\
& \geq& \log\frac{n+r-1}{n}\cdot\frac{n+s-2}{n-1}\cdot\frac
{n-2}{n+r+s-m-1}
\\
& \geq& \log\frac{(n+r-1)(n+s-2)}{2n(n+r+s-m)}
\end{eqnarray*}
since $\frac{n-2}{n-1}\geq\frac{1}{2}$ by (\ref{tulip}). Equivalently,
\begin{eqnarray*}
\label{Santa} \mathcal{N}_{\lambda}^{2}(n) & \leq& \sqrt{
\frac{2n(n+r+s-m)}{(n+r-1)(n+s-2)}}
\nonumber
\\
& \leq& 2\sqrt{\frac{n}{m}}\cdot\sqrt{\frac
{n+r+s-m}{n+s-2}}
\nonumber
\\
& \leq& 2\sqrt{\frac{n}{m}}
\end{eqnarray*}
since $n+r-1\geq n+(m-n)/2\geq m/2$ and $\frac{n+s-(m-r)}{n+s-2}\leq
1$ by (\ref{tulip}).
\end{pf}

\begin{lemma}\label{cucumber} Let $m, r, s$ be positive integers such that
$r \ge s \ge1$ and $m>r+s$. Set $\mu=(m)$ and $\lambda=(r,s,
1^{m-r-s})$. Then there exists a universal constant $K>0$ such that
\begin{eqnarray*}
\label{below} & & m^2 \frac{\theta^{\lambda}_{\mu}(2)^2}{C_\lambda(2)} \leq K\cdot
\frac
{1}{m-r+s}\cdot\sqrt{\frac{r}{s}}.
\end{eqnarray*}
Further, if $r\geq2s$ then
\[
\label{top} m^2 \frac{\theta^{\lambda}_{\mu}(2)^2}{C_\lambda(2)} \geq\frac
{1}{K}\cdot
\frac{m-r-s}{(m-r+s)^2}\cdot\sqrt{\frac{r}{s}}.
\]
\end{lemma}

\begin{pf} From (\ref{lantern}), we see that
%
\begin{eqnarray}\label
{mama}
&& m^2 \frac{\theta^{\lambda}_{\mu}(2)^2}{C_\lambda(2)}
\nonumber
\\
&&\qquad= \frac{m^2 (m-r-s+1)}{(m+r-s+1)(m+r-s)(m-r+s)(m-r+s-1)}\\
&&\qquad\quad{}\times \frac{2^{2r-2s} [(r-1)!]^2
(2s-2)!
(2r-2s+1)}{[(s-1)!]^2 (2r-1)!}
\nonumber
\\
&&\qquad\leq \frac{ (m-r-s+1)}{(m-r+s-1)^2}\cdot \frac{2^{2r-2s} [(r-1)!]^2 (2s-2)! (2r-2s+1)}{[(s-1)!]^2 (2r-1)!} \label{copy}
\end{eqnarray}
since $m+r-s+1\geq m$ and $m+r-s\geq m$. Now, write
%
\begin{eqnarray}\label{shark}
&&\frac{2^{2r-2s} [(r-1)!]^2
(2s-2)!
(2r-2s+1)}{[(s-1)!]^2 (2r-1)!}
\nonumber\\
&&\qquad=  \frac{2^{2r-2s} (2s-2)! [r!]^2}{[(s-1)!]^2 (2r)!}\cdot\frac
{2(2r-2s+1)}{r}
\\
&&\qquad\leq 4 \cdot\frac{2^{2r-2s} {2s-2\choose s-1}}{{2r\choose r}}
\nonumber
\end{eqnarray}
due to the fact that $\frac{2(2r-2s+1)}{r}\leq4$. We regard
${0\choose0}=1$. The Stirling formula says that
%
\begin{equation}
\label{wolve} 1< \frac{k!}{\sqrt{2\pi k}k^k e^{-k}} < 2
\end{equation}
for all $k\geq1$. It is easy to check from (\ref{wolve}) that there
exists a universal constant $K>0$ such that
%
\begin{equation}
\label{rubber} \frac{1}{K}\cdot\frac{2^{2k}}{\sqrt{k}}\leq\pmatrix{2k\cr k} \leq
K\cdot\frac{2^{2k}}{\sqrt{k}}
\end{equation}
for all $k\geq1$. We claim that
%
\begin{equation}
\label{spring} \frac{2^{2r-2s} {2s-2\choose s-1}}{{2r\choose r}} \leq C\sqrt{\frac{r}{s}}
\end{equation}
for all $r\geq s \geq1$. In fact, if $s=1$,
\[
\frac{2^{2r-2s} {2s-2\choose s-1}}{{2r\choose r}}=\frac
{2^{2r-2}}{{2r\choose r}}\leq C\sqrt{r}=C\sqrt{
\frac{r}{s}}
\]
by (\ref{rubber}). If $r\geq s \geq2$, by (\ref{rubber}) again,
\[
\frac{2^{2r-2s} {2s-2\choose s-1}}{{2r\choose r}}\leq C\sqrt{\frac
{r}{s-1}}\leq2C \sqrt{
\frac{r}{s}}.
\]
So (\ref{spring}) holds. Hence,
this and (\ref{copy}) imply that
\begin{eqnarray*}
m^2 \frac{\theta^{\lambda}_{\mu}(2)^2}{C_\lambda(2)} &\leq& C\cdot\frac{ m-r-s+1}{(m-r+s-1)^2}\cdot\sqrt{
\frac{r}{s}}
\\
& \leq& C\cdot\frac{1}{m-r+s}\cdot\sqrt{\frac{r}{s}}
\end{eqnarray*}
since $m-r-s+1\leq m-r+s-1$ and $m-r+s-1\geq\frac{1}{2}(m-r+s)$.

Now we prove the lower bound. By the fact $r\leq m$ it is seen that
$m+r-s+1\leq2m$. Therefore, by (\ref{mama}) and (\ref{shark}),
\begin{eqnarray*}
 m^2 \frac{\theta^{\lambda}_{\mu}(2)^2}{C_\lambda(2)}
&\geq& \frac{1}{4}\cdot\frac{m-r-s}{(m-r+s)^2}\cdot \frac{2^{2r-2s} [(r-1)!]^2 (2s-2)! (2r-2s+1)}{[(s-1)!]^2 (2r-1)!}
\\
&= & \frac{1}{4}\cdot\frac{m-r-s}{(m-r+s)^2}\cdot\frac{2^{2r-2s}
{2s-2\choose s-1}}{{2r\choose r}}\cdot
\frac{2(2r-2s+1)}{r}.
\end{eqnarray*}
The condition $r\geq2s$ implies that $\frac{2(2r-2s+1)}{r}\geq2$. By
(\ref{rubber}) again,
\[
\frac{2^{2r-2s} {2s-2\choose s-1}}{{2r\choose r}} \geq C\sqrt{\frac{r}{s}}.
\]
We complete the proof.
\end{pf}

\begin{pf*}{Proof of Proposition~\ref{mouth}}
Look at (a) of Theorem~\ref{face}, $B=1$ since $\alpha=2$. It follows
that $\mathbb{E}  [ |p_{\mu}(Z_n)|^2  ] \leq2 m$ for $1\leq
m\leq n$. So, in the rest of the paper, we only need to study the case
for $m> n \geq2$.

Review (\ref{please}),
%
\begin{equation}
\label{Norton} \mathbb{E} \bigl[\bigl|p_m\bigl(Z_n^2
\bigr)\bigr|^2\bigr]= 4 m^2 \sum_{\lambda\vdash m
:l(\lambda) \le n}
\frac{\theta^\lambda_{(m)}(2)^2}{C_\lambda(2)} \mathcal{N}_{ \lambda}^2(n).
\end{equation}
To study this quantity, we will differentiate the three cases for
$\lambda$ in the sum as appeared earlier.

\textit{Case} 1: $\lambda=(m)$.
By (\ref{light}) and (\ref{wolve}),
%
\begin{equation}
\label{crack} 4 m^2 \frac{\theta^{\lambda}_{(m)}(2)^2}{C_\lambda(2)}=\frac
{2^{2m}(m!)^2}{(2m)!} <
\frac{2^{2m}(2\sqrt{2\pi m}  m^m
e^{-m})^2}{\sqrt{4\pi m}  (2m)^{2m}e^{-2m}} <C \sqrt{m}.
\end{equation}
Hence, by (i) of Lemma~\ref{mutton},
\begin{equation}
\label{fan} 4 m^2 \frac{\theta^{\lambda}_{(m)}(2)^2}{C_\lambda(2)}\mathcal {N}_{\lambda}^{2}(n)
\leq C \sqrt{n}
\end{equation}
for any $m\geq n\geq1$ and $\lambda=(m)$.

\textit{Case} 2: $\lambda=(m-r, r)$ \textit{with} $1\leq r \leq m/2$. First, by
(\ref{mouse}),
%
\begin{equation}
\label{vain} 4 m^2 \frac{\theta^{\lambda}_{(m)}(2)^2}{C_{\lambda}(2)} = \frac{{2r\choose r}}{{2(m-r)\choose m-r}}
\cdot\frac{2^{2m-4r} m^2 (2m-4r+1)}{
(m-r)^2 (2r-1)^2 (2m-2r+1)}.
\end{equation}
By using the fact $1\leq r \leq m/2$, we have that
$(m-r)^2r^2(2m-2r+1)\geq m^3r^2/4$ and $2^{2m-4r}m^2(2m-4r+1)\leq
2\cdot2^{2m-4r}m^3$. It follows that the last ratio in (\ref{vain})
is dominated by $8\cdot2^{2m-4r}/r^2$. Thus, by (\ref{rubber}),
%
\begin{eqnarray}
\label{Jerry} 4 m^2 \frac{\theta^{\lambda}_{(m)}(2)^2}{C_{\lambda}(2)} &\leq& C
\frac{\sqrt{m-r}}{2^{2m-2r}}\cdot\frac{2^{2r}}{\sqrt{r}}\cdot \frac{2^{2m-4r}}{r^2}
\nonumber
\\[-8pt]
\\[-8pt]
\nonumber
& = & C \frac{\sqrt{m-r}}{r^{5/2}} \leq\frac{C}{r^{5/2}}\sqrt{m}
\end{eqnarray}
for all $1\leq r \leq m/2$. It follows from (ii) of Lemma~\ref{mutton} that
%
\begin{eqnarray}
\label{rice} & & 4m^2\sum_{\lambda=(m-r, r),  1\leq r \leq m/2}
\frac{\theta
_{(m)}^{\lambda}(2)^2}{C_{\lambda}(2)}\mathcal{N}_\lambda ^2(n)
\nonumber
\\
&&\qquad \leq C\cdot\sum_{1\leq r \leq m/2}\frac{1}{r^{5/2}}\sqrt{m}
\cdot\sqrt{\frac{n}{m}}
\\
&& \qquad \leq C\cdot \Biggl(\sum_{r=1}^{\infty}
\frac{1}{r^{5/2}} \Biggr)\sqrt{n}.\nonumber
\end{eqnarray}

\textit{Case }3: $\lambda=(r,s,1^{m-r-s})$ \textit{with}
$r \ge s \ge1$ \textit{and} $3\leq l(\lambda) =m-r-s+2 \le n$.
From (iii) of Lemma~\ref{mutton} and the first assertion of Lemma~\ref
{cucumber}, we get that
%
\begin{eqnarray}
\label{apple} & & 4m^2\sum_{\lambda=(r,s, 1^{m-r-s})}
\frac{\theta_{(m)}^{\lambda
}(2)^2}{C_{\lambda}(2)}\mathcal{N}_\lambda^2(n)
\nonumber
\\[-8pt]
\\[-8pt]
\nonumber
&&\qquad \leq C\sqrt{n}\sum_{r,s}\frac{1}{m-r+s}\cdot
\frac{1}{\sqrt{s}},
\end{eqnarray}
where both sums are taken over all possible $r\geq s \geq1$ with
$3\leq l(\lambda) =m-r-s+2 \le n$. These restrictions imply that
$s+1\leq m-r\leq s+n$ and hence $2s+1\leq m-r+s\leq2s+n$. It follows
that the last sum in (\ref{apple}) is bounded by
%
\begin{equation}
\label{son} \sum_{s=1}^{m}\sum
_{j=2s+1}^{2s+n}\frac{1}{j}\cdot
\frac{1}{\sqrt
{s}}=\sum_{s=1}^{m}
\frac{1}{\sqrt{s}} \sum_{j=2s+1}^{2s+n}
\frac{1}{j}
\end{equation}
for all $n\geq2$. Now,
\[
\sum_{j=2s+1}^{2s+n}\frac{1}{j} \leq
\sum_{j=2s+1}^{2s+n}\int_{j-1}^{j}
\frac{1}{x} \,dx =\int_{2s}^{2s+n}
\frac{1}{x} \,dx=\log \biggl(1+\frac{n}{2s}\biggr)
\]
for all $s\geq1$. This implies that (\ref{son}) is controlled by
\[
\sum_{s=1}^{m}\frac{1}{\sqrt{s}}\log
\biggl(1+\frac{n}{s}\biggr) \leq\sum_{s=1}^{m}
\int_{s-1}^s\frac{1}{\sqrt{y}}\log\biggl(1+
\frac{n}{y}\biggr) \,dy=\int_0^m
\frac{1}{\sqrt{y}}\log\biggl(1+\frac{n}{y}\biggr) \,dy.
\]
Set $u=y/n$. Then the last integral is equal to
\[
\int_0^{m/n}\frac{1}{\sqrt{nu}}\log\biggl(1+
\frac{1}{u}\biggr)\cdot n \,du \leq\sqrt{n}\int_0^{\infty}
\frac{1}{\sqrt{u}}\log\biggl(1+\frac
{1}{u}\biggr) \,du.
\]
Trivially, $\frac{1}{\sqrt{u}}\log(1+\frac{1}{u}) \sim\frac
{1}{u^{3/2}}$ as $u\to+\infty$ and $\frac{1}{\sqrt{u}}\log(1+\frac
{1}{u})\sim-\frac{\log u}{\sqrt{u}}$ as $u\to0^+$. It follows that
$0<\int_0^{\infty}\frac{1}{\sqrt{u}}\log(1+\frac{1}{u})
\,du<\infty$. Therefore, by (\ref{apple}),
\[
4m^2\sum_{\lambda=(r,s, 1^{m-r-s})}\frac{\theta_{(m)}^{\lambda
}(2)^2}{C_{\lambda}(2)}
\mathcal{N}_\lambda^2(n) \leq Cn
\]
for all $m\geq n\geq2$, where the sum is taken over all possible
$r\geq s \geq1$ and $3\leq l(\lambda) =m-r-s+2 \le n$. Combining
this, (\ref{Norton}), (\ref{fan}) and (\ref{rice}), we arrive at
\[
\mathbb{E} \bigl[ \bigl|p_{m}(Z_n)\bigr|^2 \bigr] \leq
K n
\]
for all $m\geq n \geq2$, where $K$ is a universal constant.
\end{pf*}

\subsection{Proof of Proposition \texorpdfstring{\protect\ref{search}}{2}}\label{black_yellow}
The following result allows us to express the variance for the circular
symplectic ensembles ($\beta=4$) in terms of some familiar quantities
treated earlier in the case of the circular orthogonal ensembles and a
new quantity $\mathcal{N}_\lambda^2(-2n)$.

\begin{lemma}[(Duality lemma)] \label{with} Recall (\ref{please}). For
any $m\geq1$ and $n\geq2$, the following holds:
%
\begin{equation}
\label{rabbit} \mathbb{E} \bigl[\bigl|p_m\bigl(Z_n^{1/2}
\bigr)\bigr|^2\bigr]= m^2 \sum_{\lambda\vdash m: \lambda_1 \le n}
\frac{\theta^\lambda_{(m)}(2)^2}{C_\lambda(2)} \mathcal{N}_\lambda^{2} (-2n),
\end{equation}
where $\lambda=(\lambda_1, \lambda_2, \ldots)$ and
%
\begin{equation}
\label{Great} \mathcal{N}_\lambda^2(-2n) =\prod
_{(i,j) \in\lambda} \biggl(1+ \frac{1}{2n-2j+i} \biggr).
\end{equation}
\end{lemma}

\begin{pf}
The quantity $\theta^\lambda_\mu(\alpha)$ has the following
duality (see (10.30) from \cite{Macdonald1998}): for partitions
$\lambda,\mu$ of $m$,
\[
\theta^\lambda_\mu(\alpha)= (-\alpha)^{m-l(\mu)}
\theta^{\lambda'}_\mu(1/\alpha),
\]
where $\lambda'$ is the partition of $m$ corresponding to the Young
diagram of the transpose of $\lambda$.
From (\ref{ski}), it is easy to see the duality
\[
C_\lambda(\alpha)= \prod_{(i,j) \in\lambda} \bigl(
\alpha(\lambda_i-j) +\lambda_j'-i+1\bigr)
\bigl(\alpha(\lambda_i-j) +\lambda _j'-i+
\alpha\bigr)=\alpha^{2m} C_{\lambda'}(1/\alpha).
\]
We furthermore have
\begin{eqnarray*}
\mathcal{N}_\lambda^\alpha(n)& =& \prod
_{(i,j) \in\lambda} \frac{n+(j-1)\alpha-(i-1)}{n+j\alpha
-i}
\\
&=& \prod_{(i,j) \in\lambda'} \frac{n+(i-1)\alpha-(j-1)}{n+i\alpha-j}
\\
&=&\prod_{(i,j) \in\lambda'} \frac{-n/\alpha-(i-1)+(j-1)/\alpha
}{-n/\alpha-i+j/\alpha}
\\
&=& \mathcal{N}_{\lambda'}^{1/\alpha} (-n/\alpha),
\end{eqnarray*}
where
%
\begin{equation}
\label{shining_star} \mathcal{N}_{\mu}^{\gamma} (x):=\prod
_{(i,j)\in\mu}\frac
{x-(i-1)+\gamma(j-1)}{x-i+\gamma j}
\end{equation}
for any partition $\mu,  \gamma>0$ and $x \in\mathbb{R}$
satisfying that the denominators in the product are not equal to zero.
It follows from dualities given above and (\ref{please}) that
\begin{eqnarray*}
\mathbb{E} \bigl[\bigl|p_m\bigl(Z_n^\alpha
\bigr)\bigr|^2\bigr]&=& \alpha^2 m^2 \sum
_{\lambda\vdash m
:l(\lambda) \le n} \frac{\alpha^{2m-2} \theta^{\lambda'}_{(m)}(1/\alpha)^2}{
\alpha^{2m} C_{\lambda'}(1/\alpha)} \mathcal{N}_{\lambda
'}^{1/\alpha}
(-n/\alpha)
\\
&=&m^2 \sum_{\lambda\vdash m: \lambda_1 \le n} \frac{\theta^\lambda_{(m)}(1/\alpha)^2}{C_\lambda(1/\alpha)}
\mathcal{N}_\lambda^{1/\alpha} (-n/\alpha),
\end{eqnarray*}
where $\lambda=(\lambda_1, \lambda_2, \ldots)$. Plugging $\alpha
=1/2$ into this identity,
\[
\mathbb{E} \bigl[\bigl|p_m\bigl(Z_n^{1/2}
\bigr)\bigr|^2\bigr]= m^2 \sum_{\lambda\vdash m: \lambda_1 \le n}
\frac{\theta^\lambda_{(m)}(2)^2}{C_\lambda(2)} \mathcal{N}_\lambda^{2} (-2n).
\]
Finally, from (\ref{shining_star}),
\[
\mathcal{N}_\lambda^2(-2n) = \prod
_{(i,j) \in\lambda} \frac{-2n +2j-i-1}{-2n +2j-i} =\prod
_{(i,j) \in\lambda} \biggl(1+ \frac{1}{2n-2j+i} \biggr).
\]
The proof is completed.
\end{pf}

\begin{lemma}\label{showup} Let $m\geq n \geq1$ and $\lambda
=(\lambda_1, \lambda_2, \ldots) \vdash m$ with $\lambda_1\leq n$.
Let $N_{\lambda}^2(-2n)$ be as in (\ref{Great}).
Then there exists a universal constant $K>0$ such that:
\begin{longlist}[(ii)]
\item[(i)] $N_{\lambda}^2(-2n) \leq K\sqrt{n}$ if $m=n$ and $\lambda=(n)$;
\item[(ii)] $N_{\lambda}^2(-2n) \leq K\frac{n}{\sqrt{(n-r+1)(n-s+1)}}$ if
$\lambda=(r, s)$ with $1\leq s\leq r \leq n$ and $r+ s = m$.
\end{longlist}
\end{lemma}

\begin{pf} Let $C:=\max_{\lambda, n\leq2}N_{\lambda}^2(-2n)$, where
$\lambda$ goes over all partitions as in~(i) and (ii) with $\lambda
_1\leq2$. Since $2n-2j+i\geq i\geq1$ for all $(i,j)\in\lambda$ with
$\lambda_1\leq n$, we know $N_{\lambda}^2(-2n)>1$, and hence $C>1$.
Also, since $m=r+s\leq2n\leq4$, these partitions are only of finitely
many. Thus, $1<C<\infty$. Then (i) and (ii) hold by taking $K=C$. From
now on, we assume, without loss of generality, that $n\geq3$.

(i) In this case,
\[
N_{\lambda}^2(-2n) = \prod_{j=1}^n
\biggl(1+\frac{1}{2n-2j+1} \biggr)=\prod_{k=1}^n
\biggl(1+\frac{1}{2k-1} \biggr) \leq\exp \Biggl(\sum
_{k=1}^n\frac{1}{2k-1} \Biggr).
\]
Now,
\[
\sum_{k=1}^n\frac{1}{2k-1} \leq1+
\sum_{k=2}^n\int_{k-1}^k
\frac
{1}{2x-1} \,dx=1+\int_1^n
\frac{1}{2x-1} \,dx=1+\frac{1}{2}\log(2n-1).
\]
The desired result then follows.

(ii) By the same argument as in the proof of (i),
%
\begin{eqnarray}
\label{project} \log N_{\lambda}^2(-2n) &\leq& \sum
_{j=1}^r\frac{1}{2n-2j+1} + \sum
_{j=1}^s\frac{1}{2n-2j+2}
\nonumber
\\[-8pt]
\\[-8pt]
\nonumber
& \leq& \sum_{k=n-r+1}^n\frac{1}{2k-1}
+ \sum_{k=n-s+1}^n\frac{1}{2k-1}.
\end{eqnarray}
Similar to (i),
\begin{eqnarray*}
\sum_{k=n-r+1}^n\frac{1}{2k-1}& \leq& 1+
\int_{n-r+1}^{n}\frac
{1}{2x-1} \,dx =  1+
\frac{1}{2}\log\frac{2n-1}{2(n-r)+1}
\\
& \leq& 1+ \frac{1}{2}\log\frac{2n}{n-r+1}.
\end{eqnarray*}
A similar inequality also holds true for the last sum in (\ref
{project}). Thus,
\[
\log N_{\lambda}^2(-2n) \leq C +\frac{1}{2}\log
\frac{n^2}{(n-r+1)(n-s+1)}.
\]
This implies (ii).
\end{pf}

\begin{lemma}\label{duckling} Let $\mathcal{N}_\lambda^2(-2n)$ be as
in (\ref{Great}). Let $\lambda=(r,s,1^{m-r-s})$ with
$1 \le s \le r \le n$, $m>r+s$ and $m\geq n$. Then there exists a
universal constant $K>0$ such that
\[
\frac{1}{K}\cdot\frac{m}{\sqrt{(n-r+1)(n-s+1)}}\leq\mathcal {N}_\lambda^2(-2n)
\leq K\cdot\frac{m}{\sqrt{(n-r+1)(n-s+1)}}.
\]
\end{lemma}

\begin{pf} We prove the upper bound and lower bound in two steps.

\textit{Step \textup{1:} Upper bound}. First,
%
\begin{eqnarray} \label{Paul}
\log N_{\lambda}^2(-2n) &= & \sum_{j=1}^r
\log \biggl(1+ \frac
{1}{2n-2j+1} \biggr) + \sum_{j=1}^s
\log \biggl(1+\frac{1}{2n-2j+2} \biggr)
\nonumber
\\[-8pt]
\\[-8pt]
\nonumber
& &{} + \sum_{i=3}^{m-r-s+2} \log \biggl(1+
\frac{1}{2n+i-2} \biggr)
\\
\label{hen}
& \leq& \sum_{j=1}^r
\frac{1}{2n-2j+1}
\nonumber
\\[-8pt]
\\[-8pt]
\nonumber
&&{}+ \sum_{j=1}^s
\frac
{1}{2n-2j+2} + \sum_{i=3}^{m-r-s+2}
\frac{1}{2n+i-2}
\end{eqnarray}
by the inequality $\log(1+x)\leq x$ for all $x>-1$.
Easily, if $r>1$ then
\begin{eqnarray*}
\sum_{j=1}^r \frac{1}{2n-2j+1} &\leq&
1+\sum_{j=1}^{r-1}\int_{j}^{j+1}
\frac{1}{2n-2x+1} \,dx
\\
&=& 1+ \int_1^{r}\frac{1}{2n-2x+1} \,dx
\\
& = & 1+\frac{1}{2}\log\frac{2n-1}{2n-2r+1}
\end{eqnarray*}
and this assertion is evidently true for $r=1$. So the above inequality
holds for all $r\geq1$.
Thus,
\[
\sum_{j=1}^s \frac{1}{2n-2j+2} \leq
\sum_{j=1}^s \frac{1}{2n-2j+1} \leq1 +
\frac{1}{2}\log\frac{2n-1}{2n-2s+1}.
\]
Likewise,
\begin{eqnarray*}
\sum_{i=3}^{m-r-s+2} \frac{1}{2n+i-2} &
\leq& \sum_{i=3}^{m-r-s+2}\int
_{i-1}^{i}\frac{1}{2n+x-2} \,dx
\\
& = & \int_{2}^{m-r-s+2}\frac{1}{2n+x-2} \,dx=\log
\frac{m+2n-r-s}{2n}.
\end{eqnarray*}
Combining the three inequalities with (\ref{hen}), we get
%
\begin{eqnarray}\label{noodles}
\log N_{\lambda}^2(-2n) &\leq& 2 + \frac{1}{2}\log
\frac
{(2n-1)^2(m+2n-r-s)^2}{(2n)^2(2n-2r+1)(2n-2s+1)} \nonumber
\\
& \leq& 2+\log3 + \frac{1}{2}\log\frac{m^2}{(2n-2r+1)(2n-2s+1)}
\\
& \leq& (2+\log3) + \log\frac{m}{\sqrt{(n-r+1)(n-s+1)}},
\nonumber
\end{eqnarray}
where the fact $\frac{2n-1}{2n} \leq1$ and the fact $m+2n-r-s \leq
m+2n\leq3m $ are used in the second inequality; the facts $2n-2r+1
\geq n-r+1$ and $2n-2s+1 \geq n-s+1$ are used in the last inequality.

\textit{Step \textup{2:} Lower bound}. Review (\ref{Paul}). Use the inequality
that $\log(1+x)\geq x-x^2$ for all $x\geq0$ to have
\begin{eqnarray*}
&&\log\mathcal{N}_\lambda^2(-2n) \\
&&\qquad\geq \sum
_{j=1}^r \frac{1}{2n-2j+1} + \sum
_{j=1}^s \frac
{1}{2n-2j+2} + \sum
_{i=3}^{m-r-s+2} \frac{1}{2n+i-2}
\\
& &\qquad\quad{} -\sum_{j=1}^r \frac{1}{(2n-2j+1)^2} -
\sum_{j=1}^s \frac
{1}{(2n-2j+2)^2} - \sum
_{i=3}^{m-r-s+2} \frac{1}{(2n+i-2)^2}.
\end{eqnarray*}
Observe that each term in the last three sums is strictly monotone in
its corresponding index. From the fact $\sum_{i=1}^{\infty}\frac
{1}{i^2}=\frac{\pi^2}{6}$, we know that the sum of the last three
sums is bounded by $\frac{\pi^2}{2}$. By the same arguments as before,
\begin{eqnarray*}
\sum_{j=1}^r \frac{1}{2n-2j+1} &\geq&
\sum_{j=1}^{r}\int_{j-1}^{j}
\frac{1}{2n-2x+1} \,dx
\\
&=& \int_0^{r}\frac{1}{2n-2x+1} \,dx
\\
& = & \frac{1}{2}\log\frac{2n+1}{2n+1-2r}.
\end{eqnarray*}
And
\[
\sum_{j=1}^s \frac{1}{2n-2j+2} \geq
\sum_{j=1}^s \frac
{1}{2(n+1)-2j+1} \geq
\frac{1}{2}\log\frac{2n+3}{2n+3-2s}.
\]
Now,
\begin{eqnarray*}
\sum_{i=3}^{m-r-s+2} \frac{1}{2n+i-2} &
\geq& \sum_{i=3}^{m-r-s+2}\int
_{i}^{i+1}\frac{1}{2n+x-2} \,dx
\\
& = & \int_{3}^{m-r-s+3}\frac{1}{2n+x-2} \,dx=\log
\frac{m+2n-r-s+1}{2n+1}.
\end{eqnarray*}
In summary,
\begin{eqnarray*}
\log\mathcal{N}_\lambda^2(-2n) &\geq& -\frac{\pi^2}{2} +
\frac{1}{2}\log\frac
{(2n+3)(m+2n-r-s+1)^2}{(2n+1)(2n-2r+1)(2n-2s+3)}
\\
& \geq& -\frac{\pi^2}{2} + \frac{1}{2}\log\frac
{m^2}{(2n-2r+1)(2n-2s+3)}
\\
& \geq& \biggl(-\frac{\pi^2}{2}-\frac{1}{2}\log6 \biggr) + \log
\frac
{m}{\sqrt{(n-r+1)(n-s+1)}},
\end{eqnarray*}
where we use the fact that $r+s\leq2n$ in the second inequality, and
the facts that $2n-2r+1 \leq2(n-r+1)$ and $2n-2s+3 \leq3(n-s+1)$ in
the last inequality.
\end{pf}

\begin{lemma}\label{bubbling} Let $\mathcal{N}_\lambda^2(-2n)$ be as
in (\ref{Great}). Then there exists a universal constant $K>0$ such
that
\begin{eqnarray*}
&&\phantom{i}\mathrm{(i)}\quad  N_{\lambda}^2(-2n) \leq K\sqrt{\frac{n}{n-m+1}}\qquad
\mbox {if } \lambda=(m) \mbox{ and } 1\leq m \leq n;
\\
&& \mathrm{(ii)}\quad N_{\lambda}^2(-2n) \leq K\frac{n}{\sqrt{(n-r+1)(n-s+1)}}
\end{eqnarray*}
if $\lambda=(r, s)$ with $1\leq s\leq r$ and $r+s = m \leq n$, or
$\lambda=(r,s,1^{m-r-s})$ with
$1 \le s \le r$ and $n\geq m>r+s$.
\end{lemma}

\begin{pf}
(i) Look at (i) in the proof of Lemma~\ref{showup}, replace ``$\prod_{j=1}^n$'' with ``$\prod_{j=1}^m$'' to have
\[
\log N_{\lambda}^2(-2n) \leq \sum_{j=1}^m
\frac{1}{2n-2j+1}=\sum_{k=n-m+1}^n
\frac{1}{2k-1} \leq1+\int_{n-m+1}^n
\frac{1}{2x-1} \,dx
\]
since $\frac{1}{2k-1} \leq\int_{k-1}^k\frac{1}{2x-1} \,dx$ for all
$k\geq2$. Thus,
\[
\log N_{\lambda}^2(-2n) \leq1+\frac{1}{2}\log
\frac{2n-1}{2n-2m+1} \leq \biggl(1+\frac{1}{2}\log2 \biggr)+
\frac{1}{2}\log\frac{n}{n-m+1}
\]
since $2n-1 < 2n$ and $2n-2m+1\geq n-m+1$.
This gives (i).

(ii) We consider the two aforementioned cases separately.

\textit{Case \textup{(a)}: $\lambda=(r, s)$ with $1\leq s\leq r$ and $r+s = m \leq n$}.
Review the proof of (ii) of Lemma~\ref{showup}.
The first paragraph is still true. The only occurrence of ``$m$'', which
is in ``$r+s=m$'', does not show up in the proof. So we obtain the same
inequality.

\textit{Case \textup{(b)}: $\lambda=(r,s,1^{m-r-s})$ with
$1 \le s \le r$ and $n\geq m>r+s$}. Review step 1 in the proof of Lemma~\ref{duckling}, no restriction on the relationship between $m$ and $n$
is used from the beginning to (\ref{noodles}). So, by (\ref
{noodles}), we have
\begin{eqnarray*}
\log N_{\lambda}^2(-2n) &\leq& 2 + \frac{1}{2}\log
\frac
{(2n-1)^2(m+2n-r-s)^2}{(2n)^2(2n-2r+1)(2n-2s+1)}
\\
& \leq& 2 + \frac{1}{2}\log\frac{9n^2}{(2n-2r+1)(2n-2s+1)}
\\
& \leq& (2+\log3) + \frac{1}{2}\log\frac{n^2}{(n-r+1)(n-s+1)}
\end{eqnarray*}
since $m+2n-r-s \leq3n$. This gives the conclusion.
\end{pf}

\begin{lemma}\label{feline} There exists a universal constant $K>0$
such that
\begin{eqnarray*}
&& \phantom{ii}\mathrm{(i)}\quad  \sqrt{\frac{mn}{n-m}}\leq Km \qquad\mbox{for all } 1\leq m <n;
\\
&& \phantom{i}\mathrm{(ii)}\quad  \log \biggl(1+\sqrt{\frac{m}{n-m}} \biggr) \leq K\sqrt{
\frac
{m}{n}} \log(m+1)\qquad \mbox{for all }1\leq m <n;
\\
&& \mathrm{(iii)}\quad  \sup_{m>n\geq1} \biggl\{\frac{m}{n}\cdot
\frac{1}{\sqrt
{m-n}}\tan^{-1}\sqrt{\frac{n}{m-n}} \biggr\} \leq K.
\end{eqnarray*}
\end{lemma}

\begin{pf} (i) If $1\leq m \leq\frac{n}{2}$, then $n-m \geq\frac
{n}{2}$ and hence $\sqrt{\frac{mn}{n-m}}\leq\sqrt{2}m$. If $\frac
{n}{2} \leq m \leq n-1$, then $\sqrt{\frac{mn}{n-m}}\leq\sqrt{mn}
\leq\sqrt{2}m$.

(ii) If $1\leq m \leq\frac{n}{2}$ then:
\[
\log \biggl(1+\sqrt{\frac{m}{n-m}} \biggr) \leq\sqrt{\frac{m}{n-m}}
\leq2\sqrt{\frac{m}{n}} \leq\frac{2}{\log2}\cdot\sqrt{
\frac
{m}{n}}\cdot\log(m+1).
\]
If $\frac{n}{2} < m < n$, then $\frac{m}{n-m}\geq1$ and $2\sqrt
{\frac{m}{n}}\geq1$. It follows that
\begin{eqnarray*}
\log \biggl(1+\sqrt{\frac{m}{n-m}} \biggr) &\leq& \log \biggl(2\sqrt{
\frac{m+1}{n-m}} \biggr)
\\
& \leq& \log2 + \frac{1}{2}\log(m+1)
\\
& \leq& 2\log(m+1)\leq4\sqrt{\frac{m}{n}}\log(m+1).
\end{eqnarray*}

(iii) Define
\[
A_{m,n}=\frac{m}{n}\cdot\frac{1}{\sqrt{m-n}}\tan^{-1}
\sqrt{\frac
{n}{m-n}}.
\]
Obviously,
\[
\sup_{m>n\geq 1}A_{m,n}\leq\sup_{n<m\leq5n}A_{m,n}
+ \sup_{m>
5n}A_{m,n} \leq\frac{5\pi}{2} + \sup
_{m> 5n}A_{m,n}.
\]
Note that
$\tan^{-1}x < x$ for all $x>0$. It follows that
\begin{eqnarray*}
\sup_{m> 5n}A_{m,n} &\leq& \sup_{m> 5n}
\biggl\{\frac{m}{n}\cdot \frac{1}{\sqrt{m-n}}\sqrt{\frac{n}{m-n}}
\biggr\}
\\
& \leq& \sup_{m> 5n} \biggl\{\frac{1}{\sqrt{n}}\cdot
\frac
{1}{1-{n}/{m}} \biggr\}\leq\frac{5}{4}.
\end{eqnarray*}
Then (iii) follows.
\end{pf}

\begin{lemma}\label{code} Recall (\ref{rabbit}). Let $m\geq n \geq
2$. Define
\[
E_{m,n}= m^2 \sum_{\lambda}
\frac{\theta^\lambda_{(m)}(2)^2}{C_\lambda(2)} \mathcal{N}_\lambda^{2} (-2n),
\]
where the sum is taken over all $\lambda=(r,s,1^{m-r-s})$ with
$1 \le s \le r \le n$ and $m>r+s$. Then there exists a universal
constant $K>0$ such that the following hold
\begin{longlist}[(iii)]
\item[(i)] $E_{m,n} \leq K \delta^{-1}n$ for all $m\geq(1+\delta
)n$ and $\delta\in(0,1]$.

\item[(ii)] $E_{m,n} \leq K n\log n$ for all $m\geq n\geq2$.

\item[(iii)] Let $w=m-n\geq0$. Then $E_{m,n} \geq K(w+1)^{-2} n\log
n$ for all $n\geq12$.
\end{longlist}
\end{lemma}

\begin{pf}
(i) From the first assertion of Lemmas \ref{cucumber} and~\ref
{duckling}, we know
\[
\label{thrower} E_{m,n}\leq Cm\sqrt{n}\sum
_{r,s}\frac{1}{\sqrt{(n-r+1)(n-s+1)s} (m-r+s)},
\]
where the sum runs over all possible $r$ and $s$ satisfying $1 \le s
\le r \le n$, \mbox{$m-r-s \geq 1$}. Obviously, $s\leq\frac{m}{2}$. Therefore,
%
\begin{equation}
\label{kindle} E_{m,n} \leq Cm\sqrt{n}\sum
_{s=1}^{m_n}\frac{1}{\sqrt
{(n-s+1)s}}\sum
_{r=s}^{n}\frac{1}{\sqrt{n-r+1} (m-r+s)},
\end{equation}
where $m_n=n \wedge[m/2]$.

\textit{Step} 1. First, we consider the term corresponding to $s=1$
dividing by $C$, which is equal to
%
\begin{equation}
\label{horse} m\sum_{r=1}^{n}
\frac{1}{\sqrt{(n-r+1)} (m-r+1)}:=V_{m,n}^1. 
\end{equation}
%
Easily $V_{n,n}^1=n\sum_{r=1}^n\frac{1}{(n-r+1)^{3/2}}<n\sum_{j=1}^{\infty}\frac{1}{j^{3/2}}=n\zeta(3/2)$, where $\zeta(z)$ is
the Riemann zeta function. Assume now $m>n\geq1$. Then
%
\begin{equation}
\label{viola} V_{m,n}^1 = m\sum
_{j=1}^n\frac{1}{\sqrt{j}}\cdot\frac{1}{m-n+j}
\end{equation}
by setting $j=n-r+1$. Each term in the sum is bounded by $\int_{j-1}^{j}\frac{1}{\sqrt{x}}\cdot\frac{1}{m-n+x} \,dx$. Consequently,
\begin{eqnarray*}
V_{m,n}^1 &\leq& m\int_0^n
\frac{1}{\sqrt{x}}\cdot\frac{1}{m-n+x} \,dx
\\
& = & \frac{2m}{\sqrt{m-n}}\int_0^{\sqrt{n/(m-n)}}
\frac{1}{1+y^2} \,dy
\\
&=&\frac{2m}{\sqrt{m-n}}\tan^{-1}\sqrt{\frac{n}{m-n}}
\end{eqnarray*}
by defining $y=\sqrt{\frac{x}{m-n}}$. From (iii) of Lemma~\ref
{feline}, we obtain that
%
\begin{equation}
V_{m,n}^1 \leq (2n)\cdot\sup_{m>n\geq 1}
\biggl\{\frac{m}{n}\cdot\frac
{1}{\sqrt{m-n}}\tan^{-1}\sqrt{
\frac{n}{m-n}} \biggr\} \leq Cn \label{Portland}
\end{equation}
for any $m>n\geq1$.
Hence, to prove the conclusion, it suffices to show
%
\begin{equation}\qquad
\label{fund} W_{m,n}:=\sum_{s=2}^{m_n}
\frac{1}{\sqrt{(n-s+1)s}}\sum_{r=s}^{n}
\frac{1}{\sqrt{n-r+1} (m-r+s)} \leq C\delta^{-1}\frac
{\sqrt{n}}{m}.
\end{equation}

\textit{Step} 2. In this step, we prove (\ref{fund}) holds for all $m\geq
(1+\delta)n$.
Set $j=n-r+1$. Then, using the same argument as in estimating the term
in (\ref{viola}), we have
\begin{eqnarray*}
\sum_{r=s}^{n}\frac{1}{\sqrt{(n-r+1)} (m-r+s)} &= &
\sum_{j=1}^{n-s+1}\frac{1}{(m-n+s-1)+j}\cdot
\frac{1}{\sqrt{j}}
\\
& \leq& \sum_{j=1}^{n-s+1}\int
_{j-1}^{j}\frac{1}{a+x}\cdot
\frac
{1}{\sqrt{x}} \,dx
\\
& = & \int_0^{n-s+1}\frac{1}{a+x}\cdot
\frac{1}{\sqrt{x}} \,dx,
\end{eqnarray*}
where $a=m-n+s-1\geq1$ for $s\geq2$. Let $y=\sqrt{x/a}$. It follows that
\begin{eqnarray*}
\int_0^{n-s+1}\frac{1}{a+x}\cdot
\frac{1}{\sqrt{x}} \,dx &=& \frac
{2}{\sqrt{a}}\int_0^{\sqrt{(n-s+1)/a}}
\frac{1}{1+y^2} \,dy
\\
& = & \frac{2}{\sqrt{a}}\cdot\tan^{-1}\sqrt{\frac{n-s+1}{a}}.
\end{eqnarray*}
Thus,
%
\begin{equation}
\label{bird} \sum_{r=s}^{n}
\frac{1}{\sqrt{(n-r+1)} (m-r+s)} \leq\frac{2}{\sqrt
{a}}\cdot\tan^{-1}\sqrt{
\frac{n-s+1}{a}}
\end{equation}
for any $m\geq n \geq1$ and $s\geq2$ [we do not need the condition
``$m\geq(1+\delta)n$'' here].
Therefore, for any $n\geq2$,
\begin{eqnarray*}
W_{m,n} &\leq& 2\sum_{s=2}^{n}
\frac{1}{\sqrt
{s(n-s+1)(m-n+s-1)}}\cdot\tan^{-1}\sqrt{\frac{n-s+1}{m-n+s-1}}
\\
& \leq& 2\sum_{k=1}^{n-1}
\frac{1}{\sqrt{k(m-k)(n-k)}}\cdot\tan ^{-1}\sqrt{\frac{k}{m-k}}
\end{eqnarray*}
by letting $k=n-s+1$. The term corresponding to $k=n-1$ in the sum is
equal to
\[
\label{horse1} \frac{1}{\sqrt{(n-1)(m-n+1)}}\tan^{-1}\sqrt{\frac{n-1}{m-n+1}}.
\]
By the inequality $\tan^{-1} x <x$ for all $x>0$, it is seen that the
above quantity is controlled by
$\frac{1}{m-n+1} \leq\frac{2}{\delta}\cdot\frac{\sqrt{n}}{m}$
due to the fact $\frac{m}{m-n}=(1-\frac{n}{m})^{-1}\leq1+\delta
^{-1}\leq2\delta^{-1}$ from the assumption $m\geq(1+\delta)n$.
Consequently, to prove (\ref{fund}), it suffices to show
%
\begin{equation}
\label{fragrant} U_{m,n}:= \sum_{k=1}^{n-2}
\frac{1}{\sqrt{k(m-k)(n-k)}}\cdot\tan ^{-1}\sqrt{\frac{k}{m-k}} \leq C
\delta^{-1}\frac{\sqrt{n}}{m}
\end{equation}
for all $m\geq(1+\delta)n$ and $n\geq3$. In fact, since $\tan^{-1}
x <x$ for all $x>0$,
\begin{eqnarray*}
U_{m,n} &\leq& \sum_{k=1}^{n-2}
\frac{1}{(m-k)\sqrt{n-k}}
\\
& \leq& \sum_{k=1}^{n-2}\int
_k^{k+1}\frac{1}{(m-x)\sqrt{n-x}} \,dx
\\
&= & \int_1^{n-1}\frac{1}{(m-x)\sqrt{n-x}} \,dx
\end{eqnarray*}
by the obvious monotonicity. Now,
\begin{eqnarray*}
U_{m,n} & \leq& \frac{1}{m-n+1}\int_1^{n-1}
\frac{1}{\sqrt{n-x}} \,dx
\\
& = & \frac{2\sqrt{n-1}-2}{m-n+1}
\\
& \leq& \frac{2\sqrt{n}}{m-n}.
\end{eqnarray*}
By the inequality $\frac{m}{m-n}\leq2\delta^{-1}$ again,
$U_{m,n}\leq4\delta^{-1}\frac{\sqrt{n}}{m}$. We get (\ref{fragrant}).

(ii) By taking $\delta=\frac{1}{2}$ in (i), we know $E_{m,n}\leq C n$
for $m\geq\frac{3}{2}n$.
So, to prove (ii),
we assume, without loss of generality, that $n\leq m \leq\frac{3}{2}n$.
Recall (\ref{kindle}). We know $m_n=n \wedge[m/2] \leq\frac
{3}{4}n$. Then $n-s+1\geq\frac{n}{4}$ for $1\leq s \leq m_n$. 
It follows that
\begin{eqnarray*}
E_{m,n} &\leq&C V_{m,n}^1+ Cm\sum
_{s=2}^{m_n}\frac{1}{\sqrt{s}}\sum
_{r=s}^n\frac{1}{\sqrt{n-r+1}  (m-r+s)}
\\
& \leq& Cn + C n\sum_{s=2}^n
\frac{1}{\sqrt{s}}\sum_{r=s}^{n}
\frac
{1}{\sqrt{n-r+1} (n-r+s)}
\end{eqnarray*}
by (\ref{horse}) and (\ref{Portland}) since $n\leq m\leq\frac
{3}{2}n$. Thus, to complete the proof, we only need to show
%
\begin{equation}
\label{beatles} H_{n}:=\sum_{s=2}^{n}
\frac{1}{\sqrt{s}}\sum_{r=s}^{n}
\frac
{1}{\sqrt{n-r+1} (n-r+s)} \leq C \log n
\end{equation}
for all $n\geq2$. In fact, apply (\ref{bird}) to the case $m=n$ so
that $a=s-1\geq\frac{s}{2}$ for $s\geq2$. We know that
\[
\sum_{r=s}^{n}\frac{1}{\sqrt{n-r+1} (n-r+s)}\leq
\frac{\pi}{\sqrt
{a}}\leq\frac{2\pi}{\sqrt{s}}
\]
for $s\geq2$. It follows that
\[
H_n \leq(2\pi)\cdot\sum_{s=2}^{n}
\frac{1}{s}\leq(2\pi)\cdot \biggl(1+\frac{\log n}{\log2} \biggr)\leq C\log
n
\]
by (\ref{shadow}) where $C=(4\pi)(\log2)^{-1}$. This gives (\ref{beatles}).

(iii) From Lemmas \ref{cucumber} and \ref{duckling},
\[
E_{m,n} \geq Cm\sum_{r, s}
\frac{m-r-s}{(m-r+s)^2}\cdot\sqrt{\frac
{r}{s}}\cdot\frac{1}{\sqrt{(n-r+1)(n-s+1)}},
\]
where $2s \leq r \leq n$ and $m>r+s$. Since $n-s+1\leq n$ and $\sqrt
{\frac{r}{s}}\geq\sqrt{\frac{n}{2s}}$ if $r\geq\frac{n}{2}$. Then
\begin{eqnarray*}
E_{m,n} &\geq& Cn\sum_{(r, s)\in T_1}
\frac{m-r-s}{(m-r+s)^2}\cdot \frac{1}{\sqrt{s}}\cdot\frac{1}{\sqrt{n-r+1}}
\\
& = & Cn\sum_{(s,t)\in T_2}\frac{w+t-s}{(w+t+s)^2}\cdot
\frac
{1}{\sqrt{s}}\cdot\frac{1}{\sqrt{t+1}},
\end{eqnarray*}
where $w=m-n$ as defined in the statement of the lemma,
\[
T_1= \biggl\{(r,s)\in\mathbb{N}^2; 2s \leq r \le n,
m>r+s \mbox { and } r\geq\frac{n}{2} \biggr\},
\]
$t=n-r$ and
\[
T_2= \biggl\{(s, t)\in\mathbb{Z}^2; s\geq1, t\geq0, 2s
+t \le n, t\leq\frac{n}{2}\mbox{ and } w+t-s\geq1 \biggr\},
\]
where $\mathbb{N}$ is the set of positive integers and $\mathbb{Z}$
is the set of real integers.
Easily,
\[
T_2 \supset T_3:= \biggl\{(s, t)\in\mathbb{N}^2;
1\leq s\leq\frac
{t}{2} \mbox{ and } 2\leq t \leq\frac{n}{2} \biggr
\}.
\]
Consequently,
\begin{eqnarray*}
E_{m,n} &\geq& Cn\sum_{2\leq t \leq{n}/{2}}\sum
_{1\leq s\leq
{t}/{2}}\frac{w+t-s}{(w+t+s)^2}\cdot\frac{1}{\sqrt{s}}\cdot
\frac
{1}{\sqrt{t+1}}
\\
& \geq& \frac{Cn}{(1+w)^2}\sum_{2\leq t \leq{n}/{2}}\sum
_{1\leq s\leq{t}/{2}}\frac{t-s}{(t+s)^2}\cdot\frac{1}{\sqrt
{s}}\cdot
\frac{1}{\sqrt{t}}
\end{eqnarray*}
since $w+t+s\leq(w+1)(t+s)$ and $t+1\leq2t$. Note that $\frac
{t-s}{(t+s)^2}$ is strictly decreasing in $s\in[1, \frac{t}{2}]$, it
is bounded below by $\frac{2}{9}\cdot\frac{1}{t}$ for all $1\leq
s\leq\frac{t}{2}$. Thus,
\begin{eqnarray*}
E_{m,n} & \geq& \frac{Cn}{(w+1)^2}\sum_{2\leq t \leq{n}/{2}}
\sum_{1\leq s\leq{t}/{2}}\frac{1}{t^{3/2}}\cdot\frac{1}{\sqrt
{s}}
\\
& \geq& \frac{Cn}{(w+1)^2}\sum_{2\leq t \leq{n}/{2}}
\frac{1}{t}
\end{eqnarray*}
because $\frac{1}{t^{3/2}}\cdot\frac{1}{\sqrt{s}}\geq\frac
{1}{t^2}$ for $1\leq s \leq t$. Finally, by (\ref{shadow}),
\[
\sum_{2\leq t \leq{n}/{2}}\frac{1}{t} \geq\log \biggl(
\frac
{1}{2} \biggl[\frac{n}{2} \biggr] \biggr) \geq C\log n
\]
for $n \geq12$, where $C=\inf_{n\geq12} \{(\log n)^{-1}\log
(\frac{1}{2} [\frac{n}{2} ] ) \} \in(0, \infty)$.
In summary,
\[
E_{m,n} \geq C\cdot\frac{n\log n}{(w+1)^2}
\]
for all $n \geq12$.
\end{pf}

\begin{lemma}\label{infinitely} Recall (\ref{rabbit}). There exists a
universal constant $K>0$ such that
\[
\mathbb{E} \bigl[\bigl|p_m\bigl(Z_n^{1/2}
\bigr)\bigr|^2\bigr] \leq K m \log(m+1)
\]
for all $1\leq m< n$.
\end{lemma}

\begin{pf}
By Lemma~\ref{with},
%
\begin{equation}
\label{draft} \mathbb{E} \bigl[\bigl|p_m\bigl(Z_n^{1/2}
\bigr)\bigr|^2\bigr]= m^2 \sum_{\lambda\vdash m: \lambda_1 \le n}
\frac{\theta^\lambda_{(m)}(2)^2}{C_\lambda(2)} \mathcal{N}_\lambda^{2} (-2n).
\end{equation}
Since $m<n$, the restriction $\lambda_1 \le n$ automatically holds.
Review (\ref{romantic}). Many of the terms in the sum are equal to
zero except the following three types of partitions: (i)~$\lambda
=(m)$; (ii) $\lambda=(m-r, r)$ with $1\leq r \leq\frac{m}{2}$;
(iii) $\lambda=(r,s,1^{m-r-s})$ with $1 \le s \le r$ and $m-r-s\geq1$.

Now let us analyze the three sums separately.

\textit{Step \textup{1:} Analysis of the sum corresponding to case \textup{(i)}}. By (\ref
{crack}), (i)~of Lemma~\ref{bubbling} and (i) of Lemma~\ref{feline},
%
\begin{equation}
m^2 \sum_{\lambda=(m)} \frac{\theta^\lambda_{(m)}(2)^2}{C_\lambda(2)}
\mathcal{N}_\lambda^{2} (-2n) \leq C\sqrt{\frac{mn}{n-m+1}}
\leq C'm,\label{Haar}
\end{equation}
where both $C$ and $C'$ are universal constants.

\textit{Step \textup{2:} Analysis of the sum corresponding to case \textup{(ii)}}. Review
(\ref{Jerry}). Replace ``$(r,s)$'' in (ii) of Lemma~\ref{bubbling} by
``$(m-r,r)$'' to obtain
\begin{eqnarray*}
& & m^2 \sum_{\lambda=(m-r, r)} \frac{\theta^\lambda_{(m)}(2)^2}{C_\lambda(2)}
\mathcal{N}_\lambda^{2} (-2n)
\nonumber
\\[-8pt]
\\[-8pt]
\nonumber
&&\qquad \leq Cn\sqrt{m}\sum_{\lambda=(m-r, r)}\frac{1}{r^{5/2}}
\cdot \frac{1}{\sqrt{(n-r+1)(n-m+r+1)}},
\end{eqnarray*}
where the sum runs over all possible $\lambda=(m-r,r)$ with $1\leq
r\leq\frac{m}{2}$ and $m < n$. Use the trivial estimate $n-r+1\geq
\frac{n}{2}$ and $n-m+r+1\geq n-m$ to see that
%
\begin{equation}\qquad
\label{blog} m^2 \sum_{\lambda=(m-r, r)}
\frac{\theta^\lambda_{(m)}(2)^2}{C_\lambda(2)} \mathcal{N}_\lambda^{2} (-2n) \leq C\sqrt{
\frac{mn}{n-m}}\sum_{r=1}^{\infty}
\frac{1}{r^{5/2}}=C\cdot\zeta \biggl(\frac{5}{2} \biggr)\cdot m
\end{equation}
by (i) of Lemma~\ref{feline}, where $\zeta(z)$ is the Riemann zeta function.

\textit{Step \textup{3:} Analysis of the sum corresponding to case \textup{(iii)}}. Consider
\[
E_{m,n}:=m^2 \sum_{\lambda=(r,s,1^{m-r-s})}
\frac{\theta^\lambda_{(m)}(2)^2}{C_\lambda(2)} \mathcal{N}_\lambda^{2} (-2n),
\]
where the sum is taken over all partition $\lambda=(r,s,1^{m-r-s})$ with
$1 \le s \le r$ and $n>m>r+s$. From the first assertion of Lemmas~\ref
{cucumber} and~\ref{bubbling}, we know
\[
\label{thrower1} E_{m,n}\leq Cn\sqrt{m}\sum
_{r,s}\frac{1}{\sqrt{(n-r+1)(n-s+1)s} (m-r+s)},
\]
where the sum runs over all possible $r$ and $s$ satisfying $1 \le s
\le r$, $m-r-s \geq 1$ and $n> m$. Clearly, $s\leq\frac{m}{2}$,
hence $n-s+1\geq\frac{n}{2}$. Further, the restriction ``$m-r-s \geq
1$'' implies that $m\geq3$. Therefore,
%
\begin{equation}
\label{star_war} E_{m,n} \leq C\sqrt{m n}\cdot\sum
_{r=1}^{m-2}\frac{1}{\sqrt
{n-r}}\sum
_{s=1}^{r}\frac{1}{\sqrt{s} (m-r+s)}.
\end{equation}
Use the inequality $\frac{1}{\sqrt{s} (m-r+s)} \leq\int_{s-1}^s\frac{1}{\sqrt{x} (m-r+x)} \,dx$ to get
\begin{eqnarray*}
\sum_{s=1}^{r}\frac{1}{\sqrt{s} (m-r+s)} &\leq&
\int_0^{r}\frac
{1}{\sqrt{x} (m-r+x)} \,dx
\\
& = & \frac{2}{\sqrt{m-r}}\int_0^{\sqrt{r/(m-r)}}
\frac{1}{1+y^2} \,dy
= \frac{2}{\sqrt{m-r}}\tan^{-1}\sqrt{\frac{r}{m-r}}
\end{eqnarray*}
by setting $y=\sqrt{\frac{x}{m-r}}$. Since $\tan^{-1}x \leq\min\{
x, \frac{\pi}{2}\}$ for all $x>0$, we have
%
\begin{eqnarray}
\label{rainbow} & & \sum_{r=1}^{m-2}
\frac{1}{\sqrt{n-r}}\sum_{s=1}^{r}
\frac
{1}{\sqrt{s} (m-r+s)}
\nonumber
\\[-8pt]
\\[-8pt]
\nonumber
&& \qquad\leq 2\sum_{1\leq r\leq{m}/{2}}\frac{\sqrt{r}}{\sqrt
{n-r}  (m-r)} + \pi\sum
_{{m}/{2}\leq r \leq m-2}\frac{1}{\sqrt
{(n-r)(m-r)}}.
\end{eqnarray}
Observe that $n-r\geq\frac{n}{2}$ and $m-r \geq\frac{m}{2}$ for
$1\leq r \leq\frac{m}{2}$. Then
%
\begin{eqnarray}
\label{strange} \sum_{1\leq r\leq{m}/{2}}\frac{\sqrt{r}}{\sqrt{n-r}  (m-r)} &\leq&
\frac{4}{m\sqrt{n}}\sum_{1\leq r\leq{m}/{2}}\sqrt{r}
\nonumber
\\
& \leq& \frac{4}{m\sqrt{n}}\sum_{1\leq r\leq{m}/{2}}\int
_r^{r+1}\sqrt{x} \,dx
\\
& \leq& \frac{4}{m\sqrt{n}}\int_1^m\sqrt{x}
\,dx \leq3\sqrt{\frac
{m}{n}}\nonumber
\end{eqnarray}
since $\int_1^m\sqrt{x} \,dx=\frac{2}{3}(m^{3/2}-1)$. On the other hand,
\begin{eqnarray*}
\sum_{{m}/{2}\leq r \leq m-2}\frac{1}{\sqrt{(n-r)(m-r)}} & \leq& \sum
_{{m}/{2}\leq r \leq m-2}\int_r^{r+1}
\frac
{1}{\sqrt{(n-x)(m-x)}} \,dx
\\
& \leq& \int_{{m}/{2}}^{m-1}\frac{1}{\sqrt{(n-x)(m-x)}} \,dx
\\
& = & \int_1^{m/2}\frac{1}{\sqrt{y(n-m+y)}} \,dy
\end{eqnarray*}
by taking $y=m-x$. Now, let $u=\sqrt{y/(n-m)}$, the above integral becomes
\begin{eqnarray*}
2\int_{1/\sqrt{n-m}}^{\sqrt{m/(2(n-m))}}\frac{1}{\sqrt{1+u^2}} \,du & \leq& 4
\int_0^{\sqrt{m/(n-m)}}\frac{1}{u+1} \,du
\\
& = & 4\log \biggl(1+\sqrt{\frac{m}{n-m}} \biggr)
\end{eqnarray*}
by using the inequality $1+u^2 \geq\frac{1}{2}(1+u)^2$. By (ii) of
Lemma~\ref{feline},
\[
\sum_{{m}/{2}\leq r \leq m-2}\frac{1}{\sqrt{(n-r)(m-r)}} \leq C\sqrt{
\frac{m}{n}}\log(m+1),
\]
which together with (\ref{rainbow}) and (\ref{strange}) gives
\[
\sum_{r=1}^{m-2}\frac{1}{\sqrt{n-r}}\sum
_{s=1}^{r}\frac{1}{\sqrt
{s} (m-r+s)} \leq C\cdot
\biggl(\sqrt{\frac{m}{n}}+\sqrt{\frac
{m}{n}}\log(m+1) \biggr).
\]
This inequality and (\ref{star_war}) conclude that
%
\begin{equation}
\label{translate} E_{m,n}\leq C\cdot m \log(m+1).
\end{equation}
At last, according to (\ref{draft}) and its following paragraph, the
desired result follows by considering (\ref{Haar}), ({\ref{blog}})
and ({\ref{translate}}) together.
\end{pf}

\begin{pf*}{Proof of Proposition~\ref{search}}
From Lemma~\ref
{infinitely}, we know that we only need to prove the theorem for the
case $m\geq n$. By (\ref{rabbit}),
%
\begin{equation}
\label{why} \mathbb{E} \bigl[\bigl|p_m\bigl(Z_n^{1/2}
\bigr)\bigr|^2\bigr]=m^2 \sum_{\lambda}
\frac{\theta^\lambda_{(m)}(2)^2}{C_\lambda(2)} \mathcal{N}_\lambda^{2} (-2n),
\end{equation}
where the sum is taken over all $\lambda\vdash m\dvtx \lambda_1\leq n$.
Review (\ref{romantic}). Many of the terms in the sum are equal to
zero except the following three types of partitions: (i) $\lambda
=(m)$; (ii) $\lambda=(m-r, r)$ with $1\leq m-r \leq n$ and $1 \leq
r\leq\frac{m}{2}$; (iii) $\lambda=(r,s,1^{m-r-s})$ with $1 \le s \le
r \le n$ and $m-r-s\geq1$.

Now let us analyze the three cases one by one.

(a): \textit{The estimate of the sum corresponding to case \textup{(i)}}. When
$\lambda=(m)$ with $\lambda_1=m\leq n$, it is seen that $m=n$, then
from (\ref{crack}),
\[
m^2 \frac{\theta^\lambda_{(m)}(2)^2}{C_\lambda(2)}=m^2\frac{\theta
^{(m)}_{(m)}(2)^2}{C_{(m)}(2)} \leq C
\sqrt{n}.
\]
By (i) of Lemma~\ref{showup}, we know
%
\begin{equation}
\label{enter} m^2 \frac{\theta^\lambda_{(m)}(2)^2}{C_\lambda(2)}N_{\lambda
}^2(-2n)
\leq Cn.
\end{equation}

(b): \textit{The estimate of the sum corresponding to case} (ii). If
$\lambda=(m-r, r)$ with $1\leq m-r \leq n$ and $1 \leq r\leq\frac
{m}{2}$, then from (\ref{Jerry}) and (ii) of
Lemma~\ref{showup} [replace ``$(r, s)$'' by ``$(m-r, r)$''],
\begin{eqnarray*}
m^2 \frac{\theta^\lambda_{(m)}(2)^2}{C_\lambda(2)}N_{\lambda
}^2(-2n) & \leq&
\frac{C}{r^{5/2}}\sqrt{m}\cdot\frac{n}{\sqrt
{(n-r+1)(n-m+r+1)}}
\\
& \leq& C\cdot\frac{n^{3/2}}{r^{5/2}}\cdot\frac{1}{\sqrt{n-r+1}}
\end{eqnarray*}
since the restrictions on $r$ imply that $m\leq 2n$ and $1\leq r \leq
n$. Therefore,
\begin{eqnarray*}
& & m^2 \sum_{\lambda}\frac{\theta^\lambda_{(m)}(2)^2}{C_\lambda
(2)}N_{\lambda}^2(-2n)
\\
&&\qquad \leq Cn^{3/2}\sum_{1\leq r \leq{n}/{2}}
\frac
{1}{r^{5/2}}\cdot\frac{1}{\sqrt{n-r+1}} + Cn^{3/2}\sum
_{
{n}/{2}\leq r \leq n}\frac{1}{r^{5/2}}\cdot\frac{1}{\sqrt{n-r+1}},
\end{eqnarray*}
where the sum is taken over all $\lambda=(m-r, r)$ with $1\leq m-r
\leq n$ and $1 \leq r\leq\frac{m}{2}$. The term in the first sum is
controlled by $\frac{2}{r^{5/2}\sqrt{n}}$; each term in the second
sum is dominated by $\frac{8}{n^{5/2}}$. Consequently,
%
\begin{eqnarray}
\label{cut}&& m^2 \sum_{\lambda}
\frac{\theta^\lambda_{(m)}(2)^2}{C_\lambda
(2)}N_{\lambda}^2(-2n)
\nonumber
\\[-8pt]
\\[-8pt]
\nonumber
&&\qquad\leq C\cdot\biggl(2\zeta
\biggl(\frac{5}{2}\biggr)n + 8\biggr)\leq C\cdot \biggl(2\zeta\biggl(
\frac{5}{2}\biggr) + 8 \biggr)n,
\end{eqnarray}
where the sum is taken corresponding to case (ii) and $\zeta(z)$ is
the Riemann zeta function.

(c): \textit{The estimate of the sum corresponding to case} (iii). Let
$m\geq n \geq2$. Define
\[
E_{m,n}= m^2 \sum_{\lambda}
\frac{\theta^\lambda_{(m)}(2)^2}{C_\lambda(2)} \mathcal{N}_\lambda^{2} (-2n),
\]
where the sum is taken over all $\lambda=(r,s,1^{m-r-s})$ with
$1 \le s \le r \le n$ and \mbox{$m>r+s$}. By Lemma~\ref{code}, there exists a
universal constant $K>0$ such that the following hold:
\begin{longlist}[(B)$'$]
\item[(A)] $E_{m,n} \leq K \delta^{-1}n$ for all $m\geq(1+\delta
)n$ and $\delta\in(0,1]$.

\item[(B)] $E_{m,n} \leq K n\log n$ for all $m\geq n\geq2$.

\item[(C)] Let $w=m-n\geq0$. Then $E_{m,n} \geq K(w+1)^{-2} n\log
n$ for all $n\geq12$.

If $12\leq n\leq m \leq2n$ then $ n\log n \geq\frac{1}{4}m\log m$.
It follows that $E_{m,n} \geq K(w+1)^{-2} m\log m$ for all $12\leq
n\leq m \leq2n$. These combined with (\ref{why}), (\ref{enter})
and~(\ref{cut}) imply that

\item[(A)$'$] $\mathbb{E} [|p_m(Z_n^{1/2})|^2] \leq K \delta
^{-1}n$ for all $m\geq(1+\delta)n$ and $\delta\in(0,1]$.

\item[(B)$'$] $\mathbb{E} [|p_m(Z_n^{1/2})|^2] \leq K m\log m$ for
all $m\geq n\geq2$.

\item[(C)$'$] $\mathbb{E} [|p_m(Z_n^{1/2})|^2] \geq E_{m,n} \geq
K(w+1)^{-2} m\log m$ for all $12\leq n\leq m \leq2n$.
\end{longlist}

Finally, (A)$'$ and (C)$'$ are identical to (i) and
(iii) in the statement of Proposition~\ref{search}, respectively. As
mentioned at the beginning of the proof, (B)$'$ and Lemma~\ref
{infinitely} implies (ii) of the proposition.
\end{pf*}

\subsection{Proofs of Theorems \texorpdfstring{\protect\ref{Nevada}}{2}, \texorpdfstring{\protect\ref
{circular_case}}{3} and \texorpdfstring{\protect\ref{symplectic_case}}{4}}\label{green_brown}

With the preparations in Sections~\ref{red-blue} and~\ref
{black_yellow}, we are now ready to prove the central limit theorems.

\begin{pf*}{Proof of Theorem~\ref{Nevada}}
For any complex numbers
$c_k$'s and $d_k$'s with $\sum_{k=1}^m(|c_k| + |d_k|)\ne0$, define
\[
X_n=\sum_{j=1}^m \bigl[
c_jp_j\bigl(Z_n^{\alpha}
\bigr)+d_j\overline{p_j\bigl(Z_n^{\alpha
}
\bigr)} \bigr] \quad\mbox{and}\quad X=\sum_{j=1}^m
[c_j\xi_j+d_j\overline{\xi _j} ].
\]
We claim that, to prove the theorem, it is enough to show
%
\begin{equation}
\label{sales} \lim_{n\to\infty}\mathbb{E}\bigl(X_n^p
\bar{X}_n^q\bigr)=\mathbb{E}\bigl(X^p
\bar{X}^q\bigr) 
\end{equation}
for any integers $p\geq0$ and $q\geq0$ with $p+q\geq1$. In fact, for
a complex random vector $U=(U_1, \ldots, U_m)\in\mathbb{C}^m$, we
treat it as the real vector $\tilde{U}\in\mathbb{R}^{2m}$ by listing
their real and imaginary parts in a column. Since the real and the
complex parts of $U_j$ are $\frac{U_j + \bar{U}_j}{2}$ and $\frac
{U_j - \bar{U}_j}{2i}$, respectively, for each $j$, then $a'\tilde
{U}$ for $a\in\mathbb{R}^{2m}$ is a linear combination of $U_j$'s and
$\bar{U}_j$'s with complex coefficients. Thus, by the Cram\'{e}r--Wold
device (see, e.g., page~176 from \cite{Durrett}), to prove the theorem,
it suffices to show $X_n$ converges weakly to $X$ as $n\to\infty$.
Trivially, $X$ has the same distribution as that of $a\eta_1+b\eta_2$
where $\eta_1, \eta_2$ are i.i.d. with distribution $N(0,1)$ and $a,
b$ are complex numbers, hence $X$ is uniquely determined by its
moments. By the moment method, we only need to check (\ref{sales}).

First, $\mathbb{E}[p_{\mu}(Z_n^{\alpha}) \overline{p_{\nu
}(Z_n^{\alpha})}]=0$
unless the weights $|\mu|$ and $|\nu|$ are equal.
This fact follows in a way similar to
the proof of Proposition~\ref{exact}. The key is Lemma~\ref{nice}:
Jack polynomials are orthogonal.
We have
%
\begin{equation}
\label{America} \mathbb{E}\bigl[p_{\mu}\bigl(Z_n^{\alpha}
\bigr)\bigr]=0 \qquad\mbox{in particular } \mathbb{E}\bigl[p_k
\bigl(Z_n^{\alpha}\bigr)\bigr]=0 
\end{equation}
for all $|\mu|\geq1$ and $k\geq1$.

Second, expand $X_n^p\bar{X}_n^q$ and $X^p\bar{X}^q$ as sums of $M$
terms, 
where the number $M$ does not depend on $n$.
In the same way, it is seen that, to prove (\ref{sales}), we only need
to show
%
\begin{equation}
\label{bell} \lim_{n\to\infty}\mathbb{E} \Biggl(\prod
_{j=1}^mp_j\bigl(Z_n^{\alpha
}
\bigr)^{l_j}\cdot\prod_{j=1}^m
\overline{p_j\bigl(Z_n^{\alpha}\bigr)}^{ l_j'}
\Biggr) = \mathbb{E} \Biggl(\prod_{j=1}^m
\xi_j^{l_j}\cdot\prod_{j=1}^m
\bar {\xi}_j^{ l_j'} \Biggr)
\end{equation}
for nonnegative integers $l_j$'s and $l_j'$'s with $\sum_{j=1}^ml_j\geq1$ or $\sum_{j=1}^ml_j'\geq1$. Set $\mu=(1^{l_1},
2^{l_2}, \ldots, m^{l_m})$ and $\nu=(1^{l_1'}, 2^{l_2'}, \ldots,
m^{l_m'})$. Then,
according to (\ref{insect}),
\[
l(\mu)=\sum_{j=1}^ml_j
\quad\mbox{and}\quad z_{\mu}=\prod_{j=1}^sj^{l_j}l_j!.
\]
The quantities $l(\nu)$ and $z_{\nu}$ are defined similarly. Hence,
by (a) of Corollary~\ref{main},
\begin{eqnarray*}
\label{cheetah} \mbox{The left-hand side of (\ref{bell})} & = & \lim
_{n\to\infty
}\mathbb{E} \bigl[p_{\mu}\bigl(Z_n^{\alpha}
\bigr)\overline{p_{\nu
}\bigl(Z_n^{\alpha}\bigr)}
\bigr]
\nonumber
\\
& = & \delta_{\mu\nu} \biggl(\frac{2}{\beta} \biggr)^{l(\mu)}
z_{\mu}.
\end{eqnarray*}
By independence and rotation-invariance, we know that the right-hand
side of~(\ref{bell}) is zero if $l_j\ne l_j'$ for some $j$, or
equivalently, $\mu\ne\nu$. If $\mu=\nu$, then
%
\begin{equation}
\label{pillar} \mathbb{E} \Biggl(\prod_{j=1}^m
\xi_j^{l_j}\cdot\prod_{j=1}^m
\bar {\xi}_j^{ l_j'} \Biggr)=\prod
_{j=1}^m \mathbb{E} \bigl(|\mathbb{\xi
}_j|^{2l_j} \bigr)= \biggl(\frac{2}{\beta}
\biggr)^{l(\mu)}\prod_{j=1}^mj^{l_j}l_j!
\end{equation}
since $|\xi_j|^2\sim\frac{2j}{\beta} W$ where $W$ is the
exponential distribution with density $e^{-x}I(x\geq0)$ and $\mathbb
{E}W^l=l!$ for all integer $l\geq1$. We then obtain (\ref{bell}).
\end{pf*}

\begin{pf*}{Proof of Corollary~\ref{campus}}
Let $\alpha=\frac{2}{\beta}$ and $Z_n^{\alpha}=(e^{i\theta_1},
\ldots, e^{i\theta_n})$. Write
\[
\label{plasma} X_n=\sum_{j=1}^n
\sum_{k=0}^mc_ke^{ik\theta_j}=
\sum_{k=0}^mc_k\sum
_{j=1}^n e^{ik\theta_j}=\mu_n+\sum
_{k=1}^mc_kp_k
\bigl(Z_n^{\alpha}\bigr)
\]
with $p_k(z)=\sum_{j=1}^n z_j^{k}$ for $z=(z_1, \ldots, z_n)$. By
Theorem~\ref{Nevada} and the continuous mapping theorem, $X_n-\mu_n$
converges weakly to $Z:=\sum_{j=1}^mc_j\xi_j$ as $n\to\infty$,
where $\xi_j$'s are independent random variables and $\xi_j \sim
\mathbb{C}N(0, \frac{2j}{\beta})$ for each $j$. It is easy to check
that $Z\sim\mathbb{C}N(0, \sigma^2)$.
\end{pf*}

\begin{lemma}\label{go_to_body} Let $X\sim\mathbb{C}N(0,1)$ and $c,
d$ be two complex numbers. Then $cX + d\bar{X}=U+iV$ where $(U,
V)'\sim N_2(\bd{0}, \bdd{\Sigma})$ with
\[
\bdd{\Sigma}=\frac{1}{2}\pmatrix{ |c+\bar{d}|^2& 2
\operatorname{Im}(cd)
\cr
2 \operatorname{Im}(cd) & |c-\bar{d}|^2 } .
\]
\end{lemma}

\begin{pf}
Let $\xi$ be a standard normal random variable and $a=a_1+a_2i$ be a
complex number, where $a_1\in\mathbb{R}$ and $a_2\in\mathbb{R}$.
Then $a\xi$, as a 2-dimensional random vector, has the same
distribution as that of $(a_1\xi, a_2\xi)\sim N_2(\bd{0},\bdd{\Sigma
}_1)$ where
\begin{eqnarray*}
\bdd{\Sigma}_1= \pmatrix{ a_1^2 &
a_1a_2
\cr
a_1a_2 &
a_2^2} =\frac{1}{4}\pmatrix{ (a+
\bar{a})^2 & \bigl(\bar{a}^2-a^2\bigr)i
\cr
\bigl(\bar{a}^2-a^2\bigr)i & -(a-\bar{a})^2 }
.
\end{eqnarray*}
Let $\xi_1, \xi_2$ be i.i.d. with distribution $N(0,1)$, and $c,d$ be
complex numbers. Then
\[
c\frac{\xi_1+i\xi_2}{\sqrt{2}} +d\frac{\xi_1-i\xi_2}{\sqrt
{2}}=\frac{c+d}{\sqrt{2}}
\xi_1+i\frac{c-d}{\sqrt{2}}\xi_2,
\]
as a sum of independent (2-dimensional) normal random vectors, has
distribution $N_2(\bd{0},\bdd{\Sigma}_2)$ where
\begin{eqnarray*}
\bdd{\Sigma}_2=\pmatrix{ \sigma_{11}&
\sigma_{12}
\cr
\sigma_{12} & \sigma_{22} } .
\end{eqnarray*}
Since the covariance matrix of the sum of two independent random
variables is the sum of their individual covariance matrices, we have
\begin{eqnarray*}
 \sigma_{11}&=&\tfrac{1}{8} \bigl((c+d+\bar{c}+
\bar{d})^2-(c-d-\bar {c}+\bar{d})^2 \bigr)
\\
&=& \tfrac{1}{8}\cdot4(c+\bar{d}) (\bar{c}+d)=\tfrac{1}{2}|c+
\bar{d}|^2
\end{eqnarray*}
by using the identity $x^2-y^2=(x+y)(x-y)$. And by the identity again,
\begin{eqnarray*}
 \sigma_{22}&=&-\tfrac{1}{8} \bigl((c+d-\bar{c}-\bar{d}
)^2 - (c+\bar{c}-d-\bar{d} )^2 \bigr)
\\
&=&-\tfrac{1}{8}\cdot4(c-\bar{d}) (d-\bar{c})=\tfrac{1}{2}|c-
\bar{d}|^2.
\end{eqnarray*}
Now,
\[
\sigma_{12}=\frac{i}{8} \bigl((\bar{c}+\bar{d})^2-(c+d)^2-(
\bar {c}-\bar{d})^2+(c-d)^2 \bigr) = \frac{i}{2} (
\overline{cd}-cd ).
\]
Thus, $c\xi+d\bar{\xi} =U+iV$ where $(U, V)'\sim N_2(\bd{0}, \bdd
{\Sigma}_3)$ and
\[
\bdd{\Sigma}_3=\frac{1}{2}\pmatrix{ |c+\bar{d}|^2&
( \overline{cd}-cd)i
\cr
( \overline{cd}-cd)i & |c-\bar{d}|^2 } .
\]
\upqed\end{pf}

\begin{lemma}\label{free_fear} Let $\{X_n;  1\leq n \leq\infty\}$
be complex normal random variables with mean zero for each $n$. Then,
$X_n$ converges to $X_{\infty}$ weakly if and only if $\lim_{n\to
\infty}\mathbb{E}(X_n^p\bar{X}_n^q)=\mathbb{E}(X_{\infty}^p\bar
{X}_{\infty}^q)$ for any integers $p\geq0$ and $q\geq0$ with
$p+q\geq1$.
\end{lemma}

\begin{pf} Write $X_n=U_n+iV_n$ for all $1\leq n \leq\infty$, where
$U_n$ and $V_n$ are real random variables. Then there exists a $2\times
2$ nonnegative definite matrix $\bdd{\Sigma}_n$ such that $(U_n,
V_n)\sim N_2(\bd{0}, \bdd{\Sigma}_n)$ for each $n$. Since both $U_n$
and $V_n$ can be expressed by linear combinations of $X_n$ and $\bar
{X}_n$ and vice versa. The lemma then can be interpreted as follows:
$(U_n, V_n)$ converges to $(U_{\infty}, V_{\infty})$ weakly if and
only if $\lim_{n\to\infty}\mathbb{E}(U_n^pV_n^q)=\mathbb
{E}(U_{\infty}^pV_{\infty}^q)$ for any integers $p\geq0$ and $q\geq
0$ with $p+q\geq1$.

The sufficiency is obtained by using the moment method and the Cram\'
{e}r--Wold device. We now show the necessity. By using characteristic
functions, it is easily seen that $(U_n, V_n)$ converges to $(U_{\infty
}, V_{\infty})$ weakly if and only if $\lim_{n\to\infty} (\bdd
{\Sigma}_n)_{ij}=(\bdd{\Sigma}_{\infty})_{ij}$ for all $1\leq i,
j\leq2$. Now, assuming $(U_n, V_n)$ converges to $(U_{\infty},
V_{\infty})$, then $U_n^pV_n^q$ converges weakly to $U_{\infty
}^pV_{\infty}^q$ by the continuous mapping theorem. So we only need to
show the uniform integrability. In fact, let $r=p+q+1$, then by the H\"
{o}lder inequality, $\mathbb{E}(U_n^{2p}V_n^{2q})\leq[\mathbb
{E}(U_n^{2r})]^{p/r}\cdot[\mathbb{E}(V_n^{2r})]^{q/r}$. We know
$\mathbb{E}(U_n^{2r})=(\bdd{\Sigma}_n)_{11}^{2r}\mathbb
{E}(N(0,1)^{2r})\to(\bdd{\Sigma}_{\infty})_{11}^{2r}\times\break \mathbb
{E}(N(0,1)^{2r})=\mathbb{E}(U_{\infty}^{2r})$ as $n\to\infty$. This
shows that $\sup_{n\geq1}\mathbb{E}(U_n^{2p}V_n^{2q})<\infty$. In
particular, $\{U_n^{p}V_n^{q};  n\geq1\}$ is uniformly integrable.
\end{pf}

\begin{pf*}{Proof of Theorem~\ref{circular_case}}
Set $m=m_n=[\log
n]+1$ for $n\geq1$, where $[x]$ is the integer part of $x\geq0$.
Review (b) of Theorem~\ref{face} and (\ref{America}). We know
\begin{eqnarray*}
\mathbb{E} \bigl[p_j\bigl(Z_n^{2}\bigr)
\overline{p_k\bigl(Z_n^{2}\bigr)} \bigr]&=&0\qquad
\mbox{for any } j \ne k \geq1;
\\
 \mathbb{E} \bigl[p_j\bigl(Z_n^{2}
\bigr)p_k\bigl(Z_n^{2}\bigr) \bigr]=
\mathbb{E}p_{\mu
}\bigl(Z_n^{2}\bigr)&=&0 \qquad\mbox{for
any } j\geq1 \mbox{ and } k \geq1,
\end{eqnarray*}
where $\mu:=(j,k)$ is a partition. In particular,
\begin{eqnarray*}
& & \mathbb{E} \bigl[ \bigl(a_{j}p_j
\bigl(Z_n^{2}\bigr)+ b_{j}\overline
{p_j\bigl(Z_n^{2}\bigr)} \bigr)\cdot
\overline{ \bigl(a_{k}p_k\bigl(Z_n^{2}
\bigr)+ b_{k}\overline{p_k\bigl(Z_n^{2}
\bigr)} \bigr)} \bigr]
\\
&&\qquad =  \delta_{jk}\cdot\bigl(|a_j|^2 +
|b_j|^2\bigr)\mathbb{E}\bigl|p_j
\bigl(Z_n^{2}\bigr)\bigr|^2
\end{eqnarray*}
for all $j\geq1$ and $k\geq1$. Set
\[
Y_n:=\sum_{j=1}^{m}(a_{j}
\xi_j + b_j\bar{\xi}_j),
\]
where $\xi_j$'s are i.i.d. random variables such that $\xi_j \sim
\mathbb{C}N(0, 2j)$ for each $j\geq1$. By the Minkowski inequality,
\begin{eqnarray*}
\mathbb{E} \Biggl|\sum_{j=1}^{\infty}(a_{j}
\xi_j + b_j\bar{\xi}_j)\Biggr|^2 & \leq&
2\mathbb{E} \Biggl|\sum_{j=1}^{\infty}a_{j}
\xi_j \Biggr|^2 +2\mathbb {E}\Biggl |\sum
_{j=1}^{m} b_j\bar{
\xi}_j\Biggr|^2
\\
&\leq& 4\sum_{j=1}^{\infty}j
\bigl(|a_j|^2+|b_j|^2\bigr)<
\infty.
\end{eqnarray*}
Therefore, $Y_n$ converges weakly to $Y:=\sum_{j=1}^{\infty}(a_{j}\xi
_j + b_j\bar{\xi}_j)$. Write $a_j\xi_j + b_j \bar{\xi}_j=U_j +
iV_j$ for each $j$ such that $(U_j,V_j)\in\mathbb{R}^2$ and
$(U_j,V_j)\sim N_2(\bd{0}, \bdd{\Sigma}_j)$. Then, by Lemma~\ref{go_to_body},
\begin{eqnarray*}
\bdd{\Sigma}_j= \pmatrix{ j|a_j+\bar{b}_j|^2
& 2j\cdot\operatorname{Im}(a_jb_j)
\vspace*{2pt}\cr
2j\cdot
\operatorname{Im}(a_jb_j) & j|a_j-
\bar{b}_j|^2 }
\end{eqnarray*}
for each $j$. Thus, $\sum_{j=1}^{\infty}(a_j\xi_j + b_j \bar{\xi
}_j)$ has the law of $U+iV$ where $(U, V)'\sim N_2(\bd{0}, \bdd{\Sigma
})$ with
\begin{eqnarray*}
\bdd{\Sigma}= \pmatrix{ \displaystyle\sum_{j=1}^{\infty}j|a_j+
\bar{b}_j|^2 & \displaystyle 2\cdot\operatorname{Im}\Biggl(\sum
_{j=1}^{\infty}j a_jb_j\Biggr)
\cr
\displaystyle 2\cdot\operatorname{Im}\Biggl(\sum_{j=1}^{\infty}j
a_jb_j\Biggr) &\displaystyle \sum_{j=1}^{\infty}j|a_j-
\bar{b}_j|^2 }
\end{eqnarray*}
since the covariance matrix of the sum of independent random variables
is the sum of their individual covariance matrices. By Lemma~\ref{free_fear},
%
\begin{equation}
\label{wash_feet} \lim_{n\to\infty}\mathbb{E}\bigl[Y_n^p
\bar{Y}_n^q\bigr]=\mathbb{E}\bigl[Y^p
\bar{Y}^q\bigr].
\end{equation}
Proposition~\ref{mouth} tells us that $\mathbb{E}
[|p_j(Z_n^{2})|^2 ] \leq Kj$ for all $j\geq1$ and $n\geq2$, where
$K$ is a universal constant.
We then have
%
\begin{eqnarray}
\label{buffet} \mathbb{E} \Biggl|\sum_{j>m}^{\infty}
\bigl(a_{j}p_j\bigl(Z_n^{2}
\bigr)+ b_{j}\overline{p_j\bigl(Z_n^{2}
\bigr)} \bigr) \Biggr|^2 & = & \sum_{j>m}
\bigl(|a_{j}|^2+|b_j|^2\bigr)
\mathbb{E}\bigl|p_j\bigl(Z_n^{2}
\bigr)\bigr|^2
\nonumber
\\[-8pt]
\\[-8pt]
\nonumber
& \leq& K\sum_{j>m}^{\infty}j
\bigl(|a_{j}|^2+|b_j|^2\bigr) \to0
\end{eqnarray}
as $n\to\infty$. This shows that $\sum_{j>m}^{\infty}
(a_{j}p_j(Z_n^{2})+ b_{j}\overline{p_j(Z_n^{2})}   )$ converges
to zero in probability as $n\to\infty$. By the Slutsky lemma, to
prove the theorem, we only need to show
%
\begin{equation}
\label{pick} X_n:=\sum_{j=1}^m
\bigl(a_{j}p_j\bigl(Z_n^{2}
\bigr)+ b_{j}\overline {p_j\bigl(Z_n^{2}
\bigr)} \bigr) \to Y
\end{equation}
weakly as $n \to\infty$. Thus from (\ref{wash_feet}),
similar to (\ref{sales}), to prove (\ref{pick}) it suffices to show that
%
\begin{equation}
\label{deer} \lim_{n\to\infty} \bigl(\mathbb{E}
\bigl[X_n^p\bar{X}_n^q\bigr]-
\mathbb {E}\bigl[Y_n^p\bar{Y}_n^q
\bigr] \bigr)=0.
\end{equation}
Recall the multinomial formula,
%
\begin{equation}
\label{university} (x_1+\cdots+x_k)^p=\sum
_{l_1+\cdots+l_k=p}\pmatrix{p\cr l_1,l_2,
\ldots , l_k}x_1^{l_1}x_2^{l_2}
\cdots x_k^{l_k}
\end{equation}
for any complex number $x_i$'s, positive integers $k\geq2$ and $p\geq
1$, where $l_i$'s are nonnegative integers. Note that $X_n$ is a sum of
$2m$ terms. Expand $X_n^p\bar{X}_n^q$ to have
\begin{eqnarray*}
 \mathbb{E}\bigl[X_n^p\bar{X}_n^q
\bigr]
&=&\sum\pmatrix{p\cr l_1, \ldots, l_{2m}}\cdot\pmatrix{q
\cr l_1',\ldots, l_{2m}'}
\cdot\prod_{j=1}^{m} \bigl(a_{j}^{l_j}b_j^{l_{m+j}}
\bigr)\cdot\prod_{j=1}^{m}
\bar{a}_{j}^{ l_j'}\bar{b}_{j}^{ l_{j+m}'}
\\
& &\quad{}\times \mathbb{E} \Biggl[\prod_{j=1}^{m}
\bigl(p_j\bigl(Z_n^{2}\bigr)^{l_j}
\overline {p_j\bigl(Z_n^{2}
\bigr)}^{l_{m+j}} \bigr)\cdot\prod_{j=1}^m
\bigl(\overline {p_j\bigl(Z_n^{2}
\bigr)}^{ l_j'}p_j\bigl(Z_n^{2}
\bigr)^{ l_{m+j}'} \bigr) \Biggr],
\end{eqnarray*}
where the sum runs over all possible nonnegative integers $l_j$'s and
$l_j'$'s with $\sum_{j=1}^{2m}l_j=p$ and $\sum_{j=1}^{2m}l_j'=q$.
Rearranging the products in the expectation, we get
%
\begin{eqnarray}\label{long1}
\mathbb{E}\bigl[X_n^p\bar{X}_n^q
\bigr] &=&\sum\pmatrix{p\cr l_1, \ldots, l_{2m}}\cdot\pmatrix{q
\cr l_1',\ldots, l_{2m}'}
\cdot\prod_{j=1}^{m} \bigl(a_{j}^{l_j}b_j^{l_{m+j}}
\bigr)\cdot\prod_{j=1}^{m}
\bar{a}_{j}^{ l_j'}\bar{b}_{j}^{ l_{m+j}'}
\nonumber
\\[-8pt]
\\[-8pt]
\nonumber
&&\quad{}\times \mathbb{E} \Biggl[\prod_{j=1}^{m}p_j
\bigl(Z_n^{2}\bigr)^{l_j+l_{m+j}'}\cdot \prod
_{j=1}^m\overline{p_j
\bigl(Z_n^{2}\bigr)}^{ l_j'+l_{m+j}} \Biggr].
\end{eqnarray}
Similarly,
%
\begin{eqnarray} \label{long2}
\mathbb{E}\bigl[Y_n^p\bar{Y}_n^q
\bigr]
&=&\sum\pmatrix{p\cr l_1, \ldots, l_{2m}}\cdot\pmatrix{q
\cr l_1',\ldots, l_{2m}'}
\cdot\prod_{j=1}^{m} \bigl(a_{j}^{l_j}b_j^{l_{m+j}}
\bigr)\cdot\prod_{j=1}^{m}
\bar{a}_{j}^{ l_j'}\bar{b}_{j}^{ l_{m+j}'}
\nonumber
\\[-8pt]
\\[-8pt]
\nonumber
& &\quad{}\times\mathbb{E} \Biggl[\prod_{j=1}^{m}
\xi_j^{l_j+l_{m+j}'}\cdot\prod_{j=1}^m
\bar{\xi}_j^{ l_j'+l_{m+j}} \Biggr].
\end{eqnarray}
We claim that
%
\begin{eqnarray}\label{growl}
&&\Biggl |\mathbb{E} \Biggl[\prod_{j=1}^{m}p_j
\bigl(Z_n^{2}\bigr)^{l_j+l_{m+j}'}\cdot\prod
_{j=1}^m\overline {p_j
\bigl(Z_n^{2}\bigr)}^{ l_j'+l_{m+j}} \Biggr]
\nonumber\\
&&\hspace*{32pt}\quad{}- \mathbb{E} \Biggl[\prod_{j=1}^{m}\xi
_j^{l_j+l_{m+j}'}\cdot\prod_{j=1}^m
\bar{\xi}_j^{ l_j'+l_{m+j}} \Biggr]\Biggr |
\\
&&\qquad\leq C_{p,q}\cdot\frac{m}{n}\cdot \prod
_{j=1}^mj^{(l_j+l_{m+j})/2}\cdot\prod
_{j=1}^mj^{(l_j'+l_{m+j}')/2}
\nonumber
\end{eqnarray}
uniformly for all possible $l_j$'s and $l_j'$'s in the two sums, where
$C_{p,q}$ is constant depending on $p$ and $q$ only. In fact, let $\mu
$ and $\nu$ be two partitions so that
\begin{eqnarray*}
\mu&=&\bigl(1^{l_1+l_{m+1}'}, 2^{l_2+l_{m+2}'}, \ldots, m^{l_m+l_{2m}'}
\bigr);
\\
 \nu&=&\bigl(1^{l_1'+l_{m+1}}, 2^{l_2'+l_{m+2}}, \ldots, m^{l_m'+l_{2m}}
\bigr).
\end{eqnarray*}
Then $l(\mu)=\sum_{j=1}^m(l_j+l_{m+j}')\leq p+q$ and similarly $l(\nu
)\leq p+q$, and
\[
\label{Minnesota} K:=|\mu| \vee|\nu| =\sum_{j=1}^mj
\bigl(l_j+l_{m+j}'\bigr) \vee\sum
_{j=1}^mj\bigl(l_j'+l_{m+j}
\bigr) \leq m(p+q).
\]
According to this notation,
%
\begin{equation}
\label{grey} \mathbb{E} \Biggl[\prod_{j=1}^{m}p_j
\bigl(Z_n^{2}\bigr)^{l_j+l_{m+j}'}\cdot\prod
_{j=1}^m\overline{p_j
\bigl(Z_n^{2}\bigr)}^{ l_j'+l_{m+j}} \Biggr]=\mathbb{E}
\bigl[p_{\mu}\bigl(Z_n^{2}\bigr)
\overline{p_{\nu}\bigl(Z_n^{2}\bigr)} \bigr].
\end{equation}
By (\ref{pillar}),
%
\begin{eqnarray}\label{brown}
& & \mathbb{E} \Biggl[\prod_{j=1}^{m}
\xi_j^{l_j+l_{m+j}'}\cdot\prod_{j=1}^m
\bar{\xi}_j^{ l_j'+l_{m+j}} \Biggr]
\nonumber
\\[-8pt]
\\[-8pt]
\nonumber
&&\qquad=\delta_{\mu\nu} \biggl(\frac{2}{\beta} \biggr)^{l(\mu)}\prod
_{j=1}^mj^{l_j+l_{m+j}'}(l_j+l_{m+j'})!=
\delta_{\mu\nu}\alpha ^{l(\mu)}z_{\mu},
\end{eqnarray}
where $\alpha=\frac{2}{\beta}=2$ and $z_{\mu}$ is as in (\ref
{insect}). Since $\sum_{j=1}^m(l_j+l_{m+j'}) \leq p+q$, then
\[
\operatorname{Card} \{1\leq j \leq m; l_j+l_{m+j'} \geq2 \} \leq
\frac{p+q}{2}.
\]
Using $l_j+l_{m+j'} \leq p+q$ for all $1\leq j \leq m$, we get
\[
0<\alpha^{l(\mu)}z_{\mu}\leq \bigl(2^{p+q}\bigl((p+q)!
\bigr)^{(p+q)/2} \bigr)\cdot\prod_{j=1}^mj^{l_j+l_{m+j}'}.
\]
A similar inequality also holds for $\alpha^{l(\nu)}z_{\nu}$. From
(a) and (b) of Corollary~\ref{theory}, we see that
\[
\bigl|\mathbb{E} \bigl[p_{\mu}\bigl(Z_n^{2}\bigr)
\overline{p_{\nu
}\bigl(Z_n^{2}\bigr)} \bigr]-
\delta_{\mu\nu}\alpha^{l(\mu)}z_{\mu}\bigr | \leq
C_{p,q}\cdot\frac{m}{n}\cdot \prod
_{j=1}^mj^{(l_j+l_{m+j}')/2}\cdot\prod
_{j=1}^mj^{(l_j'+l_{m+j})/2},
\]
where $C_{p,q}$ is a constant depending on $p$ and $q$ only. This
together with (\ref{grey}) and (\ref{brown}) yields (\ref{growl}).

Now, combining (\ref{long1}), (\ref{long2}) and (\ref{growl}), we
arrive at
\begin{eqnarray*}
& & \bigl|\mathbb{E}\bigl[X_n^p\bar{X}_n^q
\bigr]-\mathbb{E}\bigl[Y_n^p\bar {Y}_n^q
\bigr]\bigr |
\\
&&\qquad \leq C_{p,q}\cdot\frac{m}{n}\cdot\sum\pmatrix{p\cr
l_1, \ldots, l_{2m}}\cdot\pmatrix{q\cr l_1',
\ldots, l_{2m}'} \cdot\prod_{j=1}^{m}
\bigl(\sqrt{j}|a_{j}|\bigr)^{l_j}\bigl(\sqrt{j}|b_j|\bigr)^{l_{m+j}}
\\
& &\qquad\quad {}\times\prod_{j=1}^{m} \bigl(\sqrt
{j}|a_{j}|\bigr)^{l_j'}\bigl(\sqrt{j}|b_j|\bigr)^{l_{m+j}'}
\\
&&\qquad =  C_{p,q}\cdot\frac{m}{n}\cdot \Biggl(\sum
_{j=1}^m\bigl(\sqrt{j} |a_{j}|+
\sqrt{j}|b_j|\bigr) \Biggr)^{p}\cdot \Biggl(\sum
_{j=1}^m\bigl(\sqrt{j} |a_{j}|+
\sqrt{j}|b_j|\bigr) \Biggr)^{q}
\\
&&\qquad =  C_{p,q}\cdot\frac{m}{n}\cdot \Biggl(\sum
_{j=1}^m\bigl(\sqrt{j} |a_{j}|+
\sqrt{j}|b_j|\bigr) \Biggr)^{p+q},
\end{eqnarray*}
where (\ref{university}) is used in the first identity. From the
inequality $(x_1+\cdots+ x_{2m})^2\leq2m (x_1^2+\cdots+ x_{2m}^2)$
for any real number $x_i$'s we see that
\[
\bigl|\mathbb{E}\bigl[X_n^p\bar{X}_n^q
\bigr]-\mathbb{E}\bigl[Y_n^p\bar{Y}_n^q
\bigr] \bigr| \leq \bigl(C_{p,q}'\sigma^{p+q}\bigr)
\frac{1}{n}m^{1+(p+q)/2}\to0
\]
as $n\to\infty$ since $m=[\log n]$ for $n\geq3$, where $C_{p,q}'$ is
a constant depending on $p$ and $q$ only. This confirms (\ref{deer}).
\end{pf*}

\begin{pf*}{Proof of Theorem~\ref{symplectic_case}}
From the assumption
that $\sum_{j=1}^{\infty}(j\log j)(|a_j|^2+|b_j|^2)\in(0, \infty)$,
we know $\sigma^2\in(0, \infty)$. Take $m=m_n=[\log n]$ for $n\geq
3$. By (ii) of Proposition~\ref{search}, the assumption $\sum_{j=1}^{\infty}(j\log j)(|a_j|^2+|b_j|^2)\in(0, \infty)$ and the
same argument as the derivation of (\ref{buffet}), to prove the
theorem, it is enough to show that
\[
\sum_{j=1}^m \bigl(a_{j}p_j
\bigl(Z_n^{1/2}\bigr)+ b_{j}\overline
{p_j\bigl(Z_n^{1/2}\bigr)} \bigr) \to\sum
_{j=1}^{\infty}(a_j\xi_j
+ b_j \bar {\xi}_j)
\]
weakly as $n\to\infty$, where $\{\xi_j;  j\geq1\}$ are independent
random variables with $\xi_j \sim\mathbb{C}N(0, \frac{1}{2}j)$ for
each $j$. Write $a_j\xi_j + b_j \bar{\xi}_j=U_j+ i V_j$ for each $j$
with $(U_j, V_j)\in\mathbb{R}^2$. By Lemma~\ref{go_to_body}, $(U_j,
V_j)$ has the distribution $N_2(\bd{0}, \bdd{\Sigma}_j)$ where
\[
\bdd{\Sigma}_j=\frac{1}{4} \pmatrix{ j|a_j+
\bar{b}_j|^2 & 2j\cdot\operatorname{Im}(a_jb_j)\vspace*{2pt}
\cr
2j\cdot\operatorname{Im}(a_jb_j) & j|a_j-
\bar{b}_j|^2 } .
\]
It follows from the independence that $\sum_{j=1}^{\infty}(a_j\xi_j
+ b_j \bar{\xi}_j)$ has the law of $U+iV$ where $(U, V)'\sim N_2(\bd
{0}, \bdd{\Sigma})$ with
\[
\bdd{\Sigma}=\frac{1}{4} \pmatrix{ \displaystyle\sum_{j=1}^{\infty}j|a_j+
\bar{b}_j|^2 & \displaystyle 2\cdot\operatorname{Im}\Biggl(\sum
_{j=1}^{\infty}ja_jb_j\Biggr)
\cr
\displaystyle 2\cdot\operatorname{Im}\Biggl(\sum_{j=1}^{\infty}ja_jb_j
\Biggr) & \displaystyle \sum_{j=1}^{\infty
}j|a_j-
\bar{b}_j|^2 } .
\]
Then the rest proof will be completed by following the same arguments
as in the corresponding parts in the proof of Theorem~\ref{circular_case}.
\end{pf*}

\begin{appendix}
\section*{Appendix}\label{noodle_eat}

In this section, we calculate some moments for the circular $\beta
$-ensembles. The first result below is an independent check of the
second moment of the trace of a COE given in (\ref{glue}). The
derivation does not depend on the Jack function as used in Section~\ref{proofs}. It only uses the distribution of the entries of the COE.
%
\setcounter{lemma}{0}
\begin{lemma}\label{blue} Let $W_n$ be an $n\times n$ circular
orthogonal ensemble (COE), that is, $W_n=U_n^T U_n$ for some
Haar-invariant unitary matrix $U_n$.
Then $\mathbb{E} [|\Tr(W_n)|^2  ]=2n/(n+1)$ for all $n\geq2$.
\end{lemma}

\begin{pf*}{First Proof of Lemma~\ref{blue}}
We prove the lemma in
three steps.

\textit{Step} 1. Write $U_n=(u_{rs})$. First, we claim that
%
\setcounter{equation}{0}
\begin{equation}
\label{snow} \mathbb{E}\bigl[u^2_{rs}
\bar{u}_{pq}^2\bigr]=0
\end{equation}
if $r\ne p$ or $s\ne q$. In fact, since $U_n$ is Haar-invariant
unitary, the distributions of $UU_n$ and $U_nU$ are the same as that of
$U_n$ for any unitary matrix $U$. In particular, take $U=\operatorname
{diag}(e^{i\theta_k})_{1\leq k \leq n}$ to obtain that
%
\begin{eqnarray}
\label{grass3} \mathcal{L} \bigl( \bigl(e^{i\theta_r}u_{rs}
\bigr)_{1\leq r, s\leq
n} \bigr)&=&\mathcal{L} \bigl( \bigl(e^{i\theta_s}u_{rs}
\bigr)_{1\leq r, s\leq n} \bigr)
\nonumber
\\[-8pt]
\\[-8pt]
\nonumber
&=&\mathcal{L} \bigl( (u_{rs}
)_{1\leq r, s\leq n} \bigr)
\end{eqnarray}
for any $\theta_1, \ldots, \theta_n \in\mathbb{R}$, where
$\mathcal{L}(X)$ is the joint distribution of the entries of random
matrix $X$.
If $r\ne p$, taking $\theta_r-\theta_p=\pi/2$, then by (\ref
{grass3}), we have that
\[
\mathbb{E}\bigl[u^2_{rs}\bar{u}_{pq}^2
\bigr]=e^{2i(\theta_r-\theta_p)}\mathbb {E}\bigl[u^2_{rs}
\bar{u}_{pq}^2\bigr]=-\mathbb{E}\bigl[u^2_{rs}
\bar{u}_{pq}^2\bigr],
\]
which means (\ref{snow}). The case for $s=q$ can be proved similarly.

\textit{Step} 2. Recall notation
$(2m-1)! !=(2m-1)(2m-3)\cdots3\cdot1$ for any integer $m\geq1$,
and $(-1)!!=1$ by convention. We have the following fact (Lemma~2.4
from \cite{Jiang09b}):
%
\begin{equation}
\label{sushi} \mathbb{E}\bigl[\xi_1^{a_1}
\xi_2^{a_2}\cdots \xi_{n}^{a_{n}}\bigr]=
\frac{\prod_{i=1}^{n}(2a_i-1)! !}{\prod_{i=1}^{a}(n+2i-2)},
\end{equation}
where $a_1,\ldots, a_{n}$ are nonnegative integers with $a=\sum_{i=1}^{n}a_i$, $\xi_i=X_i^2/(X_1^2+\cdots+ X_n^2)$ and $X_1, \ldots
, X_n$ are i.i.d.
random variables with $X_1\sim N(0,1)$.

\textit{Step} 3. Evidently, $\Tr(W_n)=\sum_{1\leq i, j\leq n}u_{ij}^2$.
Notice, from the invariant property, by exchanging some rows and some
columns of $U_n$, we see that the distributions of $u_{rs}$ and
$u_{11}$ are identical for any $1\leq r, s\leq n$. By (\ref{snow}),
%
\begin{eqnarray}
\label{morning} \mathbb{E} \bigl[\bigl|\Tr(W_n)\bigr|^2 \bigr] &=&
\mathbb{E} \biggl[ \biggl(\sum_{r,s}u_{rs}^2
\biggr) \biggl(\sum_{p,q}\bar{u}_{p,q}^2
\biggr) \biggr]
\nonumber
\\[-8pt]
\\[-8pt]
\nonumber
& =&\mathbb{E} \biggl[\sum_{r,s}|u_{rs}|^4
\biggr]=n^2E \bigl[|u_{11}|^4 \bigr].
\nonumber
\end{eqnarray}
It is known (e.g., Lemma~2.1 in \cite{Jiang09a,Jiang09b}) that the
probability distribution of $|u_{11}|^2$ is the same as that of $(X_1^2
+ X_2^2)/\sum_{i=1}^{2n}{X_i}^2$.
By (\ref{sushi}),
\[
\mathbb{E}\bigl[\xi_1^2\bigr]=\frac{3}{2n(2n+2)}
\quad\mbox{and}\quad \mathbb {E}[\xi_1\xi_2]=\frac{1}{2n(2n+2)}.
\]
Then
\[
\mathbb{E}\bigl[|u_{11}|^4\bigr]=\mathbb{E}\bigl[(
\xi_1 + \xi_2)^2\bigr]=2\mathbb {E}\bigl[
\xi_1^2\bigr] + 2\mathbb{E}[\xi_1
\xi_2]=\frac{2}{n(n+1)}.
\]
Substitute this into (\ref{morning}) to see that $\mathbb{E}
[|\Tr(W_n)|^2 ]=2n/(n+1)$.\vadjust{\goodbreak}
\end{pf*}

\begin{pf*}{Second Proof of Lemma~\ref{blue}}
We use the following formula due to Collins \cite{Collins} (see also
\cite{Matsumoto}):
let $(u_{ij})_{1 \le i,j \le n}$ be an $n \times n$ CUE matrix
(or equivalently, an Haar-distributed unitary matrix) and
let $i_1,\ldots,i_k$, $j_1,\ldots,j_k$, $i_{1}',\ldots,i_k'$,
$j_1',\ldots,j_k'$
be elements in $\{1,2,\ldots,n\}$. Then
%
\begin{eqnarray}
\label{WgFormula} &&\mathbb{E} [ u_{i_1 j_1} \cdots u_{i_k j_k}
\overline{u_{i_1'j_1'} \cdots u_{i_k'j_k'}} ]
\nonumber
\\[-8pt]
\\[-8pt]
\nonumber
&&\qquad= \sum_{\sigma, \tau\in\mathfrak{S}_k} \operatorname{Wg}_{n,k}\bigl(
\sigma ^{-1} \tau\bigr) \Biggl( \prod_{p=1}^k
\delta_{i_p, i_{\sigma(p)}'} \Biggr) \Biggl( \prod_{q=1}^k
\delta_{j_q, j_{\tau(q)}'} \Biggr).
\end{eqnarray}
Here, $\mathfrak{S}_k$ is the symmetric group and $\operatorname{Wg}_{n,k}$
is a
class function on $\mathfrak{S}_k$,
called the Weingarten function for the unitary group.
For our purpose, we do not need the explicit definition of $\operatorname
{Wg}_{n,k}$ but
use the case for $k=2$. In fact, for $n \ge2$, we know (see (5.2) of
\cite{Collins})
%
\begin{eqnarray}\label{Wg}
\operatorname{Wg}_{n,2}(\mathrm{id}_2)& = &\frac{1}{n^2-1}\quad
\mbox {and}
\nonumber
\\[-8pt]
\\[-8pt]
\nonumber
\operatorname{Wg}_{n,2}\bigl((1 \ 2)\bigr) &=& -\frac{1}{n(n^2-1)},
\end{eqnarray}
where $\mathrm{id}_2$ and $(1 \ 2)$ are the identity permutation and
the transposition on $\{1,2\}$, respectively.

We have $|\Tr(W_n)|^2= \sum_{r,s,p,q} u_{rs}^2 \overline{u}_{pq}^2$.
By \eqref{WgFormula}, $\mathbb{E} [ u_{rs}^2 \overline{u}_{pq}^2]$
is zero unless
$r=p$ and $s=q$. Moreover, $\mathbb{E} [ u_{rs}^2 \overline
{u}_{rs}^2] =
\mathbb{E}[|u_{11}|^4]$ for all $1\leq r,  s\leq n$.
Therefore, using~\eqref{WgFormula} and \eqref{Wg}, we obtain
\[
\mathbb{E} \bigl[\bigl |\Tr(W_n)\bigr|^2 \bigr] = n^2
\mathbb{E} \bigl[|u_{11}|^4\bigr] =2n^2 \bigl\{
\operatorname{Wg}_{n,2}( \mathrm{id}_2)+ \operatorname{Wg}_{n,2}
\bigl( (1\ 2)\bigr) \bigr\}= \frac{2n}{n+1}.
\]
\upqed\end{pf*}

Lemma~\ref{blue} corresponds to the conclusion for $\beta=1$ in (\ref
{cry}), which is derived through Proposition~\ref{exact} by the Jack functions.
Now we apply the same proposition to derive some other moments for the
circular $\beta$-ensembles. Let $p_k$ and $Z_n$ be as in Theorem~\ref
{face}.

\begin{example*}
Assume $\alpha=2/\beta>0$. For $n\geq2$,
\begin{eqnarray}\label{right}
\mathbb{E} \bigl[\bigl|p_1(Z_n)\bigr|^4\bigr] &= &
\frac{2n\alpha^2(n^2+ 2(\alpha-1)n-\alpha)}{(n+\alpha
-1)(n+\alpha-2)(n+ 2\alpha-1)}
\nonumber
\\[-8pt]
\\[-8pt]
\nonumber
& = & \cases{ \displaystyle\frac{8(n^2+2n-2)}{(n+1)(n+3)},& \quad$\mbox{if $\beta=1$;}$
\nonumber
\vspace*{2pt}
\cr
2, &\quad \mbox{if $\beta=2$};
\nonumber
\vspace*{2pt}
\cr
\displaystyle\frac{2n^2-2n-1}{(2n-1)(2n-3)}, &\quad $\mbox{if $\beta=4$.}$}
\end{eqnarray}
\end{example*}

\begin{example*} Assume $\alpha=2/\beta>0$. For $n\geq2$,
%
\begin{eqnarray}\label{left}
\mathbb{E} \bigl[\bigl|p_2(Z_n)\bigr|^2\bigr]&=&
\frac{2\alpha n(n^2 + 2(\alpha-1)n +
\alpha^2-3\alpha+1)}{(n+\alpha-1)(n+2\alpha-1)(n+\alpha-2)}
\nonumber
\\[-8pt]
\\[-8pt]
\nonumber
&= & \cases{ \displaystyle\frac{4(n^2+2n-1)}{(n+1)(n+3)},& \quad $\mbox{if $\beta=1$;}$
\nonumber
\vspace*{2pt}
\cr
2,&\quad $ \mbox{if $\beta=2$;}$\vspace*{2pt}
\cr
\displaystyle\frac{4n^2-4n-1}{(2n-1)(2n-3)},&\quad $\mbox{if $
\beta=4$.}$}
\end{eqnarray}
\end{example*}

\begin{example*}
Assume $\alpha=2/\beta>0$. For $n\geq2$,
%
\begin{eqnarray}\label{upper}
 \mathbb{E} \bigl[p_2(Z_n)\overline{p_1(Z_n)^2}
\bigr] &= &\mathbb{E} \bigl[ \overline{p_2(Z_n)}p_1(Z_n)^2
\bigr]
\nonumber
\\
&= & \frac{2\alpha^2(\alpha-1)n}{(n+\alpha-1)(n+2\alpha-1)(n+\alpha
-2)}
\\
& = & \cases{ \displaystyle\frac{8}{(n+1)(n+3)},&\quad $\mbox{if $\beta=1$;}$\vspace*{2pt}
\cr
0,&\quad $
\mbox{if $\beta=2$;}$\vspace*{2pt}
\cr
\displaystyle\frac{-1}{(2n-1)(2n-3)},& \quad$\mbox{if $
\beta=4$.}$}\nonumber
\end{eqnarray}
In particular, if $\beta\ne2$, as $n\to\infty$,
%
\begin{equation}
\label{Chinese} \mathbb{E} \bigl[p_2(Z_n)
\overline{p_1(Z_n)^2} \bigr] \sim2\alpha
^2(\alpha-1)n^{-2}.
\end{equation}

\begin{pf*}{Proofs of (\ref{cry}), (\ref{right}), (\ref{left}) and
(\ref{upper})}
Let $n\geq2$, $\mu$ and $\nu$ be partitions of $2$.
Set $\alpha=2/\beta$. By Proposition~\ref{exact}
and (\ref{piano}), we have
\begin{eqnarray}\label
{long}
&&\mathbb{E} \bigl[p_{\mu}(Z_n)\overline{p_{\nu}(Z_n)}
\bigr]\nonumber\\
&&\qquad=\alpha ^{l(\mu)+ l(\nu)} \biggl(\frac{4\alpha^{2-l(\mu)}\alpha^{2-l(\nu
)}}{2\alpha^2(\alpha+1)}\frac{n(n+\alpha)}{(n+\alpha-1)(n+2\alpha
-1)}
\nonumber
\\[-8pt]
\\[-8pt]
\nonumber
&&\qquad\quad\hspace*{48pt}{}+ \frac{4(-1)^{2-l(\mu)}(-1)^{2-l(\nu
)}}{2\alpha(\alpha+1)}\frac{n(n-1)}{(n+\alpha-1)(n+\alpha-2)} \biggr)
\\
&&\qquad=\frac{2\alpha^{l(\mu)+ l(\nu)}n}{\alpha^2(\alpha+1)(n+\alpha
-1)} \biggl(\frac{\alpha^{4-l(\mu)-l(\nu)}(n+\alpha)}{n+2\alpha-1} +\frac{(-1)^{l(\mu) + l(\nu)}\alpha(n-1)}{n+\alpha-2}
\biggr).\nonumber
\nonumber
\end{eqnarray}
\begin{longlist}[(iii)]
\item[(i)] Take $\mu=\nu=(1)$ in Proposition~\ref{exact}. Since
$\theta_{(1)}^{(1)}(\alpha)=1,  C_{(1)}(\alpha)=\alpha$ for any
$\alpha>0$, we obtain (\ref{cry}).

\item[(ii)] Taking $\mu=\nu=(1,1)$ in (\ref{long}), (\ref{right}) follows.

\item[(iii)] Taking $\mu=\nu=(2)$ in (\ref{long}), (\ref{left}) follows.

\item[(iv)] Taking $\mu=(2)$ and $\nu=(1,1)$ in (\ref{long}), we get the
identity for the first expectation in (\ref{upper}). Since the value
of the expectation is real, the identity for the second expectation
follows. With the earlier conclusion, (\ref{Chinese}) is obvious.\quad\qed
\end{longlist}
\noqed\end{pf*}
\end{example*}
\end{appendix}

\section*{Acknowledgements}
The first author thanks Professor Persi Diaconis very much for
introducing to him the moment problem studied in this paper.
We thank Drs Beno\^{i}t Collins, Ming Gao, Yongcheng Qi, Ke Wang,
Gongjun Xu and Lin Zhang for very helpful communications and checks of
our proofs. We thank the referees' suggestions to study the central
limit theorems in Section~\ref{CLT_introduction}.






\printaddresses
\end{document}